\documentclass[12pt]{amsart} 

\usepackage{amsmath}
\usepackage{amsthm}
\usepackage{amssymb}
\usepackage{amsfonts}
\usepackage{amscd}

\makeatletter

\begin{document}

\bibliographystyle{amsalpha}
\makeatletter

\newcommand{\e}{\epsilon}
\newcommand{\0}{{\bold 0}}
\newcommand{\w}{{\bold w}}
\newcommand{\y}{{\bold y}}
\newcommand{\p}{{\bold p}}
\newcommand{\q}{{\bold q}}
\newcommand{\ba}{{\bold a}}
\newcommand{\bc}{{\bold c}}
\newcommand{\be}{{\bold e}}

\newcommand{\balpha}{{\boldsymbol \alpha}}
\newcommand{\bbeta}{{\boldsymbol \beta}}
\newcommand{\bmu}{{\boldsymbol \mu}}
\newcommand{\blambda}{{\boldsymbol \lambda}}
\newcommand{\bt}{{\bf t}}
\newcommand{\bl}{{\bf l}}
\newcommand{\z}{{\bold z}}
\newcommand{\x}{{\bold x}}
\newcommand{\bH}{{\mathbf H}}
\newcommand{\N}{{\bold N}}
\newcommand{\Z}{{\bold Z}}
\newcommand{\F}{{\bold F}}
\newcommand{\R}{{\bold R}}
\newcommand{\Q}{{\bold Q}}
\newcommand{\C}{{\bold C}}
\newcommand{\BP}{{\mathbf P}}
\newcommand{\BF}{{\mathbf F}}
\newcommand{\cA}{{\mathcal A}}
\newcommand{\cK}{{\mathcal K}}
\newcommand{\cG}{{\mathcal G}}
\newcommand{\cO}{{\mathcal O}}
\newcommand{\cX}{{\mathcal X}}
\newcommand{\cV}{{\mathcal V}}
\newcommand{\cH}{{\mathcal H}}
\newcommand{\cM}{{\mathcal M}}
\newcommand{\cD}{{\mathcal D}}
\newcommand{\cB}{{\mathcal B}}
\newcommand{\cC}{{\mathcal C}}
\newcommand{\cT}{{\mathcal T}}
\newcommand{\cI}{{\mathcal I}}
\newcommand{\cJ}{{\mathcal J}}
\newcommand{\cS}{{\mathcal S}}
\newcommand{\cF}{{\mathcal F}}
\newcommand{\cE}{{\mathcal E}}
\newcommand{\cP}{{\mathcal P}}
\newcommand{\cR}{{\mathcal R}}
\newcommand{\cRP}{{\mathcal R \mathcal P}}
\newcommand{\cIF}{{\mathcal I \mathcal F}}
\newcommand{\cPF}{{\mathcal P \mathcal F}}

\newcommand{\E}{{\mathbf E}}
\newcommand{\V}{{\mathbf V}}

\newcommand{\cL}{{\mathcal L}}
\newcommand{\cY}{{\mathcal Y}}
\newcommand{\sB}{{\sf B}}
\newcommand{\sE}{{\sf E}}

\newcommand{\sA}{{\sf A}}
\newcommand{\ga}{{\sf a}}
\newcommand{\es}{{\sf s}}
\newcommand{\m}{{\bold m}}
\newcommand{\bS}{{\bold S}}

\newcommand{\ovf}{{\overline{f}}}

\newcommand{\ihra}{\stackrel{i}{\hookrightarrow}}
\newcommand\rank{\mathop{\rm rank}\nolimits}
\newcommand\im{\mathop{\rm Im}\nolimits}
\newcommand\coker{\mathop{\rm coker}\nolimits}
\newcommand{\Aut}{\mathop{\rm Aut}\nolimits}
\newcommand{\Tr}{\mathop{\rm Tr}\nolimits}
\newcommand\Li{\mathop{\rm Li}\nolimits}
\newcommand\NS{\mathop{\rm NS}\nolimits}
\newcommand\Hom{\mathop{\rm Hom}\nolimits}
\newcommand\Isomd{\mathop{\rm Isomd}\nolimits}
\newcommand\Ext{\mathop{\rm Ext}\nolimits}
\newcommand\End{\mathop{\rm End}\nolimits}
\newcommand\Pic{\mathop{\rm Pic}\nolimits}
\newcommand\Spec{\mathop{\rm Spec}\nolimits}
\newcommand\Hilb{\mathop{\rm Hilb}\nolimits}
\newcommand\pardeg{\mathop{\rm pardeg}\nolimits}
\newcommand\Ker{\mathop{\rm Ker}\nolimits}
\newcommand\RH{\mathop{\bf RH}\nolimits}
\newcommand{\len}{\mathop{\rm len}\nolimits}
\newcommand{\res}{\mathop{\sf res}\nolimits}
\newcommand\Quot{\mathop{\rm Quot}\nolimits}
\newcommand\Grass{\mathop{\rm Grass}\nolimits}
\newcommand\Proj{\mathop{\rm Proj}\nolimits}
\newcommand{\codim}{\mathop{\rm codim}\nolimits}

\newcommand\lra{\longrightarrow}
\newcommand\ra{\rightarrow}
\newcommand\la{\leftarrow}
\newcommand\JG{J_{\Gamma}}
\newcommand{\wvskp}{\vspace{1cm}}
\newcommand{\vskp}{\vspace{5mm}}
\newcommand{\nvskp}{\vspace{1mm}}
\newcommand{\nid}{\noindent}
\newcommand{\new}{\nvskp \nid}
\newtheorem{Assumption}{Assumption}[section]
\newtheorem{Theorem}{Theorem}[section]
\newtheorem{Lemma}{Lemma}[section]
\newtheorem{Remark}{Remark}[section]
\newtheorem{Corollary}{Corollary}[section]
\newtheorem{Conjecture}{Conjecture}[section]
\newtheorem{Proposition}{Proposition}[section]
\newtheorem{Example}{Example}[section]
\newtheorem{Definition}{Definition}[section]
\newtheorem{Question}{Question}[section]
\newtheorem{Fact}{Fact}[section]

\title[Moduli of stable parabolic connections]
{Moduli of stable parabolic connections, 
Riemann-Hilbert correspondence and geometry of 
Painlev\'{e} equation of type VI, part II}

%\keywords{Stable parabolic connections , Representation of fundamental 
%groups, Riemann-Hilbert correspondences, Symplectic structure, Painlev\'e equat%ions, Garnier equations}

\author[M. Inaba, K. Iwasaki]{Michi-aki Inaba$^1$, Katsunori Iwasaki$^2$}

\thanks{$^1$Partly supported by Grant-in Aid for Scientific Research (Wakate-B-15740018).}
\thanks{$^2$Partly supported by Grant-in Aid for Scientific Research  (B-12440043).}
%\author[K. Iwasaki]{Katsunori Iwasaki$^{2}$
%\footnote{$^{2}$Partly supported by Grant-in Aid for Scientific Research  (B-12%440043).}}

\author[M.-H. Saito]{Masa-Hiko Saito$^3$}
\thanks{$^3$Partly supported by Grant-in Aid for Scientific Research  (B-16340009), (Houga-16654004), and JSPS-NWO exchange program.}

\dedicatory{Dedicated to Professor Masaki Maruyama \\ on his 60th birthday}

\address{{\rm Michi-aki Inaba}\\
Department of Mathematics, Faculty of Science, Kyoto University, \\
Kyoto, 606-8502, Japan}
\email{inaba@math.kyoto-u.ac.jp}
\address{{\rm Katsunori Iwasaki}\\
Faculty of Mathematics, Kyushu University, \\
 6-10-1, Hakozaki, Higashi-ku, Fukuoka 812-8581, Japan}
\email{iwasaki@math.kyushu-u.ac.jp}
\address{{\rm Masa-Hiko Saito} \\
 Department of Mathematics, Faculty of Science,  
Kobe University, \\
Kobe, Rokko, 657-8501, Japan}
\email{mhsaito@math.kobe-u.ac.jp}

\begin{abstract}
In this paper, we show that the family of moduli 
spaces of  $\balpha'$-stable 
$(\bt, \blambda)$-parabolic $\phi$-connections of rank $2$ over $\BP^1$ with 
$4$-regular singular points and 
the fixed determinant bundle of degree $-1$ is isomorphic to the 
family of Okamoto--Painlev\'e pairs introduced by Okamoto \cite{O1} and 
\cite{STT02}. We also discuss about the generalization of our theory to 
the case where the rank of the connections and  genus of the  base curve 
are arbitrary.  Defining   isomonodromic flows on the family of 
moduli space of stable parabolic connections via the Riemann-Hilbert 
correspondences, we will show that a property  of the Riemann-Hilbert 
correspondences implies the Painlev\'e property of  isomonodromic flows. 
\end{abstract}
\maketitle

\section{Introduction}
In part I  \cite{IIS-1}, we established  
a complete geometric background for Painlev\'e equations of type VI or
more generally for Garnier systems from view points of  
moduli spaces of rank $2$ stable parabolic connections, moduli spaces of 
$SL_2$-representations of 
$\pi_1(\BP^1\setminus D(\bt))$ 
and the Riemann-Hilbert correspondences between them.

In this formulation, Painlev\'e equations of type  VI  or Garnier systems 
are vector fields or systems of vector fields 
on each corresponding  family   of moduli spaces of stable parabolic 
connections arising from isomonodromic deformations of linear connections. 
 Most notably, we can give a complete geometric  proof of the Painlev\'e 
property of Painlev\'e  equations of type  VI  
and Garnier systems by proving that 
the Riemann-Hilbert correspondences  
are {\it bimeromorphic proper surjective holomorphic maps}.  
  Moreover, one can prove that the Riemann-Hilbert 
correspondences give analytic resolutions of singularities of 
moduli spaces of the $SL_2$-representations. Then on the 
inverse image of each singular point, which is a family of 
compact subvarieties in the family of moduli spaces of connections,  
the vector fields admit classical solutions such as Riccati solutions 
in  Painlev\'e  VI  case. 
See \cite{Iwa02-1}, \cite{Iwa02-2}, \cite{SU01}, \cite{IIS-0}, \cite{STe02} and \cite{IIS-3}, for further applications of our approach to 
explicit dynamics of the Painlev\'e  VI  equations such as the
classification of Riccati solutions and rational solutions, nonlinear 
monodromy, and  B\"aklund transformations as well as the relation 
with the former results \cite{Miwa}, \cite{Mal}  on the 
Painlev\'e property.  

In this paper, with the notation in \S \ref{sec:moduli-spc}, 
we study in detail the moduli space
$
\overline{M_4^{\balpha'}}(\bt, \blambda, -1)
$ 
of $\balpha'$-stable 
$(\bt, \blambda)$-parabolic $\phi$-connections of rank $2$ over $\BP^1$ with 
the fixed determinant bundle of degree $-1$ as well as 
the moduli space 
$
M_4^{\balpha}(\bt, \blambda, -1)
$ 
of corresponding 
$\balpha$-stable $(\bt, \blambda)$-parabolic connections of rank $2$ 
over $\BP^1$.  
From a general result ([Theorem 1.1,  \cite{IIS-1}] or [Theorem 
\ref{thm:existence}, \S \ref{sec:moduli-spc}]) which is also 
valid for $n \geq 5$, we can show that  
\begin{itemize}
\item $\overline{M_4^{\balpha'}}(\bt, \blambda, -1)$ is a projective 
surface, 
\item $
M_4^{\balpha}(\bt, \blambda, -1)
$ is a smooth irreducible algebraic surface with a holomorphic 
symplectic structure and   
\item there exists a natural embedding 
$
M_4^{\balpha}(\bt, \blambda, -1) \hookrightarrow \overline{M_4^{\balpha'}}(\bt, \blambda, -1).  
$
\end{itemize}

In Theorem \ref{thm:n=4}, which is the main 
theorem in this paper,  we will 
show that  the moduli space $\overline{M_4^{\balpha'}}(\bt, \blambda, -1)$ is 
isomorphic to a {\em smooth}  projective rational surface $\overline{\cS}_{\bt,\blambda}$.  
Moreover we can show that 
there exists a unique effective anti-canonical divisor 
 $\cY_{\bt, \blambda} \in |- K_{\overline{\cS}_{\bt,\blambda}}|$ of 
 $\overline{\cS}_{\bt,\blambda}$ 
 such that 
 $ \overline{\cS}_{\bt, \blambda} \setminus \cY_{\bt, \blambda, red} \simeq M_4^{\balpha}(\bt, \blambda, -1)$.  Moreover 
 $(\overline{\cS}_{\bt, \blambda}, Y_{\bt, \blambda})$ is a non-fibered rational Okamoto--Painlev\'e pairs of type $D_{4}^{(1)}$  which is defined in \cite{STT02} (cf. \cite{Sakai}).  Note that $\overline{\cS}_{\bt, \blambda} \setminus \cY_{\bt, \blambda, red} $ is isomorphic to the space of initial conditions for 
 Painlev\'e equations of type VI constructed by Okamoto \cite{O1}. 

We should mention here that an algebraic moduli space of parabolic connections 
without stability conditions was essentially considered by D. Arinlin and S. Lysenco in \cite{AL1}, \cite{AL2} and \cite{A} and they constructed a nice moduli space  for generic $\blambda$.  
However for special $\blambda$, we should consider  certain stability condition to construct a nice moduli space.  
There are also different approaches   \cite{N}, \cite{Ni} for  constructions of moduli spaces of logarithmic connections with or without parabolic structures.  

The rough plan of this paper is as follows.  
In \S \ref{sec:op}, we will explain about motivation of 
this paper and the theory of Okamoto--Painlev\'e pairs in \cite{STa} and 
\cite{STT02}.     
In \S \ref{sec:moduli-spc}, we review results in 
part I \cite{IIS-1}.  In \S \ref{sec:n=4}, we will state 
Theorem \ref{thm:n=4} and the rest of the section will be devoted to 
show this theorem.  In \S \ref{sec:mod-conn}, we 
give a formulation of moduli theory of stable parabolic connection 
with regular singularities of any rank over any smooth curve. We also 
define the moduli space of  representations of the fundamental group of 
 $n$-punctured curve of genus $g$.   Then 
we state  the existence theorem of moduli space due to Inaba \cite{Ina} 
without proof.  In \S \ref{ss:RH},  
we define the Riemann-Hilbert correspondence and 
state, also without proof,  
Theorem \ref{thm:rh} which says that the Riemann--Hilbert 
correspondence is a proper surjective bimeromorphic analytic 
morphism. 
 In \S \ref{sec:isomon}, we will define isomonodromic 
flows on the family of the moduli spaces of $\balpha$-stable 
parabolic connections.  Assuming that Theorem \ref{thm:rh} is true, 
we will show  that isomonodromic flows satisfy
the Painlev\'e property. (Note that,  if rank $r = 2$ and over $\BP^1$, 
a proof of Theorem \ref{thm:rh} is found in \cite{IIS-1}).

Throughout in this paper, we will work over the field $\C$ of complex numbers.

\section{Motivation--Painlev\'e equations of type VI 
and  Okamoto--Painlev\'e pairs}
\label{sec:op}

Let us 
recall the theory of space of initial conditions of  Painlev\'e equation of type  VI.  Fix $\blambda = (\lambda_1, \cdots, \lambda_4) \in \Lambda_4 = \C^4$ and 
consider the following ordinary differential equation 
of Painlev\'e VI type $P_{VI}(\blambda)$ parameterized by  $\blambda$:
%\small % \hidewidth inserted  by S.M. 4/27/06
\begin{equation}\label{eq:painleve}
\begin{array}{lcl}
{P_{VI}(\blambda)} : \hidewidth & &  \\
\displaystyle{\frac{d^2 x}{d t^2}}  &= & 
\displaystyle{\frac{1}{2}\left( \frac{1}{x} + 
\frac{1}{x-1}+ \frac{1}{x-t} \right)
\left( \frac{d x}{d t} \right)^2 - } \\
 & & \displaystyle{ \left( \frac{1}{t} + \frac{1}{t-1}+ 
\frac{1}{x-t} \right)\left( \frac{d x}{d t} \right)} + \frac{x(x-1)(x-t)}{t^2(t-1)^2}\times   \\
         & &   
      \displaystyle{\left[2 (\lambda_4 -\frac{1}{2})^2 - 2 \lambda_1^2 \frac{t}{x^2} + 
         2\lambda_2^2  
         \frac{t-1}{(x-1)^2} + \left(\frac{1}{2}-  
       2\lambda_3^2 \right) \frac{t(t-1)}{(x-t)^2} \right]} . 
\end{array}
\end{equation}
\normalsize
It is known that this algebraic differential equation $P_{VI}(\blambda)$ 
is equivalent to the following nonautonomous Hamiltonian system:
\begin{equation}\label{eq:hamilton}
\renewcommand{\arraystretch}{2.2}
(H_{VI}(\blambda)): \left\{
\begin{array}{ccl}
\displaystyle{ \frac{d x}{d t}} & = &
\displaystyle{ \frac{ \partial H_{VI}}{ \partial y} }, \\
\displaystyle{ \frac{ d y}{d t}} & = & -
\displaystyle{ \frac{ \partial H_{VI}}{\partial x} },
\end{array}
\right.
\end{equation}
where the Hamiltonian is given as follows. %\small
\begin{eqnarray*}
H_{VI}(x, y, t) & = &  \displaystyle{ \frac{1}{t(t-1)}
\left[ x (x-1)(x-t) y^2 -
\left\{ 2 \lambda_1 (x -1) (x - t) \right. \right.} \\
 && \displaystyle{ \left. \left. + 2 \lambda_2 x(x-t) + (2\lambda_3 - 1) x (x
-1)
\right\} y + \lambda(x - t) \right]} 
\end{eqnarray*}
\begin{equation*}
\displaystyle{ \left( \lambda := 
\left\{ (\lambda_1 + \lambda_2 + \lambda_3 -1/2)^2 - (\lambda_{4}- \frac{1}{2})^2
\right\} \right) }.
\end{equation*}\normalsize
Let us set  $T = \C \setminus \{ 0, 1 \}$ and consider the 
following algebraic vector fields  on $\cS^{(0)} = \C^2 \times T \times \Lambda_4  \ni (x, y, t, \blambda)$
\begin{equation}\label{eq:alg-vf}
 v = \frac{\partial}{\partial t} +  \frac{ \partial H_{VI}}{ \partial y} \frac{ \partial}{\partial x} - 
\frac{\partial H_{VI}}{\partial x} \frac{ \partial }{\partial y}
\end{equation}
Taking  a relative compactification 
$\overline{\cS}^{(0)} = \Sigma_0 \times T \times \Lambda_4$ 
of $\cS^{(0)}$ where $\Sigma_0 = \BP^1 \times \BP^1$
and setting $\cD^{(0)}  = \overline{\cS}^{(0)} \setminus \cS^{(0)}$, we obtain 
the commutative diagram: 
\begin{equation}\label{eq:extend}
\begin{array}{ccccc}
\cS^{(0)} & \hookrightarrow & \overline{\cS}^{(0)} & \hookleftarrow & \cD^{(0)} \\
         &  \searrow \pi & \downarrow \overline{\pi}^{(0)} & \swarrow &  \\
&&  T \times \Lambda_4.  & &  \\ 
\end{array}
\end{equation}
We can extend the vector field $v$ in (\ref{eq:alg-vf}) on $ \cS^{(0)}$ 
to a rational vector field 
\begin{equation} \label{eq:orig-vf}
\tilde{v} \in H^0 ( \overline{\cS}^{(0)}, 
\Theta_{\overline{\cS}^{(0)}}(* \cD^{(0)} )). 
\end{equation}
In general, the rational vector field $ \tilde{v} $ 
has  accessible singularities
at the boundary divisor $\cD^{(0)} $. In  \cite{O1}, 
Okamoto gave explicit resolutions of accessible singularities by successive blowings-up at points on the boundary divisor.  Then finally, we obtain a smooth 
family of smooth projective rational surfaces  
\begin{equation}\label{eq:extend}
\begin{array}{ccccc}
\cS & \hookrightarrow & \overline{\cS} & \hookleftarrow & \cD \\
         &  \searrow \pi & \downarrow \overline{\pi} & \swarrow &  \\
&&  T \times \Lambda_4.  & &  \\ 
\end{array}
\end{equation}
such that $ \cD := \overline{\cS} \setminus \cS $ is a reduced normal crossing 
divisor and $ \cS $ contains $ \cS^{(0)} $ as a Zariski open set.  
Moreover one can show that 
\begin{equation}\label{eq:1-log}
\tilde{v} \in H^{0}(\overline{\cS}, \Theta_{\overline{\cS}}(- \log \cD) (\cD)), \end{equation}
where $\Theta_{\overline{\cS}}(- \log \cD)$ denotes the sheaf of germs of 
regular vector fields with logarithmic zero along $\cD$ (cf. \cite{STT02}).  
The extended rational vector field $\tilde{v}$ on $\overline{\cS}$ has 
poles of order $1$ along $\cD$ and is regular on $\cS = \overline{\cS} \setminus \cD$.  

For each fixed $(t, \blambda) \in T \times \Lambda_4$, the fiber 
$\overline{\pi}^{-1}((t, \blambda)) = \overline{\cS}_{t, \blambda}$ has 
a unique effective anti-canonical divisor 
$ \cY_{t, \blambda} \in |-K_{\overline{\cS}_{t, \blambda}}| $ with the irreducible decomposition
$$
\cY_{t, \blambda} = 2 D_0 + D_1 + D_2 + D_3 + D_4
$$
such that $\cY_{t, \blambda, red} = \sum_{i=0}^4 D_i = \cD_{t, \blambda}$.  
Moreover it satisfies the following numerical conditions 
\begin{equation}
\fbox{   $\cY_{t, \blambda} \cdot D_i = \deg(-K_{\overline{\cS}_{t, \blambda}|D_i})  = 0$ for $ i= 0, \ldots, 4$. }
\end{equation}

In \cite{STT02}, we give the following
\begin{Definition}\label{def:op}{\rm (Cf. \cite{STT02}, \cite{STa}, \cite{Sakai}). }
{\rm 
A pair $(S, Y)$ of a smooth projective rational surface with 
an anti-canonical divisor $Y \in |-K_S|$ with the irreducible
decomposition $Y= \sum_{i} m_i Y_i$ is called a rational 
{\em Okamoto--Painlev\'e pair} if it satisfies the condition 
\begin{equation}\label{eq:op-1}
\fbox{   $Y \cdot Y_i = \deg(-K_{\overline{\cS}_{t, \blambda}|Y_i})  = 0$ for all  $ i$. }
\end{equation}
A rational Okamoto--Painlev\'e pair $(S, Y)$ is called {\em of fibered-type} if 
there exists an elliptic fibration $f:S \lra \BP^1$ such that $f^{*}(\infty) = n Y$ for some $ n \geq 1 $.  }
\end{Definition}
It is easy to see that for a rational Okamoto--Painlev\'e pair 
the configuration of $Y$ is in the list of degenerate fibers of 
elliptic surfaces due to Kodaira, which was classified by 
affine Dynkin diagrams. Therefore,   
we have a classification of rational Okamoto--Painlev\'e pairs $(S, Y)$
by the Dynkin diagram of $Y$.  
For the case of Painlev\'e VI, we can say that the pair $(\overline{\cS}_{t, \blambda}, \cY_{t, \blambda})$ appeared in a fiber of the family 
(\ref{eq:extend}) is a rational Okamoto--Painlev\'e pair of type 
$D_4^{(1)}$.
The family of the complement of the divisor $\cD$ in (\ref{eq:extend}) 
 $\cS \lra T \times \Lambda_4$, where the rational 
vector field $\tilde{v}$ is regular,  should be  the family of 
the space of initial conditions of Painlev\'e equations of type VI or 
the phase space of 
the vector field  $\tilde{v}$.  Note that  $\cS \lra T \times \Lambda_4$
contains  the original family $\cS^{(0)} \lra T \times \Lambda_4$ as 
a {\em proper} Zariski open subset, that is, $\cS^{(0)} \subsetneq \cS$.  
Here we recall the following technical lemma proved in [Proposition 1.3, \cite{STT02}]. 

\begin{Lemma}\label{prop:cf}
Let $(S, Y)$ be a rational Okamoto--Painlev\'e pair. Then the following 
conditions are equivalent to each other. 
\begin{enumerate}
\item $(S, Y)$ is non-fibered type.
\item A regular algebraic functions on the complement $S \setminus Y_{red}$
must be a constant function.  
\end{enumerate}
In particular, for a non-fibered  rational Okamoto--Painlev\'e pair 
$(S, Y)$, the complement $S \setminus Y_{red}$ is never an affine variety.   
\end{Lemma}

Since one can show that an 
Okamoto--Painlev\'e pair 
$(\overline{\cS}_{t, \blambda}, \cY_{t, \blambda})$ which appeared in 
a fiber of $\overline{\pi}$ in (\ref{eq:extend}) 
is non-fibered type, 
we obtain the following 
\begin{Corollary}\label{cor:not-affine}
As for the family (\ref{eq:extend}) for Painlev\'e equations of type VI 
constructed by Okamoto \cite{O1}, each fiber 
$\cS_{t, \blambda} = \overline{\cS}_{t, \blambda} \setminus \cD_{t, \blambda} $
is not an affine variety.  
\end{Corollary}

In Theorem \ref{thm:n=4}, we will show that 
the family (\ref{eq:extend}) $\overline{\cS} \lra T \times \Lambda_4$ 
 constructed by Okamoto in \cite{O1} is 
isomorphic to the family of moduli spaces 
$$\overline{M_4^{\balpha'}}(-1)
\lra T_4 \times \Lambda_4
$$
of $\balpha'$-stable parabolic 
$\phi$-connections of rank $2$  over $\BP^1$ with 
$4$ regular singular points.  (In order to 
identify, we need to normalize $4$ points $(t_1, t_2, 
t_3, t_4) $ to $(0, 1, t, \infty)$). 

In \cite{IIS-1}, for $\ba= (a_1, \cdots, a_4) \in \cA_4 \simeq \C^4$, 
we can also consider
the moduli space $ \cR(\cP_{4,t})_{\ba}$ of 
$SL_2(\C)$-representations $\rho$ of $\pi_1(\BP^1 \setminus D(\bt))$ 
with the conditions $\Tr[\rho(\gamma_i)] = a_i$.  Then 
we can define the Riemann-Hilbert correspondence 
\begin{equation}\label{eq:rh-1}
\RH_{t, \blambda}: S_{t,\blambda} \simeq M_4^{\balpha}(t, \blambda, -1) \lra 
\cR(\cP_{4,t})_{\ba}
\end{equation} 
where $a_i = 2 \cos 2 \pi \lambda_i$. 

Note that the Riemann--Hilbert correspondence is a highly transcendental 
analytic morphism, which is  never an algebraic morphism. 
From results in \cite{IIS-1}, we can show the following Theorem, 
which shows highly transcendental nature of the Riemann--Hilbert 
correspondence $\RH_{t, \blambda}$. 

\begin{Proposition}{\rm (Cf. [Theorem 1.4, Theorem 1.3,  \cite{IIS-1}] )}
\begin{enumerate}
\item 
For all $(t, \blambda) \in T \times \Lambda_4$, the Riemann--Hilbert 
correspondence $\RH_{t, \blambda}$ is a bimeromorphic proper surjective 
analytic morphism.  If $\blambda \in \Lambda_4$ is generic, $\RH_{t, \blambda}$
is an analytic isomorphism.  
\item For all $\ba \in \cA_4$, $\cR(\cP_{4,t})_{\ba}$ is an affine variety, 
while $ S_{t,\blambda} \simeq M_4^{\balpha}(t, \blambda, -1)$ is not 
an affine variety. Hence if $\lambda \in \Lambda_4$ is generic, 
$\RH_{t, \blambda}$ gives an analytic isomorphism between a non-affine variety 
$ S_{t,\blambda} \simeq M_4^{\balpha}(t, \blambda, -1)$  and an affine 
variety $\cR(\cP_{4,t})_{\ba}$. 

\item For a generic $\blambda \in \Lambda_4$,  $ S_{t,\blambda} \simeq M_4^{\balpha}(t, \blambda, -1)$ is a Stein manifold, but not an affine variety.  

\end{enumerate}
\end{Proposition}

In \S \ref{sec:n=4}, in order to obtain  
Okamoto-Painlev\'e pairs $ (\overline{\cS}_{\bt, \blambda}, \cY_{\bt, \blambda})$, we use  a process of blowings-up which is a little bit different 
from Okamoto's in \cite{O1}.  
The process can be explained as follows.  Take $\Sigma_2 = \BP(\cO_{\BP^1}(2) 
\oplus \cO_{\BP^1}) \lra \BP^1$, which is  
the Hirzebruch surface of degree $2$. 
Let $D_0$ denote the unique infinite section with $D_0^2 = -2 $ and take 
the fibers $F_i $ over $t_i$ for $i=1, \ldots, 4$.  From the data 
$\lambda_i$, we can determine  two points $b_i^{+}$ and $b_i^{-}$ on $F_i$. 
(See \S \ref{sec:n=4} for precise definition of  $b_i^{\pm}$). 
By blowing-up of $\Sigma_2$ at 
$8$-points $\{ b_i^{\pm} \}_{i=1}^{4}$, we obtain the rational 
surface $\overline{\cS}_{\bt,\blambda}$  and 
the  unique effective anti-canonical divisor $\cY_{\bt, \blambda}$ 
can be given by $\cY_{\bt, \blambda} = 2 D_0 + D_1+ D_2 + D_3+D_4$
where $D_i$ denotes the proper transform of $F_i$, (see Fig. \ref{fig:op1}).

\begin{figure}[h]
%WinTpicVersion2.15
\unitlength 0.1in
\begin{picture}(39.25,31.25)(2.20,-34.15)
% LINE 0 0 3 0
% 2 3200 1010 3200 3010
% 
\special{pn 20}%
\special{pa 3200 610}%
\special{pa 3200 2610}%
\special{fp}%
% LINE 0 0 3 0
% 2 1400 1000 1400 3000
% 
\special{pn 20}%
\special{pa 1400 600}%
\special{pa 1400 2600}%
\special{fp}%
% LINE 0 0 3 0
% 2 2010 1000 2010 3000
% 
\special{pn 20}%
\special{pa 2010 600}%
\special{pa 2010 2600}%
\special{fp}%
% LINE 0 0 3 0
% 2 2600 1000 2600 2990
% 
\special{pn 20}%
\special{pa 2600 600}%
\special{pa 2600 2590}%
\special{fp}%
% BOX 0 0 3 0
% 2 820 1000 3810 3010
% 
\special{pn 20}%
\special{pa 820 600}%
\special{pa 3810 600}%
\special{pa 3810 2610}%
\special{pa 820 2610}%
\special{pa 820 600}%
\special{fp}%
% LINE 0 0 3 0
% 2 820 3610 3810 3610
% 
\special{pn 20}%
\special{pa 820 3210}%
\special{pa 3810 3210}%
\special{fp}%
% VECTOR 1 0 3 0
% 2 2290 3090 2290 3550
% 
\special{pn 13}%
\special{pa 2290 2690}%
\special{pa 2290 3150}%
\special{fp}%
\special{sh 1}%
\special{pa 2290 3150}%
\special{pa 2310 3083}%
\special{pa 2290 3097}%
\special{pa 2270 3083}%
\special{pa 2290 3150}%
\special{fp}%
% LINE 0 0 3 0
% 2 820 1420 3820 1420
% 
\special{pn 20}%
\special{pa 820 1020}%
\special{pa 3820 1020}%
\special{fp}%
% LINE 0 0 3 0
% 2 1400 3560 1400 3670
% 
\special{pn 20}%
\special{pa 1400 3160}%
\special{pa 1400 3270}%
\special{fp}%
% LINE 0 0 3 0
% 2 3190 3550 3190 3660
% 
\special{pn 20}%
\special{pa 3190 3150}%
\special{pa 3190 3260}%
\special{fp}%
% LINE 0 0 3 0
% 2 2000 3560 2000 3670
% 
\special{pn 20}%
\special{pa 2000 3160}%
\special{pa 2000 3270}%
\special{fp}%
% LINE 0 0 3 0
% 2 2610 3560 2610 3670
% 
\special{pn 20}%
\special{pa 2610 3160}%
\special{pa 2610 3270}%
\special{fp}%
% LINE 0 2 3 0
% 2 1810 1610 2200 2020
% 
\special{pn 20}%
\special{pa 1810 1210}%
\special{pa 2200 1620}%
\special{dt 0.054}%
\special{pa 2200 1620}%
\special{pa 2200 1620}%
\special{dt 0.054}%
% LINE 0 2 3 0
% 2 2400 1600 2790 2010
% 
\special{pn 20}%
\special{pa 2400 1200}%
\special{pa 2790 1610}%
\special{dt 0.054}%
\special{pa 2790 1610}%
\special{pa 2790 1610}%
\special{dt 0.054}%
% LINE 0 2 3 0
% 2 3010 1620 3400 2030
% 
\special{pn 20}%
\special{pa 3010 1220}%
\special{pa 3400 1630}%
\special{dt 0.054}%
\special{pa 3400 1630}%
\special{pa 3400 1630}%
\special{dt 0.054}%
% LINE 0 2 3 0
% 2 1200 1610 1590 2020
% 
\special{pn 20}%
\special{pa 1200 1210}%
\special{pa 1590 1620}%
\special{dt 0.054}%
\special{pa 1590 1620}%
\special{pa 1590 1620}%
\special{dt 0.054}%
% LINE 0 2 3 0
% 2 1600 2410 1200 2800
% 
\special{pn 20}%
\special{pa 1600 2010}%
\special{pa 1200 2400}%
\special{dt 0.054}%
\special{pa 1200 2400}%
\special{pa 1200 2400}%
\special{dt 0.054}%
% LINE 0 2 3 0
% 2 2190 2410 1790 2800
% 
\special{pn 20}%
\special{pa 2190 2010}%
\special{pa 1790 2400}%
\special{dt 0.054}%
\special{pa 1790 2400}%
\special{pa 1790 2400}%
\special{dt 0.054}%
% LINE 0 2 3 0
% 2 2800 2400 2400 2790
% 
\special{pn 20}%
\special{pa 2800 2000}%
\special{pa 2400 2390}%
\special{dt 0.054}%
\special{pa 2400 2390}%
\special{pa 2400 2390}%
\special{dt 0.054}%
% LINE 0 2 3 0
% 2 3390 2390 2990 2780
% 
\special{pn 20}%
\special{pa 3390 1990}%
\special{pa 2990 2380}%
\special{dt 0.054}%
\special{pa 2990 2380}%
\special{pa 2990 2380}%
\special{dt 0.054}%
% STR 2 0 3 0
% 3 340 1940 340 2040 2 0
% $\varSigma_2$
\put(3.4000,-16.4000){\makebox(0,0)[lb]{$\overline{\cS}_{\bt, \blambda}$}}%
% STR 2 0 3 0
% 3 1390 680 1390 780 5 0
% $E_1$
\put(13.9000,-3.8000){\makebox(0,0){$D_1$}}%
% STR 2 0 3 0
% 3 2030 680 2030 780 5 0
% $E_2$
\put(20.3000,-3.8000){\makebox(0,0){$D_2$}}%
% STR 2 0 3 0
% 3 2610 680 2610 780 5 0
% $E_3$
\put(26.1000,-3.8000){\makebox(0,0){$D_3$}}%
% STR 2 0 3 0
% 3 3210 680 3210 780 5 0
% $E_4$
\put(32.1000,-3.8000){\makebox(0,0){$D_4$}}%
% STR 2 0 3 0
% 3 4130 1320 4130 1420 5 0
% $E_0$
\put(41.3000,-10.2000){\makebox(0,0){$D_0$}}%
% STR 2 0 3 0
% 3 1400 3780 1400 3880 5 0
% $t_1$
\put(14.0000,-34.8000){\makebox(0,0){$t_1$}}%
% STR 2 0 3 0
% 3 2000 3790 2000 3890 5 0
% $t_2$
\put(20.0000,-34.9000){\makebox(0,0){$t_2$}}%
% STR 2 0 3 0
% 3 2610 3790 2610 3890 5 0
% $t_3$
\put(26.1000,-34.9000){\makebox(0,0){$t_3$}}%
% STR 2 0 3 0
% 3 3190 3800 3190 3900 5 0
% $t_4$
\put(31.9000,-35.0000){\makebox(0,0){$t_4$}}%
% STR 2 0 3 0
% 3 2450 3180 2450 3280 5 0
% $\pi$
\put(24.5000,-28.8000){\makebox(0,0){$\pi$}}%
% STR 2 0 3 0
% 3 490 3520 490 3620 5 0
% $\BP^1$
\put(4.9000,-32.2000){\makebox(0,0){$\BP^1$}}%
% STR 2 0 3 0
% 3 3810 760 3810 860 2 0
% $\infty$-section
\put(38.1000,-4.6000){\makebox(0,0)[lb]{$\infty$-section}}%
% VECTOR 2 0 3 0
% 2 4140 920 4140 1280
% 
\special{pn 8}%
\special{pa 4140 520}%
\special{pa 4140 880}%
\special{fp}%
\special{sh 1}%
\special{pa 4140 880}%
\special{pa 4160 813}%
\special{pa 4140 827}%
\special{pa 4120 813}%
\special{pa 4140 880}%
\special{fp}%
% CIRCLE 0 0 0 0
% 4 3200 1820 3200 1855 3200 1855 3200 1855
% 
\special{pn 20}%
\special{sh 0.600}%
\special{ar 3200 1420 35 35  0.0000000 6.2831853}%
% CIRCLE 0 0 0 0
% 4 3210 2580 3210 2615 3210 2615 3210 2615
% 
\special{pn 20}%
\special{sh 0.600}%
\special{ar 3210 2180 35 35  0.0000000 6.2831853}%
% CIRCLE 0 0 0 0
% 4 2600 2600 2600 2635 2600 2635 2600 2635
% 
\special{pn 20}%
\special{sh 0.600}%
\special{ar 2600 2200 35 35  0.0000000 6.2831853}%
% CIRCLE 0 0 0 0
% 4 2610 1830 2610 1865 2610 1865 2610 1865
% 
\special{pn 20}%
\special{sh 0.600}%
\special{ar 2610 1430 35 35  0.0000000 6.2831853}%
% CIRCLE 0 0 0 0
% 4 2010 1840 2010 1875 2010 1875 2010 1875
% 
\special{pn 20}%
\special{sh 0.600}%
\special{ar 2010 1440 35 35  0.0000000 6.2831853}%
% CIRCLE 0 0 0 0
% 4 2010 2600 2010 2635 2010 2635 2010 2635
% 
\special{pn 20}%
\special{sh 0.600}%
\special{ar 2010 2200 35 35  0.0000000 6.2831853}%
% CIRCLE 0 0 0 0
% 4 1400 1830 1400 1865 1400 1865 1400 1865
% 
\special{pn 20}%
\special{sh 0.600}%
\special{ar 1400 1430 35 35  0.0000000 6.2831853}%
% CIRCLE 0 0 0 0
% 4 1400 2600 1400 2635 1400 2635 1400 2635
% 
\special{pn 20}%
\special{sh 0.600}%
\special{ar 1400 2200 35 35  0.0000000 6.2831853}%
\end{picture}%
\label{fig:op1}
\caption{Okamoto--Painlev\'e pair of type $D_4^{(1)}$}
\label{fig:op1}
\end{figure}
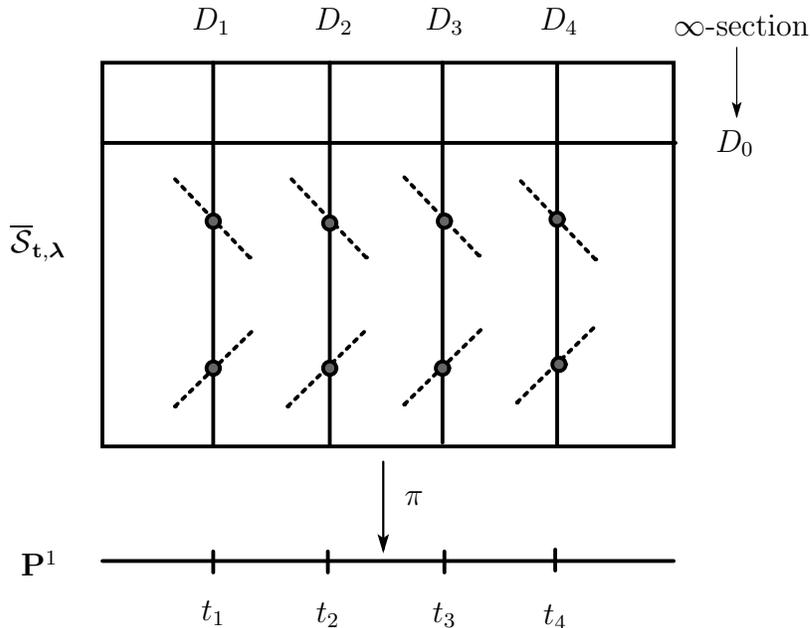

%%%%%%%%%%%%%%%%%%% sec:moduli-spc %%%%%%%%%%%%%%

\section{Moduli spaces of rank 2 stable parabolic 
connections on $\BP^1$ and their compactifications. A review of Part I. }
\label{sec:moduli-spc}
%%%%%%%%%%%%%%%%%%%%%%

In this section,  we reproduce basic notation and definition in part I \cite{IIS-1} for  reader's convenience. 

\subsection{Parabolic connections on $\BP^1$.}

Let $n \geq 3$ and set
\begin{equation}\label{eq:config-space}
T_n = \{ ( t_1, \ldots, t_n) \in (\BP^1)^n \quad | \quad t_i \not=t_j, (i \not= j) 
\}, 
\end{equation}
\begin{equation}\label{eq:exponents-space}
\Lambda_{n}  = \{ \blambda = (\lambda_1, \ldots, \lambda_n) \in  \C^{n} \}. 
\end{equation}  
Fixing a data  $(\bt, {\blambda}) = (t_1, \ldots, t_n, \lambda_1, \ldots, \lambda_n)  
\in T_n \times \Lambda_n$, we define a reduced divisor on $\BP^1$ as 
\begin{equation} \label{eq:divisor}
D(\bt) = t_1 + \cdots + t_n.
\end{equation}  
Moreover we fix a line bundle $L$ on $\BP^1$ with a logarithmic 
connection $\nabla_L: L 
\lra L \otimes \Omega^1_{\BP^1}(D(\bt))$.

\begin{Definition}\label{def:parabolic} \rm 
A (rank $2$) $(\bt, \blambda)$-parabolic connection on $\BP^1$ with the determinant 
$(L, \nabla_L)$ is a quadruplet
$(E, \nabla, \varphi, \{l_i\}_{1 \leq i \leq n}  )$ 
which consists of  
\begin{enumerate}
\item 
a rank 2 vector bundle 
$E$ on $\BP^1$, 
\item 
a logarithmic connection $\nabla:E \lra E 
\otimes \Omega^1_{\BP^1}(D(\bt)) $

\item 
a bundle  isomorphism $ \varphi: \wedge^2 E \stackrel{\simeq}{\lra} L  $ 

\item  one dimensional subspace $l_i$ of the fiber $ E_{t_i}$  of $E$ 
at $t_i$,   $l_i \subset E_{t_i}$,   $ i = 1, \ldots, n$,   
such that 

\begin{enumerate}
\item  for any local sections $s_1, s_2$ of $E$, 
$$
\varphi \otimes id ( \nabla s_1 \wedge s_2 + s_1 \wedge \nabla s_2)  = 
\nabla_L(\varphi(s_1 \wedge s_2)), 
$$

\item  
$l_i \subset  \Ker (\res_{t_i}(\nabla) - \lambda_i)$, that is, $\lambda_i$ is an eigenvalue of the residue $\res_{t_i}(\nabla)$ of $\nabla$ at $t_i$ and 
$l_i$ is a one-dimensional eigensubspace of $\res_{t_i}(\nabla)$.   
\end{enumerate}
\end{enumerate}
\end{Definition}

\begin{Definition}\label{def:isom-spb} \rm 
  Two $(\bt, \blambda)$-parabolic connections 
$$
(E_1, \nabla_1, \varphi, \{l_i\}_{1 \leq i \leq n}  ), \quad  
(E_2, \nabla_2, \varphi', \{l'_i\}_{1 \leq i \leq n}  )
$$
 on $\BP^1$ with the determinant 
$(L, \nabla_L)$ are isomorphic to each other if 
there is an isomorphism $\sigma: E_1 \stackrel{\sim}{\lra} E_2 $ and $c \in \C^{\times}$
 such that the diagrams
 \begin{equation}
\begin{CD}
  E_1 @>\nabla_1>> E_1\otimes\Omega^1_{\BP^1}(D(\bt)) \\
 @V\sigma V\cong V  @V\cong V \sigma \otimes\mathrm{id}V \\
  E_2 @>\nabla_2>> E_2\otimes\Omega^1_{\BP^1}(D(\bt)) 
 \end{CD}
 \hspace{50pt}
 \begin{CD}
  \bigwedge^2 E_1 @>\varphi>\cong > L \\
  @V\wedge^2\sigma V\cong V @VcV\cong V \\
  \bigwedge^2 E_2 @>\varphi'>\cong > L 
 \end{CD}
\end{equation}
commute and $(\sigma)_{t_i}(l_i)=l'_i$ for $i=1,\ldots,n$. 
\end{Definition}

\subsection{The set of local exponents  $\blambda \in \Lambda_n$}

Note that a data 
$ \blambda = (\lambda_1, \ldots, \lambda_n) 
\in \Lambda_n \simeq \C^n $
specifies the set of eigenvalues of the residue 
matrix of a connection $\nabla$ at $\bt=(t_1, \ldots, t_n)$, 
which will be called a set of   {\em local exponents} of $\nabla$. 

\begin{Definition}\label{def:exponents}
{\rm 
A set of local exponents 
 $ \blambda =(\lambda_1, \ldots, \lambda_n) \in \Lambda_n $ is called {\em special} if  
\begin{enumerate}
\item $\blambda$ is {\em resonant}, that is, for some $ 1 \leq i \leq n$, 
\begin{equation} \label{eq:special1}
  2 \lambda_i \in \Z, 
\end{equation}  
\item or $\blambda$ is {\em reducible}, that is, for some 
$ (\epsilon_1, \ldots, \epsilon_n) \in \{ \pm 1 \}^n $ 
\begin{equation} 
\label{eq:special2}\sum_{i=1}^{n} \epsilon_i  \lambda_i  \in \Z. 
\end{equation}
\end{enumerate}

If $\blambda \in \Lambda_n$ is not special,  
$\blambda$ is said to be  
{\em generic}. 
}
\end{Definition}

\subsection{Parabolic degrees and $\balpha$-stability}

Let us fix a series of positive rational  
numbers $\balpha = (\alpha_1, \alpha_2, 
\ldots, \alpha_{2n} )$, which is called  
{\em a weight},   such that 
\begin{equation}
0 \leq \alpha_1 < \alpha_2 < \cdots< \alpha_{i} < \cdots < \alpha_{2n} < \alpha_{2n+1}=1.  
\end{equation} 
For a $(\bt, \blambda)$-parabolic connection on $\BP^1$ with the determinant 
$(L, \nabla_L)$, we can define the parabolic degree of 
$E=(E, \nabla, \varphi, l)$ with respect to the weight $\balpha$ by  
\begin{eqnarray}\label{eq:pardeg}
\pardeg_{\balpha} E  &  = &  \deg E + 
\sum_{i=1}^n \left(\alpha_{2i -1} \dim E_{t_i}/l_{i} + 
\alpha_{2i} \dim l_i \right) \\ 
& =  &\deg L + \sum_{i=1}^n 
(\alpha_{2i-1} + \alpha_{2i} ).  \nonumber
\end{eqnarray}
Let $F \subset E$ be a rank 1 subbundle of $E$ such that $ 
\nabla F \subset F \otimes \Omega^1_{\BP^1}(D(\bt))$.  
We define the parabolic degree of $(F, \nabla_{|F})$ by
\begin{equation}\label{eq:subpardeg} 
\pardeg_{\balpha} F  =   
\deg F +  \sum_{i=1}^n \left(\alpha_{2i-1} \dim F_{t_i}/l_{i} \cap F_{t_i} + 
\alpha_{2i} \dim l_i \cap F_{t_i}
\right).
\end{equation}

\begin{Definition} {\rm Fix a weight $\balpha$.  A $(\bt, \blambda)$-parabolic 
connection $(E, \nabla, \varphi, l)$ on $\BP^1$ with the determinant $(L, \nabla_L)$  
  is said to be  {\em $\balpha$-stable}
(resp. {\em $\balpha$-semistable} ) if 
for every rank-1 subbundle $F$ with $\nabla(F) \subset F \otimes \Omega^1_{\BP^1}
(D(\bt) )$ 
\begin{equation}\label{eq:stability}
\pardeg_{\balpha} F < \frac{\pardeg_{\balpha} E}{2}, \quad (\mbox{resp.} 
\pardeg_{\balpha} F \leq  \frac{\pardeg_{\balpha} E}{2} ).
\end{equation}}
(For simplicity, ``$\balpha$-stable" will be abbreviated to ``stable").
\end{Definition}

We define the coarse moduli space  by  

\begin{equation} \label{eq:modulispace}
M_n^{\balpha} (\bt, \blambda, L)= 
\left\{ (E, \nabla, \varphi, l ); 
\begin{array}{l}
\mbox{an $\balpha$-stable 
$(\bt, \blambda)$-parabolic } \\
\mbox{connection  with } \\
\mbox{the determinant $(L, \nabla_L)$ }
\end{array}
\right\}  / \mbox{isom.}  
\end{equation}

\subsection{Stable parabolic $\phi$-connections}

If $n \geq 4$, the moduli 
space $M_n^{\balpha} (\bt, \blambda, L)$ never becomes projective 
nor complete. 
In order to obtain a compactification of the moduli 
space $M_n^{\balpha} (\bt, \blambda, L)$, we will introduce  
the notion of a stable parabolic $\phi$-connection, or equivalently, 
a stable parabolic $\Lambda$-triple.  
Again, let us fix $(\bt,\blambda)\in T_n\times\Lambda_n$
and a line bundle $L$ on ${\bf P}^1$ with a connection
$\nabla_L:L\ra L\otimes\Omega^1_{\BP^1}(D(\bt))$.

%%%%%%%%%%%%%%%%%%%%%%%%%%%%%% Definition %%%%%%%%%%%%%%%%%%%%%%%%%%%%%%%%%%%
\begin{Definition}\rm The data 
 $(E_1,E_2,\phi,\nabla,\varphi,\{l_i\}_{i=1}^n)$
 is said to be a $(\bt,\blambda)$-parabolic $\phi$-connection
 of rank $2$ with the determinant $(L,\nabla_L)$ if
 $E_1,E_2$ are rank $2$ vector bundles on $\BP^1$
 with $\deg E_1=\deg L$,
 $\phi:E_1\ra E_2$,
 $\nabla:E_1\ra E_2\otimes\Omega^1_{\BP^1}(D(\bt))$
 are morphisms of sheaves,
 $\varphi:\bigwedge^2 E_2\stackrel{\sim}\lra L$
 is an isomorphism
 and $l_i\subset (E_1)_{t_i}$ are one dimensional subspaces
 for $i=1,\ldots,n$ such that
\begin{enumerate}
 \item $\phi(fa)=f\phi(a)$ and
  $\nabla(fa)=\phi(a)\otimes df+f\nabla(a)$
  for $f\in\cO_{\BP^1}$, $a\in E_1$,
 \item $(\varphi\otimes\mathrm{id})
 (\nabla(s_1)\wedge \phi(s_2)+\phi(s_1)\wedge\nabla(s_2))
 =\nabla_L(\varphi(\phi(s_1)\wedge\phi(s_2)))$
 for $s_1,s_2\in E_1$ and
 \item $(\res_{t_i}(\nabla)-\lambda_i\phi_{t_i})|_{l_i}=0$
  for $i=1,\ldots,n$.
\end{enumerate}
\end{Definition}

%%%%%%%%%%%%%%%%%%%%%%%%%%% Remark %%%%%%%%%%%%%%%%%%%%%%%%%%%%%%%%%%%%%%%%%%%
\begin{Definition}\rm
\item Two $(\bt,\blambda)$ parabolic $\phi$-connections
$$(E_1,E_2,\phi,\nabla,\varphi,\{l_i\}), \quad 
(E'_1,E'_2,\phi',\nabla',\varphi',\{l'_i\})
$$
 are said  to be isomorphic to each other  if
 there are isomorphisms $\sigma_1:E_1\stackrel{\sim}\lra E'_1$,
 $\sigma_2:E_2\stackrel{\sim}\lra E'_2$ and $c\in\C\setminus\{0\}$
 such that the diagrams
$$
 \begin{CD}
  E_1 @>\phi>> E_2 \\
  @V\sigma_1V\cong V  @V\cong V\sigma_2V \\
  E'_1 @>\phi'>> E'_2
 \end{CD}
 \hspace{50pt}
 \begin{CD}
  E_1 @>\nabla>> E_2\otimes\Omega^1_{\BP^1}(D(\bt)) \\
  @V\sigma_1V\cong V  @V\cong V\sigma_2\otimes\mathrm{id}V \\
  E'_1 @>\nabla'>> E'_2\otimes\Omega^1_{\BP^1}(D(\bt)) 
 \end{CD}
 $$
 $$
 \begin{CD}
  \bigwedge^2 E_2 @>\varphi>\cong > L \\
  @V\wedge^2\sigma_2V\cong V @VcV\cong V \\
  \bigwedge^2 E'_2 @>\varphi'>\cong > L 
 \end{CD}
$$
commute and $(\sigma_1)_{t_i}(l_i)=l'_i$ for $i=1,\ldots,n$.
\end{Definition}

\begin{Remark} \rm 
Assume that two vector bundles $E_1,E_2$ and 
 morphisms $\phi:E_1\ra E_2$,
 $\nabla:E_1\ra E_2\otimes\Omega^1_{\BP^1}(D(\bt))$ satisfying
 $\phi(fa)=f\phi(a)$, $\nabla(fa)=\phi(a)\otimes df+f\nabla(a)$
 for $f\in\cO_{\BP^1}$, $a\in E_1$ are given.
 If $\phi$ is an isomorphism, then
 $(\phi\otimes\mathrm{id})^{-1}\circ\nabla:
 E_1\ra E_1\otimes\Omega^1_{\BP^1}(D(\bt))$
 becomes a connection on $E_1$.
\end{Remark}

Fix rational numbers $ \alpha'_1, \alpha'_2, \ldots, \alpha'_{2n}, \alpha'_{2n+1} $
satisfying
$$
 0 \leq \alpha'_1 <  \alpha'_2  <  \cdots  <   \alpha'_{2n} <  \alpha'_{2n+1}=1
$$
and positive  integers $\beta_1, \beta_2$. 
 Setting $\balpha'=(\alpha'_1, \ldots, \alpha'_{2n}), \bbeta =(\beta_1, \beta_2)$, 
we obtain  a {\em weight}  $(\balpha', \bbeta) $ for  parabolic $\phi$-connections.   

%%%%%%%%%%%%%%%%%%%%%%%%%% Definition %%%%%%%%%%%%%%%%%%%%%%%%%%%%%%%%%%%%%

\begin{Definition}\rm Fix a sufficiently large integer $\gamma $. 
Let 
$$(E_1,E_2,\phi,\nabla,\varphi,\{l_i\}_{i=1}^n)$$ be a parabolic 
$\phi$-connection. For any  subbundles
 $F_1\subset E_1$, $F_2\subset E_2$ satisfying $\phi(F_1)\subset F_2$,
 $\nabla(F_1)\subset F_2\otimes\Omega^1_{\BP^1}(D(\bt))$, 
 we define 
%\small
\begin{multline*}
\mu((F_1,F_2))_{\balpha'\bbeta}
 = 
\frac{1}{\beta_1\rank(F_1)+\beta_2 \rank(F_2)}
(\beta_1(\deg F_1(-D(\bt)))  \\
+ \beta_2 ( \deg F_2-\gamma\rank(F_2))+
  \sum_{i=1}^n \beta_1 (\alpha'_{2i-1}d_{2i-1}(F_1)+ 
  \alpha'_{2i}d_{2i}(F_1))
 \end{multline*}
 where $d_{2i-1}(F)=\dim ((F_1)_{t_i}/l_i\cap (F_1)_{t_i})$,
 $d_{2i}(F_1)=\dim ((F_1)_{t_i}\cap l_i)$. 

A parabolic $\phi$-connection 
$(E_1,E_2,\phi,\nabla,\varphi,\{l_i\}_{i=1}^n)$
is said to be $(\balpha', \bbeta)$-stable (resp. $(\balpha', \bbeta)$-semistable) if for any subbundles
 $F_1\subset E_1$, $F_2\subset E_2$ satisfying $\phi(F_1)\subset F_2$,
 $\nabla(F_1)\subset F_2\otimes\Omega^1_{\BP^1}(D(\bt))$ and
 $(F_1,F_2)\neq (E_1,E_2),(0,0)$,
 the inequality
%\begin{equation}\label{eq:phi-stable}
\begin{eqnarray}\label{eq:phi-stable}
\mu((F_1, F_2))_{\balpha' \bbeta}  < \mu((E_1,E_2))_{\balpha'\bbeta},
\\ \nonumber
\quad ( \mbox{resp.} \quad \mu((F_1, F_2))_{\balpha' \bbeta}  \leq  
\mu((E_1,E_2))_{\balpha'\bbeta} .)
%\end{equation}
\end{eqnarray}
\end{Definition}

We define the coarse moduli space of 
 $(\balpha',\bbeta)$-stable  $(\bt,\blambda)$-parabolic 
 $\phi$-connections with the determinant $(L,\nabla_L)$ by 
\begin{equation}\label{eq:coarse-moduli-compact}
\overline{M_n^{\balpha' \bbeta}}(\bt,\blambda,L):=
 \left\{ (E_1,E_2,\phi,\nabla,\varphi,\{l_i\}) \right\}/\mathrm{isom}.
\end{equation}
For a given weight $(\balpha', \bbeta)$ and $1 \leq i \leq 2n$, 
define a rational number  $\alpha_i $ by 
\begin{equation}\label{eq:rel-weight}
\alpha_{i}=  \frac{\beta_1}{\beta_1 + \beta_2}\alpha'_i.
\end{equation}
Then $\balpha=(\alpha_i) $ satisfies the condition
\begin{equation}\label{eq:inequality}
0 \leq \alpha_1 < \alpha_2 < \cdots < \alpha_{2n}<  \frac{\beta_1}{(\beta_1+ \beta_2)} < 1, 
\end{equation}
hence $\balpha$ defines a weight for parabolic connections.
It is easy to see that if we take $\gamma $ sufficiently large
 $(E,\nabla,\varphi,\{l_i\})$ is $\balpha$-stable if and only if the 
 associated parabolic $\phi$-connection $(E,E,\mathrm{id}_E,\nabla,\varphi,\{l_i\})$
 is stable with respect to $(\balpha', \bbeta)$.  
 Therefore we  see that the natural map 
\begin{equation}\label{eq:map}
(E,\nabla,\varphi,\{l_i\})  \mapsto 
 (E,E,\mathrm{id}_E,\nabla,\varphi,\{l_i\})
\end{equation}
induces  an injection
\begin{equation}\label{eq:injection}
M_n^{\balpha}(\bt,\blambda,L)  \hookrightarrow 
\overline{M_n^{\balpha' \bbeta}}(\bt,\blambda,L).
\end{equation}
Conversely, assuming  that $\bbeta =(\beta_1, \beta_2)$ are given,  
for a weight $\balpha=(\alpha_i)$  satisfying the condition (\ref{eq:inequality}), we can define 
$\alpha'_i = \alpha_i \frac{\beta_1 + \beta_2}{\beta_1}$ for $1 \leq i \leq 2n$.  Since  
$0 \leq \alpha'_1 < \alpha'_2 < \cdots < \alpha'_{2n} = \alpha_{2n} \frac{\beta_1 + \beta_2}{\beta_1}< 1$, 
 $(\balpha', \bbeta)$ give a weight for parabolic $\phi$-connections.

Moreover, considering the relative setting over 
$T_n \times \Lambda_n$, we can  define  two families of 
the moduli spaces
\begin{equation}\label{eq:fam-moduli}
\overline{\pi}_n:
\overline{M_n^{\balpha'\bbeta}}(L) \lra T_n \times \Lambda_n, \quad 
\pi_n: M_n^{\balpha}(L) \lra T_n \times \Lambda_n
\end{equation}
such that the following diagram commutes;
\begin{equation}\label{eq:family-compact}
\begin{CD} 
M_n^{\balpha}(L) & \stackrel{\iota}{\hookrightarrow} & \overline{M_n^{\balpha' \bbeta }}(L)   \\
 @V \pi_n VV   @VV \overline{\pi}_n V  \\
 T_n \times \Lambda_n  @=    T_n \times \Lambda_n.  \\  
\end{CD}
\end{equation}
Here the fibers of $\pi_n$ and $\overline{\pi}_n $ 
over $(\bt, \blambda) \in T_n \times 
\Lambda_n $ are 
\begin{equation}
\pi_n^{-1}(\bt, \blambda) = M^{\balpha} (\bt, \blambda, L), \quad 
\overline{\pi}_n^{-1}(\bt, \blambda) = \overline{M^{\balpha' \bbeta}} (\bt, \blambda, L). 
\end{equation}

\subsection{The existence of moduli spaces and their properties}
\quad 

The following  theorem was proved in \cite{IIS-1}.

%%%%%%%%%%%%%%%%%% MAIN THEOREM %%%%%%%%%%%%%%%%%%%%%%%%%%%%%%%%%%%%%%%%%%%%%%
\begin{Theorem}$($ {\rm [Theorem 2.1, \cite{IIS-1}]}$)$.  \label{thm:fund}
\begin{enumerate}
\item Fix a weight $\bbeta =(\beta_1, \beta_2)$.  
For a generic weight $\balpha'$,
$$
\overline{\pi_n}:\overline{M_n^{\balpha' \bbeta }}(L)  \lra T_n \times \Lambda_n$$
is a {\em projective} morphism. In particular, the moduli space \linebreak
$\overline{M^{\balpha'\bbeta}}(\bt, \blambda, L)$ is a {\em projective} 
algebraic scheme for all 
$(\bt, \blambda) \in T_n \times \Lambda_n$.

\item For a generic  weight $\balpha$, 
$\pi_n:M_n^{\balpha}(L)  \lra T_n \times \Lambda_n$
is a {\em smooth morphism} of relative dimension $2n-6$ with irreducible closed fibers. 
Therefore, the moduli space 
$M^{\balpha}_n(\bt, \blambda, L)$ is a {\em smooth, irreducible} algebraic 
variety of dimension $2n -6$  for all 
$(\bt, \blambda) \in T_n \times \Lambda_n$.  
\end{enumerate}
\end{Theorem}
%%%%%%%%%%%%%%%%%% END OF  MAIN THEOREM %%%%%%%%%%%%%%%%%%%%%%%%%%%%%%%%%%%%%%%%%%%%%%

\begin{Remark}\label{rem:fund-thm} 
\rm 
\begin{enumerate}
\item The structures of 
moduli spaces $M^{\balpha}_n(L)$ and $\overline{M^{\balpha' \bbeta}_n}(L)$ may 
depend on the weights $\balpha$, $(\balpha', \bbeta)$ and $\deg L$.  

\item The moduli spaces $M^{\balpha}_n(L)$ is a fine moduli space. 
In fact, we have the universal families over these moduli 
spaces.

\item   The moduli space 
$M^{\balpha}_n(\bt, \blambda, L)$ admits a natural holomorphic 
symplectic structure. (See [Proposition 6.2, \cite{IIS-1}). 
This fact is a part of the reason why 
Painlev\'e  VI  and Garnier systems can be written in 
nonautonomous  Hamiltonian systems.  

\item In case of $n=4$, we can show that $ 
\overline{M^{\balpha'\bbeta}_4}(\bt, \blambda, L)$ is  smooth (cf. Proposition \ref{prop:smooth} ). 
However we do not know whether 
$ \overline{M^{\balpha'\bbeta}_n}(\bt, \blambda, L)$  is smooth or not for $ n \geq 5$.  

\end{enumerate}
\end{Remark}

When we describe the explicit algebraic or geometric 
structure of the moduli spaces 
$M^{\balpha}_n(L)$ and $\overline{M^{\balpha' \bbeta}_n}(L)$, 
it is convenient to fix a 
determinant line bundle $(L, \nabla_L)$.  
As a typical example of the determinant bundle is 
\begin{equation}\label{eq:typical}
(L, \nabla_L) = (\cO_{\BP^1}(-t_n),  d) 
\end{equation}
where the connection is given by 
\begin{equation}
\nabla_L (z - t_n) = d(z- t_n) = (z - t_n) \otimes \frac{dz}{z - t_n}. 
\end{equation}
Here   $z$  is an inhomogeneous coordinate of $\BP^1 = \Spec \C[z] \cup \{ \infty\}$.  
For this $(L, \nabla_L) = (\cO_{\BP^1}(-t_n), d)$, 
we set 
$$
M_n^{\balpha} (\bt, \blambda, -1) = M_n^{\balpha} (\bt, \blambda, L), \quad  (\mbox{resp.} \  
\overline{M_n^{\balpha'\bbeta}} (\bt, \blambda,-1)= \overline{M_n^{\balpha'\bbeta}}
 (\bt, \blambda, L) \quad ). 
$$

%%%%%%%%%%%%%%%%%subsection n=4 %%%%%%%%%%%%%%%%
\section{Explicit construction of moduli spaces 
for the case of  $n = 4$ (Painlev\'e  VI  case).}
\label{sec:n=4}

In this section, we will deal with the case of $n=4$ in detail.  
Let us fix a sufficiently large integer 
$\gamma$ and take a weight $(\balpha', \bbeta)$ 
for parabolic $\phi$-connections where  $\balpha' = 
(\alpha'_1, \ldots, \alpha'_8)$, $\bbeta =(\beta_1, \beta_2), \gamma$  and fix
$(\bt, \blambda) = (t_1, \ldots, t_4, \lambda_1, \ldots,  \lambda_4)
\in T_4 \times \Lambda_4$.  

Then the corresponding weight 
$\balpha = (\alpha_1, \ldots, \alpha_8)$ for parabolic 
connections can be given by 
$$
\alpha_i = \alpha'_i \frac{\beta_1}{\beta_1 + \beta_2} \quad 1 \leq i  \leq 8. 
$$
For simplicity,  
we will assume that $\beta_1 = \beta_2=1$, hence 
$ \balpha = \balpha'/2$.  We also assume  $(L, \nabla_l) =   (\cO_{\BP^1}(-t_n), d)$
and  set  
$$
\overline{M^{\balpha'}_4}(\bt, \blambda, -1) = \overline{M^{\balpha' \bbeta}_4}(\bt, \blambda, L), 
\quad 
\overline{M^{\balpha'}_4}( -1) = \overline{M^{\balpha' \bbeta}_4}(L). 
$$

{F}rom  Theorem \ref{thm:fund}, 
we can obtain the commutative diagram:
%%%%%%%%%%%%%%%%%%%%%%%%%%%%%%%%%%%%%
\begin{equation}\label{eq:n=4family}
\begin{CD} 
M_4^{\balpha}(-1) & \stackrel{\iota}{\hookrightarrow} & 
\overline{M_4^{\balpha'}}(-1)   \\
 @V \pi_4 VV   @VV \overline{\pi}_4 V  \\
 T_4 \times \Lambda_4  @=    T_4 \times \Lambda_4,   \\  
\end{CD}
\end{equation}
%%%%%%%%%%%%%%%%%%%%%%%%%%%%%%%%%%%%%
such that  
$\pi_4^{-1}((\bt, \blambda)) \simeq M_4^{\balpha}(\bt, \blambda, -1 )$ 
and 
$ \overline{\pi}_4^{-1}(\bt, \blambda) \simeq \overline{M_4^{\balpha'}}(\bt, \blambda, -1 )$. 
(Note that $\balpha = \balpha'/2$). 
{F}rom Theorem \ref{thm:fund}, we see that for a generic weight 
$\balpha'$, $\overline{\pi}_4$ is a 
projective morphism and $\pi_4$ is a smooth morphism of relative dimension 
$2$.

\subsection{Main Theorem (Explicit description for $n=4$ case). }

Putting  $\beta_1=\beta_2=1$, we further assume that $|\alpha'_j|\ll 1$ for $i=1,\ldots,8$. Let $\tilde{t}_1,\ldots,\tilde{t}_4\subset \BP^1\times\Lambda_4\times T_4$ be the pull-back of the universal sections on $\BP^1\times T_4$ over $T_4$.  Put $D(\tilde{\bt}):=\tilde{t}_1+\cdots+\tilde{t}_4$
and consider the projective bundle
$$
\pi:
\BP\left(\Omega^1_{\BP^1\times T_4\times\Lambda_4/T_4\times\Lambda_4}
(D(\tilde{\bt}))\oplus \cO_{\BP^1\times T_4 \times \Lambda_4}\right)
\lra \BP^1\times T_4 \times \Lambda_4.
$$
Note that since $ \Omega^1_{\BP^1}
(D(\bt)) \simeq \cO_{\BP^1}(2) $ the fiber of $ p_{23} \circ \pi $ over 
$(\bt, \blambda) \in T_4 \times \Lambda_4 $ is isomorphic to 
$$
\BP(\cO_{\BP^1}(2) \oplus \cO_{\BP^1}) \simeq \Sigma_2 
$$
where $\Sigma_2$ is the Hirzebruch surface of degree $2$. 

Let $\tilde{D}_i\subset
\BP\left(\Omega^1_{\BP^1 \times T_4 \times\Lambda_4/T_4\times\Lambda_4}
(D(\tilde{\bt}))\oplus\cO_{\BP^1\times T_4\times\Lambda_4}\right)$
be the inverse image of $\tilde{t}_i$.
Since the residue map induces an isomorphism
$$
 \Omega^1_{\BP^1\times T_4\times\Lambda_4/T_4\times\Lambda_4}
 (D(\tilde{\bt}))|_{\tilde{t}_i}
 \stackrel{\sim}\lra
 \cO_{\tilde{t}_i},
$$
we have a canonical isomorphism
$\tilde{D}_i\stackrel{\sim}\ra\BP^1\times T_4\times\Lambda_4$.
Let $\tilde{b}_i^+\subset\tilde{D}_i$ (resp.\ $\tilde{b}_i^-\subset\tilde{D}_i$)
be the inverse image of
$[\lambda_i^+:1]\subset\BP^1\times T_4\times\Lambda_4$
(resp.\ $[\lambda_i^-:1]\subset\BP^1\times T_4\times\Lambda_4$).
We denote by $B^+$ (resp.\ $B^-$) the reduced induced structure on
$\tilde{b}_1^+\cup\cdots\cup \tilde{b}_4^+$
(resp.\ $\tilde{b}_1^-\cup\cdots\cup \tilde{b}_4^-$)
and we consider the reduced induced structure on $B=B^+\cup B^-$.
Let
\[
 g:Z\ra
 \BP\left(\Omega^1_{\BP^1\times T_4\times\Lambda_4/T_4\times\Lambda_4}
 (D(\tilde{\bt}))\oplus\cO_{\BP^1\times T_4\times\Lambda_4}\right)
\]
be the blow-up along $B^+$ and $\overline{\cS}$ be the blow-up of
$Z$ along the closure of $g^{-1}(B^-\setminus (B^+\cap B^-))$.
(It is easy to  see that $\overline{\cS} \lra T_4 \times \Lambda_4$ 
is isomorphic to the family constructed by Okamoto \cite{O1}).  
Note that $Z$ is isomorphic to the blow-up of
$Z$ along $g^{-1}(B)$.

The main purpose of this section is to prove the following theorem:

%%%%%%%%%%%%%%%%%%%%%%%%%%%%%%%%%%%%%%%%%%%%%%%%%%%%%%%%%%%%%%%%%%%%%%%%%%%%%%%%%%%%%%%%%%%%%%%%%%%%%%%% Theorem %%%%%%%%%%%%%%%%%%%%%%%%%%%%%%%%%%%%%%%%%%%%%%%%%%%%%%%%%%%%%%%%%%%%%%%%%%%%%%%%%%%%%%%%%%%%%%%%%%%%%%%%%%%%%%%%%%%%%%%%%%%%%%

\begin{Theorem}\label{thm:n=4}
 Take $\balpha'=(\alpha'_i)_{1\leq i\leq 2n}$, $\bbeta=(\beta_1,\beta_2)$ 
 and $\gamma$ such that $\beta_1=\beta_2=1$, $\gamma\gg 0$,
 $|\alpha'_{i}|\ll 1$ for $1\leq i\leq 2n$,
 $\alpha'_{2i}-\alpha'_{2i-1}<\sum_{j\neq i}(\alpha'_{2j}-\alpha'_{2j-1})$
 for $1\leq i\leq n$ and that any $(\balpha',\bbeta)$-semistable
 parabolic $\phi$-connection is $(\balpha',\bbeta)$-stable.
\begin{enumerate}
\item
 There exists an isomorphism
\begin{equation}
  \overline{M^{\balpha'}_4}(\cO_{\BP^1}(-\tilde{t}_4))
  \stackrel{\sim}\lra \overline{S}
 \end{equation}
 over $T_4\times\Lambda_4$.
\item
 Let $\cY$ be the closed subscheme of
 $\overline{M^{\balpha'}_4}(\cO_{\BP^1}(-\tilde{t}_4))$
 defined by the condition $\wedge^2\phi=0$.
 Then 
 \begin{equation}
 M^{\balpha'/2}_4(\cO_{\BP^1}(-\tilde{t}_4))=
 \overline{M^{\balpha'}_4}(\cO_{\BP^1}(-\tilde{t}_4))\setminus\cY.
 \end{equation}
\item
 For each $(\bt,\blambda)\in T_4\times\Lambda_4$,
 the fiber $\cY_{(\bt,\blambda)}$ is the anti-canonical divisor of
 $\overline{M^{\balpha'}_4}(\bt,\blambda,\cO_{\BP^1}(-\tilde{t}_4))$
 and the pair
 \begin{equation}
 (\overline{M^{\balpha'}_4}(\bt,\blambda,\cO_{\BP^1}(-\tilde{t}_4)),
 \cY_{(\bt,\blambda)})
 \end{equation}
 is an Okamoto-Painlev\'e pair of type $D_4^{(1)}$.
\end{enumerate}
\end{Theorem}
%%%%%%%%%%%%%%%%%%%%%%%%%%%%%%%%%%%%%%%%%%%%%%%%%%%%%%%%%%%%%%%%%%%%%%%%%%%%%%%%%%%%%%%%%%%%%%%%%%%%%%%%%%%%%%%%%%%%%%%%%%%%%%%%%%%%%%%%%%%%%%%%%%%%%%%%%%%%%%%%%%%%%%%%%%%%%%%%%%%%%%%%%%%%%%%%%%%%%%%%%%%%%%%%%%%%%%%%%%%%%%%%%%%%%%%%%%%%%%

%%%%%%%%%%%%%%%%%%%%%%%%%%%%%%%%%%%%%%%%%%%%%%%%%%%%%%%%%%%%%%%%%%%%%%%%%%%%%%
%%%%%%%%%%%%%%%% Subsection %%%%%%%%%%%%%%%%%%%%%%%%%%%%%%%%%%%%%%%%%%%%%%%%%%
%%%%%%%%%%%%%%%%%%%%%%%%%%%%%%%%%%%%%%%%%%%%%%%%%%%%%%%%%%%%%%%%%%%%%%%%%%%%%%
\subsection{Construction of the morphism 
$\overline{M^{\balpha'}_4}(\bt,\blambda,-1)\ra
\Sigma_2$}

We assume that $(\alpha_i)$
satisfies the condition of Lemma \ref{subbundle} below.

Take any point
$(E_1,E_2,\phi,\nabla,\varphi,\{l_i\}) \in
\overline{M^{\balpha'}_4}(\bt,\blambda,-1)$.
There are unique trivial subbundles
$L^{(0)}_1\subset E_1$, $L^{(0)}_2\subset E_2$,
whose existence is confirmed by Proposition \ref{splitting} bellow.
Since the composite
\[
 \cO_{\BP^1}\cong L^{(0)}_1\hookrightarrow E_1
 \stackrel{\phi}\lra E_2\ra E_2/L^{(0)}_2\cong\cO_{\BP^1}(-1)
\]
is zero, the composite
\begin{equation}\label{hom-def-ap}
 u:L^{(0)}_1\hookrightarrow E_1\stackrel{\nabla}\lra
 E_2\otimes\Omega^1_{\BP^1}(D(\bt)) \ra
 E_2/L^{(0)}_2\otimes\Omega^1_{\BP^1}(D(\bt))\cong\cO_{\BP^1}(1)
\end{equation}
becomes a homomorphism.
By Proposition \ref{inj-def-ap} bellow,
there is a unique point $q\in\BP^1$ satisfying $u(q)=0$.
Put $L^{(-1)}_1:=E_1/L^{(0)}_1$, $L^{(-1)}_2:=E_2/L^{(0)}_2$ and
let $p_j:E_j \ra L^{(-1)}_j$ be the projection for $j=1,2$.
We define a homomorphism
$B:E_1 \ra L^{(-1)}_2\otimes\Omega^1_{\BP^1}(D(\bt))$ by
$B(a):=(p_2\otimes\mathrm{id})\nabla(a)-
d(p_2\phi(a))$ for $a\in E_1$,
where $d$ is the canonical connection on
$L^{(-1)}_2\cong\cO_{\BP^1}(-t_4)$.
Since $u_q=0$, $B_q$ induces
a homomorphism
$h_1:(L^{(-1)}_1)_q \ra
\left(L^{(-1)}_2\otimes\Omega^1_{\BP^1}(D(\bt))\right)_q$
which makes the diagram
%\[
% \begin{CD}
%  0 @>>> (L^{(0)}_1)_q @>>> (E_1)_q @>>> (L^{(-1)}_1)_q @>>> 0 \\
%  @. @. @V B_q VV  @Vh_1 VV \\
%  @. @. \left(L^{(-1)}_2\otimes\Omega_{\BP^1}(D(\bt))\right)_q 
%  @= \left(L^{(-1)}_2\otimes\Omega_{\BP^1}(D(\bt))\right)_q  
% \end{CD}
%\]
\small
\[
 \begin{array}{ccccccc}
  0 \ra & (L^{(0)}_1)_q  & \lra & (E_1)_q 
  & \lra & (L^{(-1)}_1)_q  & \ra 0 \\
  & & \searrow{\scriptstyle u_q=0} & {\scriptstyle B_q}\downarrow &
  {\scriptstyle \exists h_1}\swarrow &  &  \\
  & & & \left(L^{(-1)}_2\otimes\Omega^1_{\BP^1}(D(\bt))\right)_q 
  & & &
 \end{array}
\]
\normalsize
commute.
On the other hand, $\phi$ induces the following commutative diagram
\[
\begin{CD}
 0 @>>> L^{(0)}_1 @>>> E_1 @>>> L^{(-1)}_1 @>>> 0 \\
 @. @V \phi_1 VV @V\phi VV @V\phi_2 VV @. \\
 0 @>>> L^{(0)}_2 @>>> E_2 @>>> L^{(-1)}_2 @>>> \; 0. 
\end{CD}
\]
We put $h_2:=\phi_2(q)$.
Then $h_1,h_2$ determine a homomorphism
\begin{equation}\label{hom-def-hirz}
 \iota:(L^{(-1)}_1)_q \lra
 \left(L^{(-1)}_2\otimes\Omega^1_{\BP^1}(D(\bt))
 \oplus L^{(-1)}_2\right)_q;
 \quad a\mapsto (-h_1(a),h_2(a)).
\end{equation}
By Proposition \ref{subbundle},
$\iota$ is injective and the inclusion
$$
 \iota:(L^{(-1)}_1)_q\hookrightarrow
 (L^{(-1)}_2)_q\otimes(\Omega^1_{\BP^1}(D(\bt))\oplus\cO_{\BP^1})_q
$$
determines a point $p(E_1,E_2,\phi,\nabla,\varphi,\{l_i\})$
of ${\bf P}_*(\Omega^1_{\BP^1}(D(\bt))\oplus\cO_{\BP^1})$,
where ${\bf P}_*(\Omega^1_{\BP^1}(D(\bt))\oplus\cO_{\BP^1})$
means $\Proj S((\Omega^1_{\BP^1}(D(\bt))
\oplus\cO_{\BP^1})^{\vee})$.
So we can define a morphism
\begin{equation}\label{juyou-na-morphism}
\begin{array}{lccl} 
 p :& \overline{M^{\balpha'}_4}(\bt,\blambda,-1) & 
 \lra  & \BP_*(\Omega^1_{\BP^1}(D(\bt))\oplus\cO_{\BP^1}); \\
& & &  \\
& (E_1,E_2,\phi,\nabla,\varphi,\{l_i\}) & \mapsto & 
 p(E_1,E_2,\phi,\nabla,\varphi,\{l_i\}).
\end{array}
\end{equation}

%%%%%%%%%%%%%%%%%%%%%%%%%%%%%%% Proposition %%%%%%%%%%%%%%%%%%%%%%%%%%%%%%%%
\begin{Proposition}\label{splitting}
 For any member
 $$
 (E_1,E_2,\phi,\nabla,\varphi,\{l_i\})\in
 \overline{M^{\balpha'}_4}(\bt,\blambda,-1), 
 $$
 we have
$$
  E_1\cong E_2\cong \cO_{{\bf P}^1}\oplus\cO_{{\bf P}^1}(-1).
$$
 \end{Proposition}

\begin{proof}
Take decompositions
\begin{gather*}
 E_1=\cO_{\BP^1}(d_1)\oplus\cO_{\BP^1}(-d_1-1) \quad (d_1\geq 0) \\
 E_2=\cO_{\BP^1}(d_2)\oplus\cO_{\BP^1}(-d_2-1) \quad (d_2\geq 0).
\end{gather*}

Assume that $d_1+d_2>1$.
Then we have $\phi(\cO_{\BP^1}(d_1))\subset\cO_{\BP^1}(d_2)$.
The composite \small
$$
 \cO_{\BP^1}(d_1)\ra E_1 \stackrel{\nabla}\lra
 E_2\otimes\Omega^1_{\BP^1}(D(\bt))
 \ra \cO_{\BP^1}(-d_2-1)\otimes\Omega^1_{\BP^1}(D(\bt))\cong\cO_{\BP^1}(1-d_2)
$$
\normalsize
becomes a homomorphism and must be zero since
$H^0(\cO_{\BP^1}(1-(d_1+d_2)))=0$.
So we have
$\nabla(\cO_{\BP^1}(d_1))\subset\cO_{\BP^1}(d_2)\otimes 
\Omega^1(D(\bt))$.
Then the subbundles $(\cO_{\BP^1}(d_1),\cO_{\BP^1}(d_2))$ 
breaks
the stability of 
$(E_1,E_2,\phi,\nabla,\varphi,\{l_i\})$.

If $d_1=1$ and $d_2=0$, then $\phi(\cO_{\BP^1}(1))=0$ and the composite
\[
 f:\cO_{\BP^1}(1)\hookrightarrow E_1 \stackrel{\nabla}\lra
 E_2\otimes\Omega^1_{\BP^1}(D(\bt))
\]
becomes a homomorphism.

Put $L:=(\im f)\otimes\Omega^1(D(\bt))^{\vee}$.
Then $L$ is a vector bundle and either $L=0$ or
$L$ is a line bundle with $\deg L\geq -1$.
Then the subsheaves $(\cO_{\BP^1}(1),L)$ breaks
the stability of $(E_1,E_2,\phi,\nabla,\varphi,\{l_i\})$.

If $d_1=0$ and $d_2=1$, then the composite
$E_1\stackrel{\phi}\ra E_2\ra \cO_{\BP^1}(-2)$ must be zero and
the composite
$f:E_1\stackrel{\nabla}\ra E_2\otimes\Omega^1_{\BP^1}(D(\bt))
\ra \cO_{\BP^1}(-2)\otimes\Omega^1_{\BP^1}(D(\bt))$
becomes a homomorphism.
Put $L:=\ker f$.
Then we have either $L=E_1$ or $L$ is a line bundle such that
$\deg L\geq -1$.
Then the subbundles $(L,\cO_{\BP^1}(1))$ breaks
the stability of $(E_1,E_2,\phi,\nabla,\varphi,\{l_i\})$. 

Hence we have $d_1=d_2=0$ and
$E_1\cong E_2\cong \cO_{\BP^1}\oplus\cO_{\BP^1}(-1)$.
\end{proof}

%%%%%%%%%%%%%%%%%%%%%%%%%%%%%%%%% Lemma definition of q %%%%%%%%%%%%%%%%%%%%%%
\begin{Lemma}\label{inj-def-ap}
 For any $(E_1,E_2,\phi,\nabla,\varphi,\{l_i\})\in
 \overline{M^{\balpha'}_4}(\bt,\blambda,-1)$,
 the homomorphism $u$ defined in (\ref{hom-def-ap}) is injective.
\end{Lemma}

\begin{proof}
Assume that $u=0$.
Then the subbundles $(L^{(0)}_1,L^{(0)}_2)$
breaks the stability
of $(E_1,E_2,\phi,\nabla,\varphi,\{l_i\})$.
Thus $u\neq 0$ and $u$ is injective.
\end{proof}

%%%%%%%%%%%%%%%%%%%%%%%%% Lemma %%%%%%%%%%%%%%%%%%%%%%%%%%%%%
\begin{Lemma}\label{subbundle}
 Assume $\alpha'_{2i}-\alpha'_{2i-1}<
 \sum_{j\neq i}(\alpha'_{2j}-\alpha'_{2j-1})$
 for any $1\leq i\leq n$.
 Then the homomorphism $\iota$ defined above is injective.
\end{Lemma}

\begin{proof}
If $\phi$ is isomorphic, then
$h_2:(L^{(-1)}_1)_q\ra (L^{(-1)}_2)_q$
is isomorphic, and so $\iota$ is injective.
So we assume that $\phi$ is not isomorphic,
that is, $\wedge^2\phi=0$.

First consider the case $\rank\phi=1$.
Take decompositions
$E_1=\cO_{\BP^1}\oplus\cO_{\BP^1}(-1)$,
$E_2=\cO_{\BP^1}\oplus\cO_{\BP^1}(-1)$.
Then the homomorphism $\phi$ can be represented by a matrix
\[
 \begin{pmatrix}
  \phi_1 & \phi_3 \\
  0 & \phi_2
 \end{pmatrix}
 \quad ( \phi_1,\phi_2\in H^0(\cO_{\BP^1}),
 \phi_3\in H^0(\cO_{\BP^1}(1)) ),
\]
where the composite
$E_1\stackrel{\phi}\ra E_2\stackrel{p_2}\ra \cO_{\BP^1}(-1)$
is represented by $(0,\phi_2)$ and
$E_1\stackrel{\phi}\ra E_2\ra \cO_{\BP^1}$ by $(\phi_1,\phi_3)$.

Now assume that $p_2\circ\phi=0$.
Then $\phi_2=0$.
If moreover $\phi_1=0$, then $\phi_3\neq 0$
since $\rank\phi=1$.
Take local bases $e_1$ of $\cO_{\BP^1}\subset E_1$
and $e_2$ of $\cO_{\BP^1}(-1)\subset E_1$.
Then the condition
$\nabla(e_1)\wedge\phi(e_2)+\phi(e_1)\wedge\nabla(e_2)=0$
implies that
$\nabla(e_1)\in\cO_{\BP^1}\otimes\Omega_{\BP^1}(D(\bt))$,
which contradicts the result of Lemma \ref{inj-def-ap}.
Thus we have $\phi_1\neq 0$.
Then, by multiplying an automorphism of $E_1$ given by
\[
 \begin{pmatrix}
  c_1 & c_3 \\
  0 & c_2  
 \end{pmatrix}
 \quad
 \left( c_1,c_2\in H^0(\cO_{\BP^1}^{\times}),
 c_3\in H^0(\cO_{\BP^1}(1)) \right),
\]
the matrix representing $\phi$ changes into the form
\[
  \begin{pmatrix}
  \phi_1 & \phi_3 \\
  0 & 0
 \end{pmatrix}
 \begin{pmatrix}
  c_1 & c_3 \\
  0 & c_2
 \end{pmatrix}
 =
 \begin{pmatrix}
  c_1\phi_1 & c_3\phi_1+c_2\phi_3 \\
  0 & 0
 \end{pmatrix}.
\]
For a suitable choice of $c_1,c_2$ and $c_3$, we have
$c_1\phi_1=1$ and $c_3\phi_1+c_2\phi_3=0$.
So we may assume without loss of generality that $\phi_3=0$
and $\phi_1=1$.

The homomorphism
$B:E_1 \ra L^{(-1)}_2\otimes\Omega^1_{\BP^1}(D(\bt))
=\Omega^1_{\BP^1}(D(\bt))(-1)$ defined by
$B(a):=(p_2\otimes\mathrm{id})\nabla(a)-
d(p_2\phi(a))$ for $a\in E_1$
can be represented by a matrix
$(\omega_3,\omega_4)$ where
$\omega_3\in H^0(\Omega^1_{\BP^1}(D(\bt))(-1))$
and $\omega_4\in H^0(\Omega^1_{\BP^1}(D(\bt)))$.
Define a homomorphism
$A:E_1 \ra \Omega^1_{\BP^1}(D(\bt))$ by
$A(a):=(q_2\otimes\mathrm{id})\nabla(a)-d(q_2\phi(a))$
for $a\in E_1$, where $q_2:E_2\ra \cO_{\BP^1}$
is the projection with respect to the given decomposition
of $E_2$ and
$d$ is the trivial connection on $\cO_{\BP^1}$.
Then $A$ can be represented by a matrix
$(\omega_1,\omega_2)$, where
$\omega_1\in H^0(\Omega^1_{\BP^1}(D(\bt)))$ and
$\omega_2\in H^0(\Omega^1_{\BP^1}(D(\bt))(1))$.
Roughly speaking $\nabla$ is represented by the matrix
\[
 \begin{pmatrix}
  \omega_1 & \omega_2 \\
  \omega_3 & \omega_4
 \end{pmatrix}.
\]
Since $\phi(e_2)=0$ and $\phi(e_1)\in\cO_{\BP^1}$, the condition
$\nabla(e_1)\wedge\phi(e_2)+\phi(e_1)\wedge\nabla(e_2)=0$
implies that
$\nabla(e_2)\in\cO_{\BP^1}\otimes\Omega^1_{\BP^1}(D(\bt))$.
Thus we have $\omega_4=0$.
Take a nonzero vector
$v^{(i)}\in l_i\subset(E_1)_{t_i}$.
Then we must have
\begin{equation}\label{residue-condition}
 (\res_{t_i}\nabla)(v^{(i)})=\lambda_i\phi_{t_i}(v^{(i)}).
\end{equation}
Since $E_1=\cO_{\BP^1}\oplus\cO_{\BP^1}(-1)$,
we can write 
$v^{(i)}=\left(
\begin{array}{c} v^{(i)}_1 \\ 
v^{(i)}_2
\end{array}
\right)
$
with $v^{(i)}_1\in(\cO_{\BP^1})_{t_i}$
and $v^{(i)}_2\in(\cO_{\BP^1}(-1))_{t_i}$.
Then we have
\begin{equation*}
\begin{array}{l}
 (\res_{t_i}\nabla)
 \left(
\begin{array}{c} v^{(i)}_1 \\ 
v^{(i)}_2
\end{array}
\right)
 =\left(
 \begin{array}{l}
 \res_{t_i}(\omega_1)v^{(i)}_1
 +\res_{t_i}(\omega_2)v^{(i)}_2 \\
 \res_{t_i}(\omega_3)v^{(i)}_1
 \end{array}
 \right),  \\ 
  \\
 \quad \phi_{t_i} \left(
\begin{array}{c} v^{(i)}_1 \\ 
v^{(i)}_2
\end{array}
\right)
 = \left(
\begin{array}{c} v^{(i)}_1 \\ 
0
\end{array}
\right)
\end{array}
\end{equation*}
Thus the equality (\ref{residue-condition}) is equivalent to
the equalities
\[
 \res_{t_i}(\omega_1)v^{(i)}_1+\res_{t_i}(\omega_2)v^{(i)}_2
 =\lambda_i v^{(i)}_1, \quad
 \res_{t_i}(\omega_3)v^{(i)}_1=0.
\]
Since $u$ is injective by Lemma \ref{inj-def-ap},
$\omega_3\neq 0$.
So there is at most one point $t_i$ which satisfies
$\res_{t_i}(\omega_3)=0$, because
$\omega_3\in H^0(\Omega^1_{\BP^1}(D(\bt))(-1))\cong
H^0(\cO_{\BP^1}(1))$.
Thus, for some $i$, we have
$\res_{t_j}(\omega_3)\neq 0$ for $j\neq i$.
Then we have $v^{(j)}_1=0$ for $j\neq i$.
So we have
$l_j\subset(\cO_{\BP^1}(-1))_{t_j}$ for $j\neq i$.
Recall that the image of $\nabla|_{\cO_{\BP^1}(-1)}$
is contained in $\cO_{\BP^1}\otimes\Omega^1_{\BP^1}(D(\bt))$
because $\omega_4=0$.
Let $F_*(E_1)$ be the filtration of $E_1$
corresponding to $\{l_j\}$.
Then $(\cO_{\BP^1}(-1),\cO_{\BP^1},
\Phi|_{\cO_{\BP^1}(-1)},F_*(E_1)\cap\cO_{\BP^1}(-1))$
is a parabolic $\phi$-subconnection of
$(E_1,E_2,\Phi,F_*(E_1))$.
Since $2(\alpha'_{2i-1}+\sum_{j\neq i}\alpha'_{2j})>
\sum_{j=1}^8 \alpha'_j$ by the assumption of the lemma, we have
\begin{equation*}
\begin{array}{l}
 \mu((\cO_{\BP^1}(-1),\cO_{\BP^1},
\Phi|_{\cO_{\BP^1}(-1)},F_*(E_1)\cap\cO_{\BP^1}(-1))) \\
 \quad \geq \frac{-1-4-1-\gamma+\alpha'_{2i-1} 
 +\sum_{j\neq i}\alpha'_{2j}}{2} \\
\quad  > \frac{-2-8-2-2\gamma+\sum_{j=1}^4 (\alpha'_{2j-1}+\alpha'_{2j})}{4} 
 = \mu((E_1,E_2,\Phi,F_*(E_1))),
\end{array}
\end{equation*}
which breaks the stability of $(E_1,E_2,\Phi,F_*(E_1))$.
Therefore $p_2\circ\phi\neq 0$ and
the homomorphism $L^{(-1)}_1\ra L^{(-1)}_2$ induced by $\phi$
is an isomorphism.
Hence $h_2:(L^{(-1)}_1)_q\ra(L^{(-1)}_2)_q$
is bijective and so $\iota$ is injective.

Next consider the case $\phi=0$.
In this case, $\nabla:E_1\ra E_2\otimes\Omega^1_{\BP^1}(D(\bt))$
is a homomorphism.
If we choose a decomposition
$E_1=\cO_{\BP^1}\oplus\cO_{\BP^1}(-1)$,
$E_2=\cO_{\BP^1}\oplus\cO_{\BP^1}(-1)$,
$\nabla$ is represented by a matrix
\begin{equation*}
\begin{pmatrix}
  \omega_1 & \omega_2 \\
  \omega_3 & \omega_4
 \end{pmatrix} \quad    
\left\{
\begin{array}{l} 
\omega_1,\omega_4\in H^0(\Omega^1_{\BP^1}(D(\bt))), \\
\omega_2\in H^0(\Omega^1_{\BP^1}(D(\bt))(1)), \\ 
\omega_3\in H^0(\Omega^1_{\BP^1}(D(\bt))(-1)). \\
\end{array} \right.
\end{equation*}
Notice that $\omega_3$ corresponds to the homomorphism
$u:L^{(0)}_1\ra E_2/L^{(0)}_2\otimes\Omega_{\BP^1}(D(\bt))$
and so $\omega_3\neq 0$.
Let $q$ be the point of ${\bf P}^1$ satisfying $\omega_3(q)=0$.
Assume that $\omega_4(q)=0$.
Multiplying an automorphism of $E_1$ given by
\[
 \begin{pmatrix}
  c_1 & c_3 \\
  0 & c_2
 \end{pmatrix}
 \quad
 \left( c_1,c_2\in H^0(\cO_{\BP^1}^{\times}),
 c_3\in H^0(\cO_{\BP^1}(1)) \right),
\]
the matrix representing $\nabla$ changes into the form
\[
  \begin{pmatrix}
  \omega_1 & \omega_2 \\
  \omega_3 & \omega_4
 \end{pmatrix}
 \begin{pmatrix}
  c_1 & c_3 \\
  0 & c_2
 \end{pmatrix}
 =
 \begin{pmatrix}
  c_1\omega_1 & c_3\omega_1+c_2\omega_2 \\
  c_1\omega_3 & c_3\omega_3+c_2\omega_4
 \end{pmatrix}.
\]
For a suitable choice of $c_2,c_3$,
we have $c_3\omega_3+c_2\omega_4=0$.
So we may assume without loss of generality that $\omega_4=0$.
Take  a nonzero element $v^{(i)}$ of
$l_i\subset(E_1)_{t_i}$.
We can write $v^{(i)}=
\left(
\begin{array}{c} v^{(i)}_1 \\ 
v^{(i)}_2
\end{array}
\right)$
with $v^{(i)}_1\in(\cO_{\BP^1})_{t_i}$ and
$v^{(i)}_2\in(\cO_{\BP^1}(-1))_{t_i}$.
Then we have
\begin{align*}
 (\res_{t_i}\nabla)(v^{(i)})&=
 (\res_{t_i}\nabla)\left(
\begin{array}{c} v^{(i)}_1 \\ 
v^{(i)}_2
\end{array}
\right) \\
 &=
 \left(
\begin{array}{l} 
\res_{t_i}(\omega_1)v^{(i)}_1+\res_{t_i}(\omega_2)v^{(i)}_2 \\ 
\res_{t_i}(\omega_3)v^{(i)}_1
\end{array}
\right)
\end{align*}
Since $(\res_{t_i}\nabla)(v^{(i)})=
\lambda_i\phi_{t_i}(v^{(i)})=0$,
we have $\res_{t_i}(\omega_3)v^{(i)}_1=0$
for $i=1,\ldots,4$.
There is at most one $i$ satisfying
$\res_{t_i}(\omega_3)=0$ because
$\omega_3\in H^0(\Omega^1_{\BP^1}(D(\bt))(-1))$.
So we may assume that for some $i$,
$\omega_3(t_j)\neq 0$ for $j\neq i$.
Then we have $v^{(j)}_1=0$ for $j\neq i$ and
$l_j\subset \cO_{\BP^1}(-1)_{t_j}$ for $j\neq i$.
Since $\omega_4=0$,
$\nabla(\cO_{\BP^1}(-1))\subset
\cO_{\BP^1}\otimes\Omega^1_{\BP^1}(D(\bt))$.
If $F_*(E_1)$ is the filtration of $E_1$ corresponding to
$\{l_j\}$, then
$(\cO_{\BP^1}(-1),\cO_{\BP^1},\Phi|_{\cO_{\BP^1}(-1)},
F_*(E_1)\cap\cO_{\BP^1}(-1))$
is a parabolic $\phi$-subconnection of
$(E_1,E_2,\Phi,F_*(E_1))$ and
\begin{equation*}
\begin{array}{l}
 \mu(\cO_{\BP^1}(-1),\cO_{\BP^1},\Phi|_{\cO_{\BP^1}(-1)},
 F_*(E_1)\cap\cO_{\BP^1}(-1)) \\ 
 \quad  \geq
 \frac{-1-4-1-\gamma+\alpha'_{2i-1}+\sum_{j\neq i}\alpha'_{2j}}{2} \\
\quad  > \frac{-2-8-2-2\gamma+\sum_{j=1}^4(\alpha'_{2j-1}+\alpha'_{2j})}{4} 
 =\mu(E_1,E_2,\Phi,F_*(E_1))
\end{array}
\end{equation*}
which contradicts the stability of
$(E_1,E_2,\Phi,F_*(E_1))$.
Therefore we have $\omega_4(q)\neq 0$, which means that
$h_1$ is bijective and so $\iota$ is injective.
\end{proof}

%%%%%%%%%%%%%%%%%%%%%%%%%%%%%%%%%%%%%%%%%%%%%%%%%%%%%%%%%%%
%%%%%%%%%%%%%%%%%% Subsection Smoothness %%%%%%%%%%%%%%%%%%
\subsection{Smoothness of
$\overline{M^{\balpha'}_4}(\bt,\blambda,-1)$}

Let $\cY$ be the closed subscheme of
$\overline{M^{\balpha'}_4}(-1)$
defined by the condition $\wedge^2\phi=0$
and $Y(\bt,\blambda)$ be the fiber of $\cY$ over $(\bt,\blambda)$.

%%%%%%%%%%%%%%%%%%%%%%%%%% Proposition %%%%%%%%%%%%%%%%%%%%%%%%%%%%%%%%
\begin{Proposition}\label{inj-on-D_4^(1)}
 Under the assumption of Lemma \ref{subbundle}, the restriction
 $Y(\bt,\blambda) \stackrel{p}\lra 
 {\bf P}_*\left(\Omega^1_{\BP^1}(D(\bt))\oplus\cO_{\BP^1} \right)$
 of the morphism $p$ defined above is injective.
\end{Proposition}

\begin{proof}
Let $D_0$ be the section of
${\bf P}_*\left(\Omega^1_{\BP^1}(D(\bt))\oplus\cO_{\BP^1} \right)$
over $\BP^1$ defined by the injection
$\Omega^1_{\BP^1}(D(\bt))\hookrightarrow
\Omega^1_{\BP^1}(D(\bt))\oplus\cO_{\BP^1}$.
Take any point $(E_1,E_2,\phi,\nabla,\varphi,\{l_i\})\in Y(\bt,\blambda)$.
From the proof of Lemma \ref{subbundle}, we can see that
$p((E_1,E_2,\phi,\nabla,\varphi,\{l_i\})\in D_0$
if and only if $\phi=0$.

First assume that $\rank\phi=1$.
As in the proof of Lemma \ref{subbundle},
We take decompositions $E_1=\cO_{\BP^1}\oplus\cO_{\BP^1}(-1)$,
$E_2=\cO_{\BP^1}\oplus\cO_{\BP^1}(-1)$ and
represent $\phi$ by a matrix
\[
 \begin{pmatrix}
  \phi_1 & \phi_3 \\
  0 & \phi_2
 \end{pmatrix}
 \quad ( \phi_1,\phi_2\in H^0(\cO_{\BP^1}),
 \phi_3\in H^0(\cO_{\BP^1}(1)) ).
\]
By the proof of Lemma \ref{subbundle}, $\phi_2\neq 0$.
Multiplying a certain automorphism of $E_2$,
we may assume that $\phi_3=0$ and $\phi_2=1$.
Since $\rank\phi=1$, we have $\phi_1=0$.
Consider the homomorphism
$B:E_1\ra\cO_{\BP^1}(-1)\otimes\Omega^1_{\BP^1}(D(\bt))$
defined by $B(a)=p_2\nabla(a)-d(p_2\phi(a))$.
Let $(\omega_3,\omega_4)$
($\omega_3\in H^0(\Omega^1_{\BP^1}(D(\bt))(-1)),
\omega_4\in H^0(\Omega^1_{\BP^1}(D(\bt)))$)
be the matrix which represents $B$.
Since $\phi_1=0$, $\phi_3=0$, the composite
$E_1\stackrel{\nabla}\ra
E_2\otimes\Omega^1_{\BP^1}(\bt)\xrightarrow{q_2\otimes 1}
\cO_{\BP^1}\otimes\Omega^1_{\BP^1}(\bt)$
becomes a homomorphism, which can be represented
by a matrix $(\omega_1,\omega_2)$ with
$\omega_1\in H^0(\Omega^1_{\BP^1}(\bt)),
\omega_2\in H^0(\Omega^1_{\BP^1}(\bt)(1))$.
Roughly speaking, $\nabla$ is represented by the matrix
\[
 \begin{pmatrix}
  \omega_1 & \omega_2 \\
  \omega_3 & \omega_4
 \end{pmatrix}.
\]
We use the same notation as in the proof of Lemma \ref{subbundle}.
Then we have
$\nabla(e_1)\wedge\phi(e_2)+\phi(e_1)\wedge\nabla(e_2)=0$.
Since $\phi(e_1)=0$ and $\phi(e_2)\in\cO_{\BP^1}(-1)$,
we have
$\nabla(e_1)\in\cO_{\BP^1}(-1)\otimes\Omega^1_{\BP^1}(D(\bt))$
and so $\omega_1=0$.
Take a nonzero element $v^{(i)}$ of
$l_i\subset(E_1)_{t_i}$ and write
$v^{(i)}=\left( \begin{array}{c} v^{(i)}_1 \\ v^{(i)}_2\end{array} \right)$
where $v^{(i)}_1\in(\cO_{\BP^1})_{t_i}$
and $v^{(i)}_2\in\cO_{\BP^1}(-1)_{t_i}$.
Then we have
\begin{equation*}
\begin{array}{l}
(\res_{t_i}\nabla)(v^{(i)})=
 (\res_{t_i}\nabla)\left( \begin{array}{c} v^{(i)}_1 \\ v^{(i)}_2\end{array} \right) \\
\quad =\left( 
\begin{array}{l}
\res_{t_i}(\omega_2)v^{(i)}_2 \\ 
 \res_{t_i}(\omega_3)v^{(i)}_1+\res_{t_i}(\omega_4)v^{(i)}_2
 +\res_{t_i}\left(\frac{dz}{z-t_4}\right)v^{(i)}_2
 \end{array}
 \right), \\
 \\
 \phi_{t_i}(v^{(i)})=\phi_{t_i}\left( \begin{array}{c} v^{(i)}_1 \\ v^{(i)}_2\end{array} \right)
 =\left( \begin{array}{c} 0 \\ v^{(i)}_2\end{array} \right)
\end{array}
\end{equation*}
Since $(\res_{t_i}\nabla)(v^{(i)})=
\lambda_i\phi_{t_i}(v^{(i)})$,
we have
\begin{equation*}
\begin{array}{l}
 \res_{t_i}(\omega_2)v^{(i)}_2=0, \\
 \res_{t_i}(\omega_3)v^{(i)}_1+
 \res_{t_i}(\omega_4)v^{(i)}_2
 +\res_{t_i}\left(\frac{dz}{z-t_4}\right)v^{(i)}_2
 =\lambda_i v^{(i)}_2.
\end{array}
\end{equation*}
If $\omega_2(t_i)=0$ for any $i$, then $\omega_2=0$ because
$\omega_2\in H^0(\Omega^1_{\BP^1}(D(\bt))(1))
\cong H^0(\cO_{\BP^1}(3))$ and
there is a decomposition
$$
(E_1,E_2,\phi,\nabla,\{l_i\})=
(E_1,\cO_{\BP^1}(-1),\phi,\nabla,\{l_i\})\oplus
(0,\cO_{\BP^1},0,0,\{0\}),
$$
which contradicts the stability of
$(E_1,E_2,\phi,\nabla,\varphi,\{l_i\})$.
On the other hand, if $\omega_2(t_i)\neq 0$, then
$v^{(i)}_2=0$, $v^{(i)}_1\neq 0$ and
$\omega_3(t_i)=0$.
However, there is at most one $i$ which satisfies
$\omega_3(t_i)=0$ because 
$\omega_3\in H^0(\Omega^1_{\BP^1}(D(\bt))(-1))
\cong H^0(\cO_{\BP^1}(1))$.
Therefore there is only one $i$ which satisfies
$\omega_2(t_i)\neq 0$.
In this case, $\omega_3(t_i)=0$ and so $q=t_i$, which means
that the image $p(E_1,E_2,\phi,\nabla,\varphi,\{l_j\})$
is contained in the fiber $D_i$ of
${\bf P}_*\left(\Omega^1_{\BP^1}(D(\bt))\oplus\cO_{\BP^1} \right)$
over $t_i$.
Applying certain automorphisms of $E_1$ and $E_2$ represented by
a matrix of the form
$$
 \begin{pmatrix}
  c & 0 \\
  0 & 1
 \end{pmatrix}
 \quad ( c\in H^0(\cO^{\times}_{\BP^1}) ),
$$
we may assume that
$$
 \omega_2=\frac{\prod_{j\neq i}(z-t_j)}{\prod_{j=1}^4(z-t_j)}dz,
 \quad \omega_3=\frac{z-t_i}{\prod_{j=1}^4(z-t_j)}dz,
$$
where $z$ is a fixed inhomogeneous coordinate of $\BP^1$.
Then giving a value $\res_{t_i}(\omega_4)$ is equivalent to giving
a point $p(E_1,E_2,\phi,\nabla,\varphi,\{l_i\})$
in the fiber $D_i$.
Applying an automorphism of $E_1$
represented by a matrix of the form
\[
 \begin{pmatrix}
  1 & c \\
  0 & 1
 \end{pmatrix}
 \quad ( c\in H^0(\cO_{\BP^1}(1)) ),
\]
we may assume that $\omega_4$ is of the form
\[
 \omega_4=\frac{adz}{\prod_{j=1}^4(z-t_j)}
\]
with $a\in\C$.
$a$ is determined by the value $\res_{t_i}(\omega_4)$.
Thus the matrices representing $\phi$ and $\nabla$
are determined uniquely, up to automorphisms of $E_1$ and
$E_2$, by the point $p(E_1,E_2,\phi,\nabla,\varphi,\{l_j\})$.
Recall that $v^{(i)}_1\neq 0$, $v^{(i)}_2=0$ and
$\res_{t_j}(\omega_3)v^{(j)}_1+\res_{t_j}(\omega_4)v_2^{(j)}
+\res_{t_j}(\frac{dz}{z-t_4})v^{(j)}_2=\lambda_j v^{(j)}_2$
for $j\neq i$.
Since $\res_{t_j}(\omega_3)\neq 0$ for $j\neq i$,
every $v^{(j)}$ (including $v^{(i)}$) is uniquely determined
up to a scalar multiplication.
Thus the parabolic structure is determined by $\phi,\nabla$.
Hence $(E_1,E_2,\phi,\nabla,\varphi,\{l_j\})$
is uniquely determined by
the point $p(E_1,E_2,\phi,\nabla,\varphi,\{l_j\})$.

Next we assume that $\phi=0$. Let
$$
 \begin{pmatrix}
  \omega_1 & \omega_2 \\
  \omega_3 & \omega_4
 \end{pmatrix}, 
 \quad 
 \left\{
 \begin{array}{l}
  \omega_1,\omega_4\in H^0(\Omega^1_{\BP^1}(D(\bt))), \\
 \omega_2\in H^0(\Omega^1_{\BP^1}(D(\bt))(1)), \\
 \omega_3\in H^0(\Omega^1_{\BP^1}(D(\bt))(-1)). \\
 \end{array} \right.
$$
be a matrix representing $\nabla$.
Let $q$ be the point of ${\bf P}^1$ satisfying $\omega_3(q)=0$.
We may assume without loss of generality that $q\neq t_i$ for $i=1,2,3$.
From the proof of Lemma \ref{subbundle},
we have $\omega_4(q)\neq 0$.
Applying an automorphism of $E_1$,
we may assume
\[
 \omega_4=\frac{(z-t_1)(z-t_2)}{\prod_{j=1}^4(z-t_j)}dz,
 \quad \omega_3=\frac{z-q}{\prod_{j=1}^4(z-t_j)}dz.
\]
For a nonzero element $v^{(i)}\in l_i$,
we have $(\res_{t_i}\nabla)(v^{(i)})=
\lambda_i\phi_{t_i}(v^{(i)})=0$ for $i=1,\ldots,4$.
Thus $\det(\nabla_{t_i})=
\omega_1(t_i)\omega_4(t_i)-\omega_2(t_i)\omega_3(t_i)=0$
for $i=1,\ldots,4$.
Since $\omega_3(t_i)\neq 0$ for $i=1,2$, we have
$\omega_2(t_i)=0$ for $i=1,2$.
We write
\[
 \omega_2=\frac{(z-t_1)(z-t_2)u}{\prod_{j=1}^4(z-t_j)}dz
\]
with $u$ a polynomial in $z$ of degree less than or equal to $1$.
Applying a certain automorphism of $E_2$ of the form
\[
 \begin{pmatrix}
  c_1 & c_2 \\
  0 & 1
 \end{pmatrix}
 \quad ( c_1\in H^0(\cO_{\BP^1}^{\times}),
 c_2\in H^0(\cO_{\BP^1}(1)) ),
\]
we may assume that $u=z-t_3$.
Note that $\nabla$ is of the form
\[
 \frac{dz}{\prod_{j=1}^4(z-t_j)}
 \begin{pmatrix}
  \alpha &
  (z-t_1)(z-t_2)(z-t_3) \\
  z-q & (z-t_1)(z-t_2)
 \end{pmatrix}
 \quad ( \alpha\in H^0(\cO_{\BP^1}(2)))
\]
Since $\det(\nabla_{t_3})=0$, we have $\alpha(t_3)=0$.
The condition $\det(\nabla_{t_4})=0$ implies that
$\alpha$ is of the form $\alpha=(z-t_3)(c(z-t_4)+t_4-q)$,
where $c\in\C$.
If $c=1$, we have
$\nabla(E_1)\subset\cO_{\BP^1}(-1)\otimes\Omega^1_{\BP^1}(D(\bt))$
after applying a certain automorphism of $E_2$.
Then there is a decomposition
$(E_1,E_2,\phi,\nabla,\{l_i\})=
(E_1,\cO_{\BP^1}(-1),\phi,\nabla,\{l_i\})\oplus
(0,\cO_{\BP^1},0,0,\{0\})$,
which contradicts the stability of
$(E_1,E_2,\phi,\nabla,\varphi,\{l_i\})$.
Thus we have $c\neq 1$.
Applying a certain automorphism of $E_2$ of the form
\[
 \begin{pmatrix}
  t & (1-t)(z-t_3) \\
  0 & 1
 \end{pmatrix}
 \quad ( t\in H^0(\cO_{\BP^1}^{\times})),
\]
we may assume that $c=0$.
Since $\nabla_{t_i}\neq 0$,
$\ker(\nabla_{t_i})=l_i$ for $i=1,\ldots,4$.
Hence $(E_1,E_2,\phi,\nabla,\varphi,\{l_i\})$
is uniquely determined by $q$ and it is determined by the point
$p(E_1,E_2,\phi,\nabla,\varphi,\{l_i\})$.
\end{proof}

%%%%%%%%%%%%%%%%%%%%%%%%% Proposition %%%%%%%%%%%%%%%%%%%%%%%%%%%%%%%%%
\begin{Proposition}\label{prop:smooth}
 Under the assumption of Lemma \ref{subbundle},
 $\overline{M^{\balpha'}_4}(-1)$
 is smooth over $T_4\times\Lambda_4$.
\end{Proposition}

\begin{proof}
Let $A$ be an artinian local ring over $T_4\times\Lambda_4$
with residue field $A/m=k$ and $I$ be an ideal of $A$
such that $mI=0$.
It is sufficient to show that
\[
 \overline{M^{\balpha'}_4}(-1)(A) \lra
 \overline{M^{\balpha'}_4}(-1)(A/I)
\]
is surjective.
Take any member
\[
 (E_1,E_2,\phi,\nabla,\varphi,\{l_i\})\in
 \overline{M^{\balpha'}_4}(-1)(A/I).
\]
Note that
$E_1\cong\cO_{{\bf P}^1_{A/I}}\oplus\cO_{{\bf P}^1_{A/I}}(-1)$ and
$E_2\cong\cO_{{\bf P}^1_{A/I}}\oplus\cO_{{\bf P}^1_{A/I}}(-1)$.
Then the homomorphism $\phi:E_1\ra E_2$ can be represented
by a matrix of the form
\[
 \begin{pmatrix}
  \phi_1 & \phi_3 \\
  0      & \phi_2
 \end{pmatrix}
 \quad ( \phi_1,\phi_2\in A/I,
 \phi_3\in H^0(\cO_{{\bf P}^1_{A/I}}(1)) ).
\]
As in the proof of Proposition \ref{inj-on-D_4^(1)},
we may assume that $\phi_3\in m\otimes H^0(\cO_{\BP^1_{A/I}}(1))$.
Put
\small
\begin{align*}
 A:&=(q_2\otimes 1)\circ\nabla-d\circ q_2\circ\phi:E_1\lra
 \cO_{{\bf P}^1_{A/I}}\otimes\Omega^1_{{\bf P}^1}(D(\bt))
 \cong\cO_{{\bf P}^1_{A/I}}(2),  \\
 B:&=(p_2\otimes 1)\circ\nabla-d\circ p_2\circ\phi:E_1\lra
 \cO_{{\bf P}^1_{A/I}}(-1)\otimes\Omega^1_{{\bf P}^1}(D(\bt))
 \cong\cO_{{\bf P}^1_{A/I}}(1),
\end{align*}
\normalsize
where $q_2:E_2\ra\cO_{{\bf P}^1_{A/I}}$,
$p_2:E_2\ra\cO_{{\bf P}^1_{A/I}}(-1)$ are projections
with respect to the decomposition of $E_2$.
Let $(\omega_1,\omega_2)$ and $(\omega_3,\omega_4)$
be the matrices representing $A$ and $B$, respectively.
We can see that the condition
\small
$$
 (\varphi\otimes 1)
 (\nabla(s_1)\wedge\phi(s_2)+\phi(s_1)\wedge\nabla(s_2))
 =d(\varphi(\phi(s_1)\wedge\phi(s_2)))
 \quad (s_1,s_2\in E_1)
$$
\normalsize
is equivalent to the equality
\[
 \omega_1\phi_2-\omega_3\phi_3+\omega_4\phi_1=0.
\]
Let $(t_1,\ldots,t_4)\in\BP^1(A)\times\cdots\times\BP^1(A)$,
$(\lambda_1,\ldots,\lambda_4)\in A\times\cdots\times A$
be the data corresponding to the structure morphism
$\Spec A\ra T_4\times\Lambda_4$.
Let $v^{(i)}$ be a basis of $l_i$.
Then we can write $v^{(i)}=\left( \begin{array}{c} v^{(i)}_1 \\ v^{(i)}_2\end{array} \right)$
with $v^{(i)}_1\in\cO_{\BP^1_{A/I}}|_{t_i}$ and
$v^{(i)}_2\in\cO_{\BP^1_{A/I}}(-1)|_{t_i}$
We must find lifts
$$
\tilde{\phi}_1,\tilde{\phi}_2,\tilde{\phi}_3,\tilde{\omega}_1,
\tilde{\omega}_2,\tilde{\omega}_3,\tilde{\omega}_4,
\left( \begin{array}{c} v^{(i)}_1 \\ v^{(i)}_2\end{array} \right)_{i=1,\ldots,4}
$$
over $A$ of
$\phi_1,\phi_2,\phi_3,\omega_1,\omega_2,\omega_3,\omega_4,
\left( \begin{array}{c} v^{(i)}_1 \\ v^{(i)}_2\end{array} \right)_{i=1,\ldots,4}$ satisfying the 
following conditions:
\begin{equation*}
\left\{
\begin{array}{l}
 \tilde{\omega}_1\tilde{\phi}_2-\tilde{\omega}_3\tilde{\phi}_3
 +\tilde{\omega}_4\tilde{\phi}_1=0, \\ 
 (\res_{t_i}(\tilde{\omega}_1)-
 \lambda_i\tilde{\phi}_1)\tilde{v}^{(i)}_1
 +(\res_{t_i}(\tilde{\omega}_2)-
 \lambda_i\tilde{\phi_3}(t_i))\tilde{v}^{(i)}_2=0, \\
 \res_{t_i}(\tilde{\omega}_3)\tilde{v}^{(i)}_1+
 \left(\res_{t_i}(\tilde{\omega}_4)
 +\left(\res_{t_i}\left(\frac{dz}{z-t_4}\right)-\lambda_i\right)
 \tilde{\phi}_2 \right)\tilde{v}^{(i)}_2=0,  \\
 \text{for $i=1,\ldots,4$}.
\end{array}\right.
\end{equation*}

\normalsize
Since we have already proved the smoothness of
$M^{\balpha/2}_4(-1)$ over $T_4\times\Lambda_4$,
we may assume that $\wedge^2\phi\in mA/I$.

Assume that $\phi_1\in mA/I$ and
$\phi_2\in (A/I)^{\times}$.
Still we may assume that $\phi_3=0$.
In this case we can see from the proof of
Proposition \ref{inj-on-D_4^(1)} that
$\res_{t_i}(\omega_3)\in mA/I$ and
$\res_{t_i}(\omega_2)\in(A/I)^{\times}$ for some $i$.
Take lifts
$\tilde{\omega}_2^{(i)}\in\Omega^1_{\BP^1_A}(D(\bt))(1)_{t_i}$,
$\tilde{\omega}_4\in H^0(\Omega^1_{\BP^1_A}(D(\bt)))$,
$\tilde{\phi}_1\in A$ and $\tilde{\phi}_2\in A$
of $\omega_2(t_i)$, $\omega_4$, $\phi_1$ and $\phi_2$,
respectively.
Put
$\tilde{\omega}_1:=-\tilde{\omega}_4
\tilde{\phi}_1\tilde{\phi}_2^{-1}$.
Then we can find a lift
$\tilde{\omega}_3\in H^0(\Omega^1_{\BP^1_A}(D(\bt))(-1))$
of $\omega_3$ satisfying 
\begin{equation*}
\begin{array}{l}
 (\res_{t_i}(\tilde{\omega}_1)-\lambda_i\tilde{\phi}_1)
 \left(\res_{t_i}(\tilde{\omega}_4)
 +\left(\res_{t_i}\left(\frac{dz}{z-t_4}\right)-\lambda_i\right)
 \tilde{\phi}_2 \right) \\
\quad \quad  -
 \res_{t_i}(\tilde{\omega}_2^{(i)})\res_{t_i}(\tilde{\omega}_3)=0.
\end{array}
\end{equation*}
Let $\tilde{\omega}_2$ be the element of 
$H^0(\Omega^1_{\BP^1_A}(D(\bt))(1))$ satisfying
\begin{equation*}
\begin{array}{l}
 (\res_{t_j}(\tilde{\omega}_1)-\lambda_j\tilde{\phi}_1)
 \left(\res_{t_j}(\tilde{\omega}_4)
 +\left(\res_{t_j}\left(\frac{dz}{z-t_4}\right)-\lambda_j\right)
 \tilde{\phi}_2 \right) \\
\quad \quad  -\res_{t_j}(\tilde{\omega}_2)\res_{t_j}(\tilde{\omega}_3)=0
\end{array}
\end{equation*}
for $j\neq i$ and $\tilde{\omega}_2(t_i)=\tilde{\omega}_2^{(i)}$.
For $j=1,\ldots,4$, we can take lifts
$\tilde{v}^{(j)}_1\in\cO_{\BP^1_A}|_{t_j}$,
$\tilde{v}^{(j)}_2\in\cO_{\BP^1_A}(-1)|_{t_j}$
of $v^{(j)}_1,v^{(j)}_2$ satisfying
\[
 (\res_{t_i}(\tilde{\omega}_1)-
 \lambda_i\tilde{\phi}_1)\tilde{v}^{(i)}_1+
 \res_{t_i}(\tilde{\omega}_2)\tilde{v}^{(i)}_2=0.
\]
and
\[
 \res_{t_j}(\tilde{\omega}_3)\tilde{v}^{(j)}_1+
 \left(\res_{t_j}(\tilde{\omega}_4)
 +\left(\res_{t_j}\left(\frac{dz}{z-t_4}\right)-\lambda_j\right)
 \tilde{\phi}_2 \right)\tilde{v}^{(j)}_2=0.
\]
for $j\neq i$.
Put $\tilde{\phi}_3:=0$.
Then $\tilde{\phi}_1,\tilde{\phi}_2,\tilde{\phi}_3,
\tilde{\omega}_1,\tilde{\omega}_2,\tilde{\omega}_3,
\tilde{\omega}_4,(\tilde{v}^{(j)}_1,\tilde{v}^{(j)}_2)_{j=1}^4$
are desired lifts.

Next assume that $\phi_2\in m/I$.
In this case, we can see from the proof of
Proposition \ref{inj-on-D_4^(1)} that $\phi_1\in m/I$ and
$\phi_2\in mH^0(\cO_{\BP^1_{A/I}}(1))$.
Take a lift
$\tilde{\omega}_3\in H^0(\Omega^1_{\BP^1_A}(D(\bt))(-1))$
of $\omega_3$ and let $q\in\BP^1(A)$ be the zero point of
$\tilde{\omega}_3$.
There exists $i\in\{1,\ldots,4\}$ such that
$\res_{t_j}(\tilde{\omega}_3)\in A^{\times}$ for $j\neq i$.
Applying a certain auotomorphism of $E_1$,
we may assume that $\res_{t_i}(\omega_4)\in (A/I)^{\times}$.
Take lifts $\tilde{\omega}_4\in H^0(\Omega_{\BP^1_A}(D(\bt)))$,
$\tilde{\omega}_2^{(i)}\in\Omega_{\BP^1_A}(D(\bt)(1))_{t_i}$
and $\tilde{\phi}_2\in A$ of $\omega_4$, $\omega_2(t_i)$ and
$\phi_2$, respectively.
We can see from Lemma \ref{subbundle} that
$\tilde{\omega}_4(q)$ is a basis of
$\Omega^1_{\BP^1_A}(D(\bt))|_q$.
Then we can find an element
$\tilde{\omega}_1\in H^0(\Omega^1_{\BP^1_A}(D(\bt)))$
such that
\begin{align*}
 &\left(\res_{t_i}(\tilde{\omega}_1)\tilde{\omega}_4(q)
 +\lambda_i\tilde{\omega}_1(q)\tilde{\phi}_2\right)
 \left(\res_{t_i}(\tilde{\omega}_4)
 +\left(\res_{t_i}\left(\frac{dz}{z-t_4}\right)-\lambda_i\right)
 \tilde{\phi}_2 \right) \\
 &=\res_{t_i}(\tilde{\omega}_3)
 \res_{t_i}(\tilde{\omega}_2^{(i)})\tilde{\omega}_4(q)
 -\lambda_i\left(\res_{t_i}(\tilde{\omega}_1)
 \tilde{\phi}_2\tilde{\omega}_4(q)-
 \res_{t_i}(\tilde{\omega}_4)\tilde{\omega}_1(q)\tilde{\phi}_2\right).
\end{align*}
We can take an element $\tilde{\phi}_1$ of $A$ such that
$\tilde{\phi}_2\tilde{\omega}_1(q)
+\tilde{\phi}_1\tilde{\omega}_4(q)=0$.
Then there is an element
$\tilde{\phi}_3\in H^0(\cO_{\BP^1_A}(1))$ such that
\[
 \tilde{\omega}_1\tilde{\phi}_2-\tilde{\omega}_3\tilde{\phi}_3
 +\tilde{\omega}_4\tilde{\phi}_1=0.
\]
Let $\tilde{\omega}_2$ be the element of
$H^0(\Omega^1_{\BP^1_A}(D(\bt))(1))$ satisfying
$\tilde{\omega}_2(t_i)=\tilde{\omega}_2^{(i)}$ and
\begin{equation*}
\begin{array}{l}
 (\res_{t_j}(\tilde{\omega}_1)-\lambda_j\tilde{\phi}_1)
 \left(\res_{t_j}(\tilde{\omega}_4)
 +\left(\res_{t_j}\left(\frac{dz}{z-t_4}\right)-\lambda_j\right)
 \tilde{\phi}_2 \right) \\ 
 = \res_{t_j}(\tilde{\omega}_3)
 (\res_{t_j}(\tilde{\omega}_2)-\lambda_j\tilde{\phi}_3(t_j))
\end{array}
\end{equation*}
for $j\neq i$.
We can take lifts $\tilde{v}^{(j)}_1\in\cO_{\BP^1_A}|_{t_j}$,
$\tilde{v}^{(j)}_2\in\cO_{\BP^1_A}(-1)|_{t_j}$
of $v^{(j)}_1,v^{(j)}_2$  such that
\begin{equation*}
 \res_{t_j}(\tilde{\omega}_3)\tilde{v}^{(j)}_1+
 \left(\res_{t_j}(\tilde{\omega}_4)
 +\left(\res_{t_j}\left(\frac{dz}{z-t_4}\right)-\lambda_j\right)
 \tilde{\phi}_2 \right)\tilde{v}^{(j)}_2=0
\end{equation*}
for $j=1,\ldots,4$.
Then $\tilde{\phi}_1,\tilde{\phi}_2,\tilde{\phi}_3,
\tilde{\omega}_1,\tilde{\omega}_2,\tilde{\omega}_3,
\tilde{\omega}_4,(\tilde{v}^{(j)}_1,\tilde{v}^{(j)}_2)_{j=1}^4$
are desired lifts.
\end{proof}

%%%%%%%%%%%%%%%%%%%%%%%%%%%%%%%%%%%%%%%%%%%%%%%%%%%%%%%%%%%%%%%%%
%%%%%%%%%%%%%%%%% Subsection Theorem  %%%%%%%%%%%%%%%%%%%%%%
%%%%%%%%%%%%%%%%%%%%%%%%%%%%%%%%%%%%%%%%%%%%%%%%%%%%%%%%%%%%%%%%%%
\subsection{Proof of Theorem \ref{thm:n=4}}

We put $\lambda_i^+:=\lambda_i$ for $i=1,\ldots,4$,
$\lambda_i^-:=-\lambda_i$ for $i=1,\ldots,3$
and $\lambda_4^-:=1-\lambda_4$.
Let $D_i$ be the fiber of
${\bf P}_*(\Omega^1_{{\bf P}^1}(D(\bt))\oplus\cO_{{\bf P}^1})$
over $t_i\in{\bf P}^1$
and $b_i^+$ (resp.\ $b_i^-$) be the point of $D_i$
corresponding to $\lambda_i^+$ (resp.\ $\lambda_i^-$).
Put $Z:=\{b_1^+,\ldots,b_4^+,b_1^-,\ldots,b_4^-\}$.

%%%%%%%%%%%%%%%%%%%%%%%%%%%%% Proposition %%%%%%%%%%%%%%%%%%%%%%%%%%%%%%%%%%
\begin{Proposition}\label{prop:gen-isom}
Under the above notation,
\begin{equation}\label{iso-open}
 \overline{M^{\balpha'}_4}(\bt,\blambda,-1)\setminus p^{-1}(Z)
 \stackrel{p}\lra
 {\bf P}_*(\Omega^1_{\BP^1}(D(\bt))\oplus\cO_{\BP^1})\setminus Z
\end{equation}
is an isomorphism.
\end{Proposition}

\begin{proof}
Let $D_0$ be the section of
${\bf P}_*(\Omega^1_{\BP^1}(D(\bt))\oplus\cO_{{\bf P}^1})$
over $\BP^1$ defined by the injection
$\Omega^1_{{\bf P}^1}(D(\bt))\hookrightarrow
\Omega^1_{{\bf P}^1}(D(\bt))\oplus\cO_{{\bf P}^1}$.
First we will show that
\begin{equation}\label{iso-conn}
 \overline{M^{\balpha}_4}(\bt,\blambda,-1)
 \setminus \bigcup_{i=0}^4p^{-1}(D_i)
 \lra {\bf P}_*(\Omega^1_{{\bf P}^1}(D(\bt))\oplus
 \cO_{{\bf P}^1})\setminus \bigcup_{i=0}^4D_i
\end{equation}
is an isomorphism.
Fix a section
\[
 \tau:(\pi_2)_*(\pi_1^*\Omega^1_{{\bf P}^1}(D(\bt))|_{\Delta})
 \lra (\pi_2)_*(\pi_1^*\Omega^1_{{\bf P}^1}(D(\bt)))
\]
of the canonical homomorphism
$$
(\pi_2)_*(\pi_1^*\Omega^1_{{\bf P}^1}(D(\bt)))\lra
(\pi_2)_*(\pi_1^*\Omega^1_{{\bf P}^1}(D(\bt))|_{\Delta}),
$$
where
$$
 \pi_1:\BP^1\times(\BP^1\setminus D(\bt))\ra\BP^1,
 \quad
 \pi_2:\BP^1\times(\BP^1\setminus D(\bt))\ra
 \BP^1\setminus D(\bt)
$$
are projections and
$\Delta\subset\BP^1\times(\BP^1\setminus D(\bt))$
is the diagonal.
Take a point $s$ of
${\bf P}_*(\Omega^1_{{\bf P}^1}(D(\bt))\oplus
\cO_{{\bf P}^1})\setminus \bigcup_{i=0}^4D_i$,
which is given by $q\in{\bf P}^1$ and an injection
$(-h_1,h_2):\C\hookrightarrow
\Omega^1_{{\bf P}^1}(D(\bt))|_q\oplus\cO_{{\bf P}^1}|_q$.
We may assume that $h_2=1$.
We put
\begin{equation*}
\begin{array}{l}
 \omega_4:=\tau_q(h_1)\in H^0(\Omega^1_{{\bf P}^1}(D(\bt))), \\
 \omega_3:=\frac{z-q}{(t_4-q)\prod_{j=1}^4(z-t_j)}dz
 \in H^0(\Omega^1_{{\bf P}^1}(D(\bt))(-1)),  \\
 \end{array}
 \end{equation*}
where $z$ is a fixed inhomogeneous coordinate of $\BP^1$.
Let $\omega_2$ be the element of
$H^0(\Omega^1_{{\bf P}^1}(D(\bt))(1))$ determined by
\small
$$
 (\res_{t_i}(\omega_4)+\lambda_i)
 \left(\res_{t_i}(\omega_4)+\res_{t_i}\left(\frac{dz}{z-t_4}\right)
 -\lambda_i\right)
 +\res_{t_i}(\omega_2)\res_{t_i}(\omega_3)=0
$$
\normalsize
for $i=1,\ldots,4$.
Define a rational connection $\nabla$ on
$\cO_{{\bf P}^1}\oplus\cO_{{\bf P}^1}(-1)$ by
\[
 \nabla
 \left(
 \begin{array}{c}
 f_1 \\ 
 f_2
 \end{array}
\right):=
\left(
 \begin{array}{c}
 df_1 \\ 
 df_2
 \end{array}
\right)
+
\left(
\begin{array}{c}
-f_1\omega_4+f_2\omega_2 \\
f_1\omega_3+f_2\omega_4
\end{array}\right)
\]
for $f_1\in\cO_{\BP^1}$ and $f_2\in\cO_{\BP^1}(-1)$.
Then $s\mapsto(\cO_{\BP^1}\oplus\cO_{\BP^1}(-1),\nabla)$
determines a morphism
\[
 {\bf P}_*(\Omega^1_{\BP^1}(D(\bt))\oplus\cO_{\BP^1})
 \setminus \bigcup_{i=0}^4D_i \lra
 \overline{M^{\balpha'}_4}(\bt,\blambda,-1)
 \setminus \bigcup_{i=0}^4p^{-1}(D_i),
\]
which is just the inverse of the morphism (\ref{iso-conn}).
Then the morphism (\ref{iso-open}) is surjective,
since it is proper and dominant.
The morphism (\ref{iso-open}) is also injective
by the above argument and Proposition \ref{inj-on-D_4^(1)}.
Thus, by Zariski's Main Theorem, the morphism (\ref{iso-open})
is an isomorphism. 
\end{proof}

%%%%%%%%%%%%%%%%%%%%%%%%%%%%%%%%%%%%%%%%%%%%%%%%%%%%%%%%%%%%%%%%%%%%
%%%%%%%%%%%%%%%%%%%%%%%%%%% Proposition (-1)-curve %%%%%%%%%%%%%%%%%%%%%
\begin{Proposition}\label{prop:(-1)-curve}
If $\lambda_i^+\neq\lambda_i^-$, then
$p^{-1}(b_i^+)\cong{\bf P}^1$,
$p^{-1}(b_i^-)\cong{\bf P}^1$
and these are $(-1)$-curves.
\end{Proposition}

\begin{proof}
We can see that $p^{-1}(b_i^+)$ is just the moduli space of
$(\bt,\blambda)$-parabolic $\phi$-connections 
$(\cO_{\BP^1}\oplus\cO_{\BP^1}(-1),\cO_{\BP^1}\oplus\cO_{\BP^1}(-1),
\phi,\nabla,\varphi,\{l_j\})$
satisfying
\begin{equation*}
\begin{array}{l}
 \phi \left(
 \begin{array}{c}
 s_1 \\ 
 s_2
 \end{array} \right)=
 \left(\begin{array}{c}\phi_1s_1 \\ 
 s_2 \end{array} \right)
 \\
 \nabla\left(
 \begin{array}{c}
 s_1 \\ 
 s_2
 \end{array} \right) 
 =
 \left(
 \begin{array}{c}
 \phi_1 s_1 \\ 
 s_2
 \end{array} \right)
 + 
 \left(
 \begin{array}{c} 
  s_1\phi_1\frac{\lambda_i^+\prod_{j\neq i}(t_i-t_j)}
 {\prod_{j=1}^4(z-t_j)}dz+s_2\omega_2 \\
 s_1\frac{(z-t_i)dz}{\prod_{j=1}^4(z-t_j)}
 -s_2\frac{\lambda_i^+\prod_{j\neq i}(t_i-t_j)}
 {\prod_{j=1}^4(z-t_j)}dz
 \end{array}
  \right)
\end{array}
\end{equation*}
\normalsize
for $s_1\in\cO_{\BP^1}$ and $s_2\in\cO_{\BP^1}(-1)$,
where $\phi_1\in\C$,
$l_j=\ker(\res_{t_j}(\nabla)-\lambda_j^{+}\phi|_{t_j})$
for $j=1,\ldots,4$
and $\omega_2\in H^0(\Omega^1_{\BP^1}(D(\bt))(1))$
satisfies the condition \small
\begin{equation*}
\begin{array}{l}
 \phi_1\left( \res_{t_k} \left(
 \frac{\lambda_i^+\prod_{j\neq i}(t_i-t_j)}
 {\prod_{j=1}^4(z-t_j)}dz \right)
 -\lambda_k^+ \right) 
 \left( \res_{t_k} \left( \frac{dz}{z-t_4}-
 \frac{\lambda_i^+\prod_{j\neq i}(t_i-t_j)}
 {\prod_{j=1}^4(z-t_j)}dz \right) -\lambda_k^+ \right) \\
 -\res_{t_k}\left(\frac{(z-t_i)dz}{\prod_{j=1}^4(z-t_j)}\right)
 \res_{t_k}(\omega_2)=0.
\end{array}
\end{equation*}
\normalsize
for $k\neq i$.
Then we can define a mapping
\small
\begin{equation*}
\begin{array}{ccc}
p^{-1}(b_i^+) &  \lra &   \BP^1 \\
  &  &   \\
(\cO_{\BP^1}\oplus\cO_{\BP^1}(-1),\cO_{\BP^1}\oplus\cO_{\BP^1}(-1),
 \phi,\nabla,\varphi,\{l_j\}) & 
 \mapsto  & [\phi_1:\res_{t_i}(\omega_2)] \\
 \end{array}
 \end{equation*}
\normalsize
which is an isomorphism.

Similarly we can see that $p^{-1}(b_i^-)\cong\BP^1$.
Since $\overline{M^{\balpha'}_4}(\bt,\blambda,-1)$
and ${\bf P}_*(\Omega^1_{\BP^1}(D(\bt))\oplus\cO_{\BP^1})$
are smooth, $p^{-1}(b_i^+),p^{-1}(b_i^-)$ must be
$(-1)$-curves.
\end{proof}

%%%%%%%%%%%%%%%%%%%%%%%%%%%%%%%%%%%%%%%%%%%%%%%%%%%%%%%%%%%%%%%%%%%%%%%
%%%%%%%%%%%%%%%%%%%%%%%% Proposition (-2) curve %%%%%%%%%%%%%%%%%%%%%%%%%%
\begin{Proposition}\label{prop:(-2)-curve}
Assume that $\lambda_i^+=\lambda_i^-$.
Put
\begin{align*}
 C_1&:=\left\{(E_1,E_2,\phi,\nabla,\varphi,\{l_j\})\in p^{-1}(b_i^+)
 \left| l_i=L_1^{(0)}|_{t_i}\right\}\right. ,  \\
 C_2&:=\left\{\left.(E_1,E_2,\phi,\nabla,\varphi,\{l_j\})\in p^{-1}(b_i^+)
 \right| \res_{t_i}(\nabla)=\lambda_i\phi_{t_i} \right\}.
\end{align*}
Then $C_1\cong{\bf P}^1$, $C_2\cong{\bf P}^1$,
$C_1\cap C_2=\{\text{one point}\}$,
$C_1\cap Y(\bt,\blambda)=\{\text{one point}\}$,
$C_2\subset M^{\balpha}_4(\bt,\blambda,-1)$,
$(C_1)^2=-1$, $(C_2)^2=-2$ and
$p^{-1}(b_i^+)=C_1\cup C_2$.
\end{Proposition}

\begin{proof}
 $p^{-1}(b_i^+)$ is the moduli space of the objects
 $$
 (\cO_{\BP^1}\oplus\cO_{\BP^1}(-1),\cO_{\BP^1}\oplus\cO_{\BP^1}(-1),
 \phi,\nabla,\varphi,\{l_j\})
 $$
 satisfying
\begin{equation*}
\begin{array}{l}
 \phi
\left(
 \begin{array}{c}
s_1 \\
s_2
\end{array}
\right)=
\left(
 \begin{array}{c}
\phi s_1 \\
s_2
\end{array}
\right)
\\
 \nabla
 \left(
 \begin{array}{c}
s_1 \\
s_2 \end{array}
\right)
=
 \left(
 \begin{array}{c}
\phi ds_1 \\
ds_2
\end{array}
\right)
+
 \left(
 \begin{array}{l}
 s_1 \phi_1\frac{\lambda_i^+\prod_{j\neq i}(t_i-t_j)}
 {\prod_{j=1}^4(z-t_j)}dz+s_2\omega_2 \\
 s_1 \frac{(z-t_i)dz}{\prod_{j=1}^4(z-t_j)}
 -s_2\frac{\lambda_i^+\prod_{j\neq i}(t_i-t_j)}
 {\prod_{j=1}^4(z-t_j)}dz 
 \end{array}
 \right)
\end{array}
\end{equation*}
for $s_1\in\cO_{\BP^1}$ and $s_2\in\cO_{\BP^1}(-1)$, where
$\phi_1\in\C$, $l_k=\ker(\res_{t_k}(\nabla)-\lambda_k\phi|_{t_k})$
for $k\neq i$ and $\omega_2$ satisfies the condition 
\small
\begin{equation*}
\begin{array}{l}
 \phi_1\left( \res_{t_k} \left(
 \frac{\lambda_i^+\prod_{j\neq i}(t_i-t_j)}
 {\prod_{j=1}^4(z-t_j)}dz \right)
 -\lambda_k^+ \right) 
  \left( \res_{t_k} \left( \frac{dz}{z-t_4}-
 \frac{\lambda_i^+\prod_{j\neq i}(t_i-t_j)}
 {\prod_{j=1}^4(z-t_j)}dz \right) -\lambda_k^+ \right) \\
  \hspace{1cm} -\res_{t_k}\left(\frac{(z-t_i)dz}{\prod_{j=1}^4(z-t_j)}\right)
 \res_{t_k}(\omega_2)=0.
 \end{array}
\end{equation*}
\normalsize
for $k\neq i$.
If $v^{(i)}=\left( \begin{array}{c} v^{(i)}_1 \\ v^{(i)}_2\end{array} \right)$ is a basis of $l_i$,
$\res_{t_i}(\omega_2)v^{(i)}_2=0$.
Thus we have
\[
 p^{-1}(b_i^+)=\left(\{v^{(i)}_2=0\}\cap p^{-1}(b_i^+)\right)
 \cup\left(\{\omega_2(t_i)=0\}\cap p^{-1}(b_i^+)\right).
\]
We can see that $\{v^{(i)}_2=0\}\cap p^{-1}(b_i^+)=C_1$ and
$\{\omega_2(t_i)=0\}\cap p^{-1}(b_i^+)=C_2$.
From the proof of Proposition \ref{inj-on-D_4^(1)},
we can see that the objects of $C_2$ satisfies the condition
$\phi_1\neq 0$.
Thus we have $C_2\cap Y(\bt,\blambda)=\emptyset$.
We can also see that
$C_1\cap C_2$ consists of one point corresponding to the object of
$p^{-1}(b_i^+)$ satisfying $\omega_2(t_i)=0$, $\phi_1=1$
and $l_i=L^{(0)}_1|_{t_i}$.
$C_1\cap Y(\bt,\blambda)$ consists of one point
corresponding to the object of $C_1$ satisfying $\phi_1=0$.
We have $C_1\cong\BP^1$ by the same proof
as Proposition \ref{prop:(-1)-curve}.
$\phi,\nabla,\varphi$ and $l_k$ for $k\neq i$
are all constant on $C_2$.
So $C_2$ is just the moduli of lines
$l_i\subset\cO_{\BP^1}|_{t_i}\oplus\cO_{\BP^1}(-1)|_{t_i}$,
which is isomorphic to $\BP^1$.

Let $N_4(\bt,\blambda,-1)$
be the moduli space of rank $2$ bundles $E$ with a connection
$\nabla:E\ra E\otimes\Omega^1_{{\bf P}^1}(D(\bt))$
and a horizontal isomorphism
$\varphi:\bigwedge^2 E\stackrel{\sim}\ra\cO_{{\bf P}^1}(-x_4)$
satisfying
\begin{enumerate}
\item $\det(\res_{t_i}(\nabla)-\lambda_i\mathrm{id}_{E|_{t_i}})=0$
 for $i=1,\ldots,4$ and
\item $(E,\nabla)$ is stable in the sense of Simpson
\cite{SimpsonI}.
\end{enumerate}
Then there is a canonical morphism
\[
 M_4^{\balpha}(\bt,\blambda,-1)
 \lra N_4(\bt,\blambda,-1),
\]
which is obtained  by forgetting parabolic structure.
We can see that the image of $C_2$ in
$N_4(\bt,\blambda,-1)$
is a singular point with $A_1$-singularity.
Thus $C_2$ is a $(-2)$-curve and we can see that
$C_1$ is a $(-1)$-curve.
\end{proof}

The morphism
$p:\overline{M^{\balpha'}_4}(\bt,\blambda,\cO_{\BP^1}(-t_4))\ra
\BP(\Omega^1_{\BP^1}(D(\bt))\oplus\cO_{\BP^1})$
defined in (\ref{juyou-na-morphism}) extends to the morphism
\[
 p:\overline{M^{\balpha'}_4}
 (\cO_{\BP^1\times T_4\times\Lambda_4}(-\tilde{t}_4))\lra
 \BP\left(\Omega^1_{\BP^1\times T_4\times\Lambda_4/T_4\times\Lambda_4}
 (D(\tilde{\bt}))\oplus\cO_{\BP^1\times T_4\times\Lambda_4}\right).
\]
We can check that the inverse image $p^{-1}(B^+)$
is a Cartier divisor on
$\overline{M^{\balpha'}_4}(\bt,\blambda,\cO_{\BP^1}(-t_4))$.
Since $Z$ is a blow up of
$$
\BP\left(\Omega^1_{\BP^1\times T_4\times\Lambda_4/T_4\times\Lambda_4}
(D(\tilde{\bt}))\oplus\cO_{\BP^1\times T_4\times\Lambda_4}\right)
$$
along $B^+$, $p$ induces a morphism
\[
 f:\overline{M^{\balpha'}_4}(\bt,\blambda,\cO_{\BP^1}(-t_4))
 \lra Z.
\]
We can also check that $f^{-1}(g^{-1}(B))=p^{-1}(B)$
is a Cartier divisor on
$\overline{M^{\balpha'}_4}(\bt,\blambda,\cO_{\BP^1}(-t_4))$.
Since $\overline{S}$ is a blow up of $Z$ along
$g^{-1}(B)$, $f$ induces a morphism
\[
 f':\overline{M^{\balpha'}_4}(\bt,\blambda,\cO_{\BP^1}(-t_4))
 \lra \overline{S}.
\]
We can see by Proposition \ref{prop:gen-isom},
Proposition \ref{prop:(-1)-curve} and
Proposition \ref{prop:(-2)-curve} that
each fiber of $f'$ over $T_4\times\Lambda_4$ is an isomorphism.
Thus $f'$ is an isomorphism and
Theorem \ref{thm:n=4} (1) is proved.

Theorem \ref{thm:n=4} (2) is easy.
It is well-known that
$K_{\overline{S}_{(\bt,\blambda)}}\equiv -(2D_0+D_1+D_2+D_3+D_4)$.
So it is sufficient to prove the following proposition
in order to prove Theorem \ref{thm:n=4} (3).

%%%%%%%%%%%%%%%%%%%%%%%%%%%%% Proposition Y=-K %%%%%%%%%%%%%%%%%%%%%%%%
\begin{Proposition}
$\cY$ is a Cartier divisor on
$\overline{M^{\balpha'}_4}(-1)$ flat over $T_4\times\Lambda_4$
and the divisor $Y(\bt,\blambda)$ on 
$\overline{M^{\balpha'}_4}(\bt,\blambda,-\cO_{\BP^1}(-t_4))$
has multiplicity $2$
along $(p|_{\cY(\bt,\blambda)})^{-1}(D_0)$ and
$1$ along $(p|_{Y(\bt,\blambda)})^{-1}(D_i)$ for
$i=1,\ldots,4$.
\end{Proposition}

\begin{proof}
Let $(\cE_1,\cE_2,\tilde{\phi},\tilde{\nabla},
\tilde{\varphi},\{\tilde{l}_i\})$
be a universal family on
$\BP^1\times\overline{M^{\balpha'}_4}(-1)$.
Then $\tilde{\phi}:\cE_1\ra\cE_2$ determines a section
$f$ of $(\pi_{M^{\balpha}_4})_*(\det(\cE_1)^{-1}\otimes\det(\cE_2))$,
whose zero scheme is $\cY$.
Since $(\pi_{M^{\balpha}_4})_*(\det(\cE_1)^{-1}\otimes\det(\cE_2))$
is a line bundle on $\overline{M^{\balpha'}_4}(-1)$,
$\cY$ is a Cartier divisor on $\overline{M^{\balpha'}_4}(-1)$.
$Y(\bt,\blambda)$ is also a Cartier divisor on
$\overline{M^{\balpha'}_4}(\bt,\blambda,-1)$
and so $\cY$ is flat over $T_4\times\Lambda_4$.

Let $U_i$ be the open subscheme of $Y(\bt,\blambda)$
whose underlying space is
$(p|_{Y(\bt,\blambda)})^{-1}(D_i\setminus(D_0\cap D_i))$.
Then $U_i$ is just the moduli space of the objects
$(\cO_{\BP^1}\oplus\cO_{\BP^1}(-1),
\cO_{\BP^1}\oplus\cO_{\BP^1}(-1),\phi,\nabla,\varphi,\{l_j\})$
satisfying
\begin{align*}
 \phi
 \left(
 \begin{array}{c}
 f_1 \\f_2
 \end{array}
 \right)
 &=
  \left(
 \begin{array}{c}
 0 \\f_2
 \end{array}
 \right), \\
 \nabla
 \left(
 \begin{array}{c} 
 f_1 \\
 f_2
 \end{array}
\right)&=
\left(
\begin{array}{c}
0 \\
df_2
\end{array}
\right)
 +\left(
 \begin{array}{l} 
 f_2\frac{\prod_{j\neq i}(z-t_j)}{\prod_{j=1}^4(z-t_j)}dz \\
 f_1\frac{(z-t_i)dz}{\prod_{j=1}^4(z-t_j)}
 +f_2\frac{adz}{\prod_{j=1}^4(z-t_j)}
 \end{array}
 \right)
\end{align*}
for $f_1\in\cO_{\BP^1}$ and $f_2\in\cO_{\BP^1}(-1)$,
where $a\in\C$ and
$l_j=\ker(\res_{t_j}(\nabla)-\lambda_j\phi_{t_j})$
for $j=1,\ldots,4$.
Thus $U_i\cong\mathbf{A}^1$ and $U_i$ is reduced.

Let $U_0$ be the open subscheme of $Y(\bt,\blambda)$
such that $p(U_0)=D_0\setminus\bigcup_{j=1}^4 D_j$ as sets.
$U_0$ is the moduli space of the objects
$(\cO_{\BP^1}\oplus\cO_{\BP^1}(-1), \linebreak
\cO_{\BP^1}\oplus\cO_{\BP^1}(-1),\phi,\nabla,\varphi,\{l_j\})$
satisfying
\begin{align*}
 \phi
 \left(
 \begin{array}{c}
 f_1
\\
 f_2
 \end{array}
 \right)&=
 \left(
 \begin{array}{l}
 f_1\phi_1+f_2\phi_3 \\
 f_2\phi_2
 \end{array}
 \right)\\
 \nabla
 \left(
 \begin{array}{c}
 f_1 
 \\
 f_2
 \end{array}
 \right)&=
 \left(
 \begin{array}{l}
 \phi_1df_1+\phi_3df_2 
 \\
 \phi_2df_2
 \end{array}
 \right)
 +
 \left(
 \begin{array}{l}
 \omega_1f_1
 \\
 \omega_3f_1+\omega_4f_2
 \end{array}
 \right)
\end{align*}
for $f_1\in\cO_{\BP^1}$ and $f_2\in\cO_{\BP^1}(-1)$
with the conditions $\phi_1\phi_2=0$ and
$\omega_1\phi_2-\omega_3\phi_3+\omega_4\phi_1=0$, where
$q\in\BP^1\setminus\{t_1,\ldots,t_4\}$,
$l_j=\ker(\res_{t_j}(\nabla)-\lambda_j\phi_{t_j})$
for $j=1,\ldots,4$ and
\small
\begin{align*}
 \omega_1&=\frac{\prod_{k=3}^4(z-t_k+(t_k-t_1)(t_k-t_2)
 \lambda_k\phi_1)}{\prod_{j=1}^4(z-t_j)}dz, \
 \omega_3=\frac{(z-q)dz}{(t_4-q)\prod_{j=1}^4(z-t_j)} \\
 \omega_4&=\frac{\prod_{k=1}^2(z-t_k+(t_k-t_3)(t_k-t_4)\lambda_k\phi_2)}
 {\prod_{j=1}^4(z-t_j)}dz.
\end{align*}
\normalsize
$\phi_2$ and $\phi_3$ are determined by
$\phi_1$ and the conditions
\[
 \omega_1(q)\phi_2+\omega_4(q)\phi_1=0,\quad
 \omega_3(t_j)\phi_3(t_j)=\omega_1(t_j)\phi_2+\omega_4(t_j)\phi_1
 \quad (j=1,2)
\]
and $\phi_2$ must satisfy the condition $\phi_1^2=0$.
Thus $U_0\cong\BP^1\setminus\{t_1,\ldots,t_4\}
\times\Spec\C[\phi_1]/(\phi_1^2)$ and
$Y(\bt,\blambda)$ has multiplicity $2$ along
$(p|_{Y(\bt,\blambda)})^{-1}(D_0)$.
\end{proof}

%%%%%%%%%%%%%%%%%%%%%%%%%%%%%%%%%%%%%%%%%%%%%%%%%%%
\section{Moduli of stable parabolic connections in general case}
\label{sec:mod-conn}
%%%%%%%%%%%%%%%%%%%%%%%%%%%%%%%%%%%%%%%%%%%%%%%%%%%

In this section, we will formulate the general  moduli theory of  
$\balpha$-stable parabolic 
connections over a curve and state 
the existence theorem of the coarse 
moduli scheme due to Inaba \cite{Ina}.  
We fix   integers $g, d, r, n $ with $g \geq 0,  r > 0, n> 0$ and 
let $(C, \bt) = (C, t_1, \ldots, t_n)$ be an 
$n$-pointed 
smooth projective curve of genus $g$,  
which consists of a smooth projective curve $C$ 
and a set of $n$-distinct points $\bt= \{ t_i \}_{1 \leq i \leq n}$ on $C$.  
We denote by  $D(\bt) = t_1 + \cdots + t_n$ the divisor associated to 
$\bt $. Define 
the set of exponents as 
\begin{equation}\label{eq:setofexp}
\Lambda_{r}^{n}(d) := \left\{ \blambda = (\lambda^{(i)}_j)^{1\leq i \leq n}_{0 \leq j \leq r-1} \in \C^{nr} \quad  \left| 
 \quad d + \sum_{1\leq i \leq n, \ 0 \leq j \leq r-1} \lambda^{(i)}_j = 0 \right\} \right..
\end{equation}

\begin{Definition}\label{def:parabolic} 
{\rm A {\em $(\bt, \blambda)$-parabolic connection of rank $r$} on 
$C$ is a collection of data $(E, \nabla, \{ l_{\ast}^{(i)} \}_{1 \leq i \leq n} )$ consisting of:
\begin{enumerate}
\item a vector bundle $E$ of rank $r$ on $C$, 
\item a logarithmic connection $\nabla:E \lra E \otimes \Omega^1_{C}(D(\bt))$, 
\item and a filtration $l_{\ast}^{(i)}: E_{|t_i} = l_0^{(i)} \supset l_{1}^{(i)} \supset \cdots 
\supset l_{r-1}^{(i)} \supset l_{r}^{(i)} =0 $ for each $i, 1 \leq i \leq n$
such that $\dim (l_j^{(i)}/l_{j+1}^{(i)}) = 1$ and $(\res_{t_i}(\nabla) - \lambda^{(i)}_j)(l_j^{(i)}) 
\subset l_{j+1}^{(i)}$ for $ j = 0, 1, \cdots, r-1$. 
\end{enumerate}
We set $\deg E = \deg \left( \wedge^r E \right)  $ as usual. 
} 

\end{Definition}

Take a sequence of rational numbers 
$
\balpha = ( \alpha^{(i)}_j )^{1\leq i \leq n}_{1 \leq j \leq r}
$ 
such that 
\begin{equation} \label{eq:weights}
 0<\alpha^{(i)}_1<\alpha^{(i)}_2<\cdots<\alpha^{(i)}_r<1
\end{equation}
for $i= 1 , \ldots , n$ and 
$\alpha^{(i)}_j\neq\alpha^{(i')}_{j'}$ for $(i,j)\neq(i',j')$.
We choose $\balpha=(\alpha^{(i)}_j)$ sufficiently generic.
Let $(E,\nabla,\{l^{(i)}_*\}_{1\leq i\leq n})$ be a 
$(\bt, \blambda)$-parabolic connection, and 
$F \subset E$ a nonzero subbundle 
satisfying $\nabla(F)\subset F\otimes\Omega_C^1(D(\bt))
$. 
We define integers $  \len (F)^{(i)}_j$ by  
\begin{equation}\label{eq:length}
\len (F)^{(i)}_j = \dim (F|_{t_i}\cap l^{(i)}_{j-1})/(F|_{t_i}\cap l^{(i)}_j).
\end{equation}
Note that 
$\len(E)^{(i)}_j = \dim (l_{j-1}^{(i)}/l_{j}^{(i)}) = 1$ 
for $1 \leq j \leq r$.

\begin{Definition}{\rm
A parabolic connection 
$(E,\nabla,\{l^{(i)}_*\}_{1\leq i\leq n})$
is $\balpha$-stable if
for any proper nonzero subbundle $F \subsetneqq E$
satisfying $\nabla(F)\subset F\otimes\Omega_C^1(D(\bt))$,
the inequality
\begin{equation}
 \frac{\deg F + \sum_{i=1}^m \sum_{j=1}^r
 \alpha^{(i)}_j \len(F)^{(i)}_j}{\rank F} 
 < \frac{\deg E + \sum_{i=1}^n \sum_{j=1}^r \alpha^{(i)}_j 
 \len(E)^{(i)}_j}{\rank E}
\end{equation}
holds.}
\end{Definition}

For a fixed  $(C, \bt)$ and $\blambda$, let 
us define the coarse moduli space by  
\begin{equation}\label{eq:connection}
\begin{array}{l}
\cM^{\balpha}_{((C, \bt), \blambda)}(r,n,d) = \\
\{ ( E, \nabla, \{ l_{\ast}^{(i)} \}_{1 \leq i \leq n} ) \ | \ 
\begin{array}{l}
\mbox{an $\balpha$-stable $(\bt, \blambda)$-parabolic connection } \\
\mbox{of rank $r$ and degree $d$ over $C$}
\end{array} \}/\simeq. 
\end{array}
\end{equation}

Varying $(C, \bt)$ and $\blambda$, we can also  consider the moduli space  
in relative setting. 
Let $\cM_{g, n}$ be the coarse moduli space of 
$n$-pointed curves of genus $g$. Here we assume that 
every  point of  $\cM_{g, n}$ corresponds to 
an $n$-pointed smooth curve $(C, \bt)$ such that $\bt =(t_1, \ldots, t_n)$ 
is a set of $n$-distinct points on $C$.  
We consider a finite covering $\cM'_{g,n} \ra \cM_{g,n}$ 
where  $\cM'_{g,n}$ is the coarse moduli space  
of $n$-pointed curves of genus $g$ 
with a suitable level structure so that  
there exists  the  universal family  
$(\cC, \tilde{\bt}) = (\cC, \tilde{t}_1, \ldots,\tilde{t}_n)$
of $n$-pointed curves (with a level structure).  
From now on, for simplicity, we set 
\begin{equation}\label{eq:n-pointed}
T = \cM'_{g,n}
\end{equation}
and let  
\begin{equation}\label{eq:n-pointed-fam}
(\cC, \tilde{\bt}) \lra T = \cM'_{g, n} 
\end{equation}
be the universal family.

We can show the existence theorem of moduli space as a smooth quasi-projective 
algebraic scheme (cf. \cite{IIS-1}, \cite{Ina}).

\begin{Theorem}\label{thm:existence} $($Cf. \cite{IIS-1},  \cite{Ina}$)$. 
Assume that $r, n, d$ are positive integers. 
There exists a relative moduli scheme
\begin{equation}\label{eq:conn-fam}
\varphi_{r,n,d}: \cM^{\balpha}_{(\cC, \tilde{\bt}) /T }(r,n,d) \lra T \times \Lambda^{(n)}_r(d)
\end{equation}
of $\balpha$-stable parabolic connections of rank $r$ and degree $d$,
which is smooth and quasi-projective over $T \times \Lambda^{(n)}_r(d)$.  
Moreover the fiber 
$\cM^{\balpha}_{((C, \bt), \blambda)}(r,n,d)$ of $\varphi_{r,n,d}$ 
over $((C, \bt),\blambda) \in  T   \times \Lambda^{(n)}_r(d)$ 
is the moduli space of $\balpha$-stable 
$(\bt,\blambda)$-parabolic connections over $C$,   
which is a {\em smooth algebraic scheme} and  
\begin{equation}\label{eq:dim-moduli}
\dim  \cM^{\balpha}_{((C, \bt), \blambda)}(r,n,d) = 2r^2(g-1)+nr(r-1)+2.
\end{equation}

\end{Theorem}

\begin{Remark}
{\rm 
\quad 
\vspace{0.2cm}
\begin{enumerate}
\item When $C = \BP^1$ and $r = 2$, Theorem \ref{thm:existence} is  proved in \cite{IIS-1}.  
\item 
Inaba \cite{Ina} showed that the moduli space 
$\cM^{\balpha}_{((C, \bt), \blambda)}(r,n,d)$ is irreducible in the following cases:
\begin{enumerate}
\item $g \geq 2, n \geq 1$, 
\item $g =1, n \geq 2$, 
\item $ g= 0, r \geq 2, r n - 2r - 2 >0 $
\end{enumerate}

\end{enumerate}
}
\end{Remark}

%%%%%%%%%%%%%%%%%% SS {The moduli space of representations}%%%%%%%%%
\subsection{The moduli space of representations}
\label{ss:moduli-of-rep}
%%%%%%%%%%%%%%%%%%%%%%%%%%%%%%%%%%%%%%%%%%%%%%%%%%%%%%%%%%%%%%%%%%%

For each $n$-pointed curve $(C, \bt) = (C, t_1, \cdots, t_n) \in T = \cM'_{g,n} $
$(g \geq 0, n \geq 1)$, set $D(\bt) = t_1 + \cdots + t_n$.  
By abuse of notation, 
we denote by  
$\pi_1(C \setminus D(\bt)*)$ the fundamental group  of 
$ C \setminus \{ t_1, \cdots, t_n \}$.   
The set  
\begin{equation}
 \Hom(\pi_1(C \setminus D(\bt), *),
 GL_r(\C))
\end{equation}
of $GL_{r}(\C)$-representations of  
$\pi_1(C \setminus D(\bt),*)$ 
is an affine variety, and  $GL_r(\C)$ naturally acts on this space by the 
adjoint action.   

We define the moduli space by 
\begin{equation}\label{eq:rep-space}
\cRP_{(C, \bt)}^r = \Hom(\pi_1(C \setminus D(\bt),*),
 GL_r(\C))//Ad(GL_r(\C)). 
\end{equation}
Here the quotient $//$ means the categorical quotient (\cite{Mum:GIT}). 
More precisely, 
it is known that $\pi_1(C \setminus D(\bt),*)$ is generated by $(2 g + n)$-elements 
$\alpha_1, \ldots, \alpha_{g}, \beta_1, \ldots, 
\beta_{g}, \gamma_1, \ldots, \gamma_n$ with one 
relation 
$$ 
\prod_{i=1}^{g} [\alpha_i, \beta_i] \gamma_1 \cdots \gamma_n = 1.
$$
Therefore  if we denote by $R$ the ring of invariants of the 
simultaneous adjoint action of $GL_r(\C)$ on the 
coordinate ring of $GL_r(\C)^{2g + n-1}$,  then we have 
an isomorphism 
\begin{equation}
\cRP_{(C, \bt)}^r  \simeq \Spec(R). 
\end{equation}
Hence the moduli space  $\cRP_{(C, \bt)}^r$  becomes 
an affine algebraic scheme.  Furthermore,  each closed point
of $\cRP_{(C, \bt)}^r$  corresponds to a Jordan equivalence class of 
a representation (cf. [Section 4, \cite{IIS-1}]).

Let us set 
\begin{equation}\label{eq:a-param}
 \cA^{(n)}_r:=\left\{
 \ba=(a^{(i)}_j)^{1\leq i\leq n}_{0\leq j\leq r-1} \in \C^{nr}
 \left| a^{(1)}_0 a^{(2)}_0 \cdots a^{(n)}_0=(-1)^{rn}
 \right\}\right..  
\end{equation}
For each $\ba = (a^{(i)}_j) \in \cA^{(n)}_r$ and $ i, 1 \leq i \leq n $, 
we set $\ba^{(i)} = (a_{0}^{(i)}, \cdots, a_{r-1}^{(i)}) $ and
define 
\begin{equation}\label{eq:char}
\chi_{\ba^{(i)}}(s) = s^r + a^{(i)}_{r-1} s^{r-1} + \cdots + a^{(i)}_{0}. 
\end{equation} 

Moreover we define a morphism
\begin{equation}\label{eq:character}
\phi_{(C, \bt)}^r:\cRP_{(C, \bt)}^r \lra
 \cA^{(n)}_r
\end{equation}
by the relation
\begin{equation}\label{eq:char-poly}
\det(sI_r-\rho(\gamma_i))=\chi_{\ba^{(i)}}(s)
\end{equation}
where $[\rho] \in \cRP_{(C, \bt)}^r$  and  $\gamma_i$ is a counterclockwise 
loop around $t_i$.

For $\ba =  (a^{(i)}_j) \in \cA^{(n)}_r$, we denote by $\cRP_{(C, \bt), \ba}^r$ 
the fiber of $\phi_{(C, \bt)}^{r}$ over $\ba$, that is, 
\begin{equation}\label{eq:fiber}
\cRP_{(C, \bt),\ba}^r =  \{ [\rho] \in \cRP_{(C, \bt)}^r \vert  
\det(sI_r-\rho(\gamma_i))=\chi_{\ba^{(i)}}(s), 1 \leq i \leq n  \}. 
\end{equation}

\normalsize
For any  covering $T' \ra T $, 
we can define a relative moduli space
$\cRP^r_{n,T'}= \coprod_{(C, \bt) \in T' } \cRP^r_{(C, \bt)}$
of representations with the natural morphism 
\begin{equation}\label{eq:natural}
\cRP^r_{n, T'} \lra T'.  
\end{equation}

As in Section 4, \cite{IIS-1}, there exists a finite covering 
$T' \lra T$ with the morphism 
\begin{equation}\label{eq:phi-r-n}
\phi^{r}_n: \cRP^r_{n, T'}    \lra  T' \times \cA^{(n)}_r, 
\end{equation}
such that 
$$
(\phi^{r}_n)^{-1}((C, \bt), \ba) = \cRP_{(C, \bt),\ba}^r.
$$

%%%%%%%%%%%%%% SS \subsection{RH correspondence} $$$$$$$$$$$
\section{The Riemann-Hilbert correspondence}
\label{ss:RH}
%%%%%%%%%%%%%%%%%%%%%%%%%%%%%%%%%%%%%%%%%%%%%%%%%%%%%%%%%%%%

Next we define the 
Riemann-Hilbert correspondence from the moduli space of $\balpha$-stable 
parabolic connections to the moduli space of the representations. 

Let us fix positive integers $r, d$, 
$\balpha = ( \alpha_j^{(i)}) $ as in (\ref{eq:weights}), 
and  $(C, \bt) \in  T' = \cM'_{g,n}$. 
For simplicity, we set   $ \cM^{\balpha}_{((C, \bt), \blambda)}= 
\cM^{\balpha}_{((C, \bt), \blambda)}(r,n,d)$ (cf. (\ref{eq:connection})). 

We define a morphism
\begin{equation}\label{eq:rh-exponent}
 rh:\Lambda^{(n)}_r(d) \lra \cA^{(n)}_r, \quad  
 rh(\blambda) = \ba
\end{equation}
by the relation 
\begin{equation}\label{eq:rel-exponents}
 \prod_{j=0}^{r-1}(s-\exp(-2\pi\sqrt{-1}\lambda^{(i)}_j))=
 s^r+a^{(i)}_{r-1}s^{r-1}+\cdots+a^{(i)}_0.
\end{equation}

For each member
$(E,\nabla,\{l^{(i)}_j\})  \in
\cM^{\balpha}_{(C,\bt),\blambda}$, 
the solution subsheaf of $E^{an}$ 
\begin{equation}\label{eq:local-sys}
\ker(\nabla^{an}|_{C \setminus D(\bt)}) \subset E^{an}
\end{equation}
becomes a local system on
$C \setminus D(\bt)$
and corresponds to a representation 
\begin{equation}
\rho:\pi_1(C \setminus \{\bt \}, *) \lra GL_{r}(\C). 
\end{equation}
Since the eigenvalues of the residue matrix of $\nabla^{an}
$ at $t_i$ 
are $\lambda_j^{(i)}$, $0 \leq j \leq r-1$, considering 
the local fundamental solutions of $\nabla^{an} = 0$ near $t_i$, 
the monodromy matrix of  $ \rho(\gamma_i) $  has eigenvalues 
$ \exp(-2\pi\sqrt{-1}\lambda^{(i)}_j)$, $ 0 \leq j \leq r-1$. 
Hence under the relation (\ref{eq:rel-exponents}),  
or $\ba = rh(\blambda)$,  we 
can define a morphism
\begin{equation}
\RH_{(C, \bt), \blambda}: \cM^{\balpha}_{((C, \bt),\blambda)}\lra
\cRP_{(C,\bt), \ba}^r. 
\end{equation}

Replacing  $T=\cM'_{g, n}$ by a certain finite \'etale covering 
$u:T' \lra T$ and  varying $((C, \bt), \blambda) \in T' \times 
\Lambda_r^{(n)}(d)$
we can define a morphism
\begin{equation}\label{eq:global-RH}
 \RH: \cM^{\balpha}_{(\cC, \bt)/T'}(r,n,d) \lra \cRP^r_{n,T'} 
\end{equation}
which makes the diagram
\begin{equation}\label{eq:grh}
\begin{CD}
\cM^{\balpha}_{(\cC, \tilde{\bt}) /T'}(r,n,d) @> \RH >>
  \cRP^r_{n,T'} \\
  @V \varphi_{r,n,d} VV @VV \phi^{r}_n V \\
  T' \times \Lambda^{(n)}_r(d) @> Id \times rh >> T' \times\cA^{(n)}_r
 \end{CD}
\end{equation}
commute. The following result is proved in \cite{Ina}.  

\begin{Theorem}\label{thm:rh} $(${\rm [Theorem 2.2,  \cite{Ina}]} $)$. 
 Assume that $\balpha$ is so generic that $\balpha$-stable $\Leftrightarrow$ $\balpha$-semistable.  Moreover we assume that $r \geq 2, r n - 2 r - 2 >0 $ if $g =0$,  $ n \geq 2 $ if $ g = 1$ and $ n \geq 1 $ if $g \geq 2$.  
 Then the  morphism
\begin{equation}\label{eq:rh1}
 \RH: \cM^{\balpha}_{(\cC, \tilde{\bt})/T'}(r,n,d) \lra
 \cRP^r_{n,T'} \times_{\cA^{(n)}_r} \Lambda^{(n)}_r
\end{equation}
induced by (\ref{eq:global-RH}) is a {\em proper surjective
bimeromorphic analytic} morphism.
In particular, for each $((C, \bt),\blambda)\in T' \times\Lambda^{(n)}_r(d)$, the restricted morphism 
\begin{equation}
 \RH_{((C, \bt), \blambda)}: \cM^{\balpha}_{((C, \bt),\blambda)}(r, n, d)
 \lra
 \cRP^r_{(C,\bt), {\ba}}
\end{equation}
gives an analytic resolution of singularities of
$\cRP^r_{(C,\bt), \ba}$ where $\ba=rh(\blambda)$.
\end{Theorem}

\begin{Remark}\label{ref:sing}
{\rm 
Take $\blambda \in \Lambda^{(n)}_r$ such that $rh(\blambda) = \ba$.  
A representation $\rho$ such that $[\rho] \in \cRP^r_{(C, \bt), \ba}$ 
is said to be {\em resonant} if 
\begin{equation}\label{eq:resonant}
\text{$\dim(\ker(\rho(\gamma_i)-\exp(-2\pi\sqrt{-1}\lambda^{(i)}_j)))
  \geq 2$ for some $i,j$}.
\end{equation}
The singular locus of $\cRP^r_{(C, \bt), \ba}$ 
is given by the set 
\begin{equation}\label{eq:rep-sing}
\left(\cRP^r_{(C, \bt), \ba} \right)^{sing} :=
 \left\{ [\rho]\in \cRP^r_{(C, \bt), \ba} \left|
 \begin{array}{l}
  \text{$\rho$ is reducible or} \\
  \text{ resonant}
 \end{array}
 \right\}\right..
\end{equation}
Moreover we denote  the smooth part of $\cRP^r_{(C, \bt), \ba}$  by 
\begin{equation}\label{eq:smooth}
  \left(\cRP^{r}_{(C,\bt), \ba}\right)^{\sharp} =  \cRP^r_{(C, \bt), \ba} 
  \setminus \left(\cRP^r_{(C, \bt), \ba} \right)^{sing}. 
\end{equation}  
Theorem \ref{thm:rh} implies that the restriction
\begin{equation}
\RH_{((C, \bt), \blambda)| 
\left(\cM^{\balpha}_{(C,\bt),\blambda}\right)^{\sharp}}:
\left(\cM^{\balpha}_{(C, \bt),\blambda}\right)^{\sharp} \stackrel{\simeq}{\lra} 
\left( \cRP^r_{(C,\bt), {\ba}}\right)^{\sharp}
\end{equation}
is an analytic isomorphism, where
$$
\left( \cM^{\balpha}_{(C, \bt),\blambda} \right)^{\sharp}=
\RH_{((C, \bt), \blambda)}^{-1}(\left(\cRP^r_{(C, \bt), {\ba}} \right)^{\sharp}).
$$
}\end{Remark}

%%%%%%%%%%%%%%%%%%%%%%%%%%%%%%%%%%%%%%
\section{Isomonodromic flows and Differential systems of Painlev\'e type}
\label{sec:isomon}
%%%%%%%%%%%%%%%%%%%%%%%%%%%%%%%%%%%%%

Consider the family of the moduli spaces 
of $\balpha$-stable parabolic connections 
\begin{equation}\label{eq:mod-fam-1}
\varphi_{r,n,d}: \cM^{\balpha}_{(\cC, \bt)/T}(r,d, n) \lra T \times \Lambda^{(n)}_r(d)
\end{equation}
where 
$
T = \cM'_{g,n}
$
as in (\ref{eq:n-pointed}). 

Fix $((C_0, \bt_0), \blambda_0) \in T \times \Lambda^{(n)}_r(d)$ and 
take an $\balpha$-stable parabolic connection $\x = (E, \nabla, \{ l^{(i)}_{*} 
\}_{1 \leq i \leq n}) \in 
\cM^{\balpha}_{((\cC_0, \bt_0), \blambda_0)}(r,d, n)$. 
Let $\Delta =\{ t \in \C | |t| < 1 \}$ be the unit disc and let 
$h:\Delta \lra T $ be a holomorphic embedding such that  $h(0) = (C_0, \bt_0)$. Then pulling back the universal family, we obtain the 
family of $n$-pointed curves $f:(\cC, \bt) \lra \Delta $ with the central fiber 
$ f^{-1}(0) = (C_0, \bt_0) $.  
An $\balpha$-stable parabolic connection  
$({\mathcal E}, {\mathbf \nabla}, { l} )$  
on the family of $n$-pointed curves 
$ (\cC, \bt)$ over $\Delta$ is called 
a {\em $($$1$-parameter$)$ deformation}  of  
$(E, \nabla, \{ l^{(i)}_{*} \}_{1 \leq i \leq n}) $ 
if we have an isomorphism  
$({\mathcal E}, {\mathbf \nabla}, { l} )_{|(C_0, \bt_0)} 
\simeq (E, \nabla, \{ l^{(i)}_{*} 
\}_{1 \leq i \leq n}) $.  Restricting the $\balpha$-stable parabolic connection $({\mathcal E}, {\mathbf \nabla}, { l} )$ to each fiber $({\mathcal C}_t, \bt_t)$, we have a family of $\balpha$-stable parabolic connections $({\mathcal E}_t, {\mathbf \nabla}_t, { l}_t )$ over $({\mathcal C}_t, \bt_t)$ which are automatically flat in the direction of each fiber.  
If the connection $\nabla$ 
on ${\mathcal E}$ is  flat on the total space ${\mathcal C}$, 
which means that the curvature $2$-form of $\nabla$ vanishes 
over the total space ${\mathcal C}$, the  associated 
representations $\rho_t: \pi_1({\mathcal C}_t \setminus \{ \bt_t \}, *) 
\lra GL_r(\C)$  is constant with respect to  $t \in \Delta$.  Moreover 
the converse is also true.  
Therefore such a deformation  
$({\mathcal E}, {\mathbf \nabla}, { l} )$ over  
${\mathcal C} \lra \Delta $  is called 
{\em an isomonodromic 
deformation} of a $\balpha$-stable parabolic connection.  
Under an isomonodromic deformation,  
local exponents $\blambda_t$ of the connection 
$({\mathcal E}_t, {\mathbf \nabla}_t, {l}_t )$
 are also constant, so we have $\blambda_t = \blambda_0$.  
Therefore an isomonodromic deformation determines a holomorphic map  
$\tilde{h}: \Delta \lra \cM^{\balpha}_{(\cC, \bt), \blambda_0 /T}(r,d, n)$ which is a lift of $h : \Delta 
\lra T $ such that $\tilde{h}(0) = \x \in \cM^{\balpha}_{((\cC_0, \bt_0), \blambda_0)}(r,d, n)$. 
$$
\begin{array}{ccc}
    &  & \cM^{\balpha}_{(\cC, \bt), \blambda_0 /T}(r, n, d) \\
    &     &   \\
    & \tilde{h} \nearrow & \downarrow  {\footnotesize \varphi_{r,n,d,\blambda_0}} \\
    &          &   \\  
 \Delta& \stackrel{h}{\lra} & T  \times \{\blambda_0\} \\
 \end{array}.
$$

Next we will define a global foliation $\cIF$ on the total space of $\cM^{\balpha}_{(\cC, \bt)/T}(r,d, n)$ from  isomonodromic deformations of the $\balpha$-stable parabolic connections. We mean that a foliation $\cIF$ is a subsheaf of 
the tangent sheaf $\Theta_{\cM^{\balpha}_{(\cC, \bt)/T}(r,d, n)}$.  
We will show that the global foliation $\cIF$ coming from isomonodromic deformations has the Painlev\'e property, whose precise meaning will be defined in 
Theorem \ref{thm:Painleve}.

Let us consider 
the universal covering map   
$u: \tilde{T} \ra T =\cM'_{g,n}$. Note that $u$
factors  thorough the morphism  $u':\tilde{T} \ra T'$. 
Pulling back the  
fibration 
$\phi^{r}_n :\cRP^r_{n, T'} \lra T' \times \cA^{(n)}_r $
in 
(\ref{eq:phi-r-n}) by $u'$, we obtain the fibration   
$\cRP^r_{n, T'} \times_{T'}   \tilde{T}  \lra \tilde{T}$, 
which becomes a trivial fibration as explained in Section 4 in \cite{IIS-1}.   
This means that if we fix 
 a point $(C_0, \bt_0) \in T $ there exists 
an isomorphism 
\begin{equation}\label{eq:isom}
\pi:\cRP^r_{n, T'} \times_{T'} \tilde{T}  \stackrel{\simeq}{\lra} \cRP^r_{(C_0, \bt_0)} \times \tilde{T}
\end{equation}
which makes the following diagram commute. 
\begin{equation}\label{eq:comm1}
\begin{CD}
\cRP^r_{n, T'} \times_{T'} \tilde{T}  @> \pi >\simeq> \cRP^r_{(C_0, \bt_0)} \times \tilde{T} 
 \\
  @V \widetilde{\phi^{r}_n} VV@VV p_2 \times \phi^r_{(C_0, \bt_0)}V \\
 \tilde{T} \times \cA^{(n)}_r    @>>> \tilde{T} \times \cA^{(n)}_r. 
 \end{CD}
\end{equation}
Fixing $\ba \in \cA^{(n)}_r$, we set  
$
\cRP^r_{n, T', \ba} = \left(\phi^r_n \right)^{-1} ( T' \times \{\ba \} ).  
$
From the morphisms (\ref{eq:character}) and  (\ref{eq:phi-r-n}), 
we also have the following commutative diagram:
\begin{equation}\label{eq:comm2}
\begin{CD}
\cRP^r_{n, T', \ba} \times_{T'} \tilde{T}  @> \pi_{\ba} >\simeq> \cRP^r_{(C_0, \bt_0), \ba} \times \tilde{T}
\\
@V \widetilde{\phi^{r}_{n, \ba}} VV @VVp_2V \\
\tilde{T} \times \{\ba\}  @> \simeq>> \tilde{T}.
\end{CD}
\end{equation}

By using the isomorphism  (\ref{eq:comm2}) we can define the smooth part of 
$\cRP^r_{n, T', \ba} \times_{T'} \tilde{T}$ by 
$$
\left(\cRP^r_{n, T', \ba} \times_{T'} \tilde{T} \right)^{\sharp} = 
\pi_{\ba}^{-1}\left( \left( \cRP^r_{(C_0, \bt_0), \ba} \right)^{\sharp} \times \tilde{T} \right)
$$ 
where $\left( \cRP^r_{(C_0, \bt_0), \ba} \right)^{\sharp}$ is the smooth locus  of $\cRP^r_{(C_0, \bt_0), \ba}$ (cf. (\ref{eq:smooth})). 
Note that for generic $\ba$  the variety $\cRP^r_{(C_0, \bt_0), \ba}$ 
is non-singular,  but for special $\ba$,   $\cRP^r_{(C_0, \bt_0), \ba}$ 
does have singularities (cf. [(\ref{eq:rep-sing}), Remark \ref{ref:sing}]). 

We also have the following commutative diagram
\begin{equation}\label{eq:comm3}
\begin{CD}
\left(\cRP^r_{n, T', \ba} \times_{T'} \tilde{T} \right)^{\sharp} @> \pi_{\ba} >\simeq>\left( \cRP^r_{(C_0, \bt_0), \ba} \right)^{\sharp} \times \tilde{T} 
 \\
  @VVV @VVp_2V \\
 \tilde{T} \times \{\ba\}  @> \simeq>> \tilde{T} 
 \end{CD}.
\end{equation}
By using this isomorphism,  for any fixed $\ba \in \cA^{(n)}_r $,  
we  define  the set of constant sections
\begin{equation}
\Isomd ( \tilde{T}, \left(\cRP^r_{n, T', \ba} 
\times_{T'} \tilde{T} \right)^{\sharp}) = 
\left\{ \sigma:\tilde{T} \ra \left( \cRP^r_{n, T', \ba} \times_{T'} 
 \tilde{T} \right)^{\sharp}, \mbox{constant} \right\}.
\end{equation}
Note that by using the isomorphism (\ref{eq:comm3}), we have a natural 
isomorphism 
\begin{equation}\label{eq:const}
\Isomd ( \tilde{T}, \left(\cRP^r_{n, T', \ba} 
\times_{T'} \tilde{T} \right)^{\sharp}) \simeq  
\left( \cRP^r_{(C_0, \bt_0), \ba} \right)^{\sharp}.
\end{equation}
A section $\sigma \in \Isomd(\tilde{T}, \left(\cRP^r_{n, T', \ba} \times_{T'} \tilde{T}\right)^{\sharp}) $ is called an {\em isomonodromic section}  by trivial reason and 
its image $\sigma(\tilde{T})$ is called an {\em isomonodromic flow}.  

Next, considering the pullback of $\varphi_{r,n,d}$ in (\ref{eq:conn-fam}) 
by $\tilde{T} \lra T$, 
we can obtain the family of moduli spaces of  $\balpha$-stable parabolic connections
\begin{equation}\label{eq:pullback-fam}
\widetilde{\varphi_{r,n,d}}: \cM^{\balpha}_{(\cC, \bt)/\tilde{T}} \lra \tilde{T} \times \Lambda^{(n)}_r(d). 
\end{equation}
Fixing  $\blambda \in \Lambda$ such that $rh(\blambda) = \ba$, we also obtain the restricted family over $\tilde{T} \times \{ \blambda \}$ 
\begin{equation}\label{eq:pullback-famr}
\widetilde{\varphi_{r,n,d, \blambda}}: \cM^{\balpha}_{((\cC, \bt), \blambda)/\tilde{T}} \lra \tilde{T} \times \{ \blambda \}. 
\end{equation}
Restricting the Riemann-Hilbert correspondence (\ref{eq:grh}) to this space, 
we obtain the following commutative diagram
\begin{equation}\label{eq:grh2}
\begin{CD}
\cM^{\balpha}_{((\cC, \tilde{\bt}), \blambda) /\tilde{T}}(r,n,d) @> \RH_{ \blambda } >>
  \cRP^r_{n,T, \ba} \times_{T} \tilde{T}  \\
  @V \widetilde{\varphi_{r,n,d, \blambda}} VV @VV \widetilde{\phi^{r}_{n, \ba}} V \\
 \tilde{T} \times \{ \blambda \}
 @> Id \times rh >> \tilde{T} \times \{ \ba \}
 \end{CD}. 
\end{equation}
Note that by Theorem \ref{thm:rh} the morphism $\RH_{\blambda}$ gives an 
analytic resolution of singularities.   
Set  
\begin{equation}\label{eq:smooth-conn}
 \left( \cM^{\balpha}_{((\cC, \tilde{\bt}), \blambda) /\tilde{T}}(r,n,d) ) \right)^{\sharp} = \RH_{\blambda}^{-1} ( (\cRP^r_{n,T, \ba} \times_{T} \tilde{T})^{\sharp} ),
\end{equation}
and
\begin{equation}\label{eq:sing-conn}
 \left( \cM^{\balpha}_{((\cC, \tilde{\bt}), \blambda) /\tilde{T}}(r,n,d) ) \right)^{sing} = \RH_{\blambda}^{-1} ( (\cRP^r_{n,T, \ba} \times_{T} \tilde{T})^{sing} ). 
\end{equation} 
(Cf. (\ref{eq:rep-sing}), (\ref{eq:smooth})). 
Then we have an analytic isomorphism 
$$
(\RH_{\blambda})^{\sharp} : \left( \cM^{\balpha}_{((\cC, \tilde{\bt}), \blambda) /\tilde{T}}(r,n,d) ) \right)^{\sharp} \stackrel{\simeq}{\lra} (\cRP^r_{n,T, \ba} \times_{T} \tilde{T})^{\sharp}. 
$$
Now we define: 
\begin{equation}\label{eq:isomod}
\Isomd(\tilde{T}, 
\left( \cM^{\balpha}_{((\cC, \tilde{\bt}), \blambda) /\tilde{T}}(r,n,d) ) \right)^{\sharp} ) 
= \RH_{\blambda}^{-1} ( \Isomd(\tilde{T}, 
(\cRP^r_{n,T, \ba} \times_{T} \tilde{T})^{\sharp}) ) .  
\end{equation}
Each section $\sigma \in \Isomd( \tilde{T}, \left(\cM^{\balpha}_{((\cC, \tilde{\bt}), \blambda) /\tilde{T}}(r,n,d) \right)^{\sharp} ) $ is called an 
{\em isomonodromic section} on $\left(\cM^{\balpha}_{((\cC, \tilde{\bt}), \blambda) /\tilde{T}}(r,n,d) \right)^{\sharp} $ and its image 
$$
\sigma(\tilde{T}) \subset \left(\cM^{\balpha}_{((\cC, \tilde{\bt}), \blambda) /\tilde{T}}(r,n,d) \right)^{\sharp}
$$ 
is called an {\em isomonodromic flow}. Note that 
since the Riemann-Hilbert correspondence $(\RH_{\blambda})^{\sharp}$ 
is a highly  non-trivial analytic isomorphism,  isomonodromic 
flows $\{ \sigma(\tilde{T}) \}$ are not  constant any more and it is known that 
they define highly transcendental analytic functions.  

From the morphism (\ref{eq:pullback-famr}) restricted to $\left(\cM^{\balpha}_{((\cC, \tilde{\bt}), \blambda) /\tilde{T}}(r,n,d) \right)^{\sharp}$, 
we obtain the natural sheaf homomorphism 
$$
\Theta_{\left(\cM^{\balpha}_{((\cC, \tilde{\bt}), \blambda) /\tilde{T}}(r,n,d) \right)^{\sharp}} \stackrel{\widetilde{\varphi_{r,n,d, \blambda}}^{*}}{\lra} \widetilde{\varphi_{r,n,d, \blambda}}^{*} ( \Theta_{\tilde{T}})_{|\left(\cM^{\balpha}_{((\cC, \tilde{\bt}), \blambda) /\tilde{T}}(r,n,d) \right)^{\sharp}} \lra 0.$$
Then the set of all isomonodromic sections defines a sheaf homomorphism 
\begin{equation}\label{eq:tangent-seq}
\cV_{\blambda}: \widetilde{\varphi_{r,n,d, \blambda}}^{*} (\Theta_{\tilde{T}})_{|\left(\cM^{\balpha}_{((\cC, \tilde{\bt}), \blambda) /\tilde{T}}(r,n,d) \right)^{\sharp}}
\lra \Theta_{\left(\cM^{\balpha}_{((\cC, \tilde{\bt}), \blambda) /\tilde{T}}(r,n,d) \right)^{\sharp}} 
\end{equation}
which gives  a splitting of  the homomorphism $\widetilde{\varphi_{r,n,d, \blambda}}^{*}$. The splitting (\ref{eq:tangent-seq}) is algebraic, because the condition of isomonodromic flows given by the vanishing of the curvature $2$-forms  
of the associated universal connections.  Since the 
exceptional locus for $\RH = \cup_{\blambda} \RH_{\blambda}$ has 
codimension at least $2$, by Hartogs' theorem,  
it is easy to see that this algebraic splitting (\ref{eq:tangent-seq}) can be extend to the whole family of moduli spaces, and we obtain an extended homomorphism  
\begin{equation}\label{eq:splitting}
{\mathcal V}_{\blambda}: \widetilde{\varphi_{r,n,d, \blambda}}^{*} (\Theta_{\tilde{T}})
\lra \Theta_{\cM^{\balpha}_{((\cC, \tilde{\bt}), \blambda) /\tilde{T}}(r,n,d)}.
\end{equation}

Under the notation above, we have the following 
\begin{Definition}{\rm \label{def:isom-fol}
\begin{enumerate} 
\item
The foliation  $ {\mathcal I \mathcal F}_{\blambda} $ defined by the subsheaf
\begin{equation}\label{eq:foliation}
\cIF_{\blambda} = {\mathcal V}_{\blambda}(\widetilde{\varphi_{r,n,d, \blambda}}^{*} (\Theta_{\tilde{T}}))
\subset
\Theta_{\cM^{\balpha}_{((\cC, \tilde{\bt}), \blambda) /\tilde{T}}(r,n,d)} 
\end{equation}
is called {\em an isomonodromic foliation on $\cM^{\balpha}_{((\cC, \tilde{\bt}), \blambda) /\tilde{T}}(r,n,d)$}. 
\item Let $h:\Delta \lra \tilde{T}$ be a holomorphic embedding such that 
$h(t)=(C_t, \bt_t)$ for $t \in \Delta$.  A holomorphic map  
$\tilde{h}:\Delta \lra \cM^{\balpha}_{((\cC, \tilde{\bt}), \blambda) /\tilde{T}}(r,n,d)$ such that 
$ \widetilde{\varphi_{r,n,d, \blambda}} \circ \tilde{h} = h $ is called a $\cIF_{\blambda}$-lift of $h$ if  $\tilde{h}$ is tangent to $\cIF_{\blambda}$, 
that is, 
$\tilde{h}_{*} (\Theta_{\Delta}) \subset \cIF_{\blambda}$. 
\end{enumerate}}
\end{Definition}

\begin{Lemma}\label{lem:sep-fol} 
Let $h:\Delta \lra \tilde{T}$ be a holomorphic embedding and 
$\tilde{h}:\Delta \lra \cM^{\balpha}_{((\cC, \tilde{\bt}), \blambda) /\tilde{T}}(r,n,d)$ a $\cIF_{\blambda}$-lift of $h$. Then the image of $\RH_{\blambda} \circ \tilde{h}$ lies  in the image of a constant section 
$\sigma \in \Isomd ( \tilde{T}, \left(\cRP^r_{n, T', \ba} \times_{T'} \tilde{T}\right))$. 
\end{Lemma}
\begin{proof}
Note that a lift $\tilde{h}$ of $h$ 
corresponds to a $1$-parameter deformation of 
$\balpha$-stable parabolic connection under a deformation of $n$-pointed curves associated to $h:\Delta \lra \tilde{T}$.  
Since 
$\left( \cM^{\balpha}_{((\cC, \tilde{\bt}), \blambda) /\tilde{T}}(r,n,d) \right)^{\sharp}$ is a Zariski dense open subset of $\left( \cM^{\balpha}_{((\cC, \tilde{\bt}), \blambda) /\tilde{T}}(r,n,d) \right)$, we see that the curvature form 
vanishes on the $\cIF$-foliation defined on the total space $\left( \cM^{\balpha}_{((\cC, \tilde{\bt}), \blambda) /\tilde{T}}(r,n,d) \right)$.   Therefore 
if $\tilde{h}$ is a $\cIF$-lift of $h$, 
we can conclude that the deformation of  connections is isomonodromic. 
Hence the associated representations  of the fundamental group of $\cC_{t} \setminus \{ \bt_t \} $ are constant, which means that $\RH_{\blambda}(\tilde{h}(\Delta))$ is contained in the image of a constant section of $\left(\cRP^r_{n, T', \ba} \times_{T'} \tilde{T}\right) \lra \tilde{T} $. 
\end{proof}

Now, we can show that the isomonodromic foliation is a 
differential system satisfying the Painlev\'e property (cf. \cite{Mal}, \cite{Miwa} and \cite{IIS-3}).

\begin{Theorem} \label{thm:Painleve}
For any $\blambda \in \Lambda^{(n)}_r(d)$, 
the isomonodromic foliation $\cIF_{\blambda}$ defined on $ \cM^{\balpha}_{((\cC, \tilde{\bt}), \blambda) /\tilde{T}}(r,n,d) $ has Painlev\'e property.  
That is, for any  holomorphic embedding  $h:\Delta  \lra \tilde{T}$ of 
the unit disc $\Delta = \{ t \in \C | |t| < 1 \}$ such that 
$h(0) = (C, \bt)$ and $ \x =( E, \nabla, \{ l^{(i)}_{\ast} \}_{1 \leq i \leq n}) \in \cM^{\balpha}_{((C, \bt), \blambda)}(r,n,d)$, there exists the unique $\cIF_{\blambda}$-lift 
$$
\tilde{h}:\Delta \lra \cM^{\balpha}_{((\cC, \tilde{\bt}), \blambda) /\tilde{T}}(r,n,d)
$$ of $h$ such that $\tilde{h}(0) = \x$.  
\end{Theorem}
\begin{proof}
If $\x \in \left(\cM^{\balpha}_{((C, \bt), \blambda)}(r,n,d)\right)^{\sharp}$, 
there is a unique isomonodromic section  $ \sigma: \tilde{T} \lra  \left(\cM^{\balpha}_{((\cC, \tilde{\bt}), \blambda) /\tilde{T}}(r,n,d)\right)^{\sharp}$ 
such that $\sigma((C, \bt)) = \x $.  The holomorphic map  
$\tilde{h} = \sigma \circ h : \Delta \lra  \left(\cM^{\balpha}_{((\cC, \tilde{\bt}), \blambda) /\tilde{T}}(r,n,d) \right)^{\sharp}$ is the unique 
$\cIF_{\blambda}$-lift of $h$.

Let us consider the case when $\x \in  \left(\cM^{\balpha}_{((C, \bt), \blambda)}(r,n,d) \right)^{sing}$. 
Pulling back the commutative diagrams (\ref{eq:grh2}) and (\ref{eq:comm2}) via 
the embedding $h:\Delta \lra \tilde{T}$, we obtain the 
commutative diagram 
\begin{equation}\label{eq:rest-cd}
\begin{CD}
\cM^{\balpha}_{((\cC, \tilde{\bt}), \blambda) /\Delta}(r,n,d) @> \pi_{\ba} \circ \RH_{ \blambda } >>
  \cRP^r_{(C_0, \bt_0), \ba}  \times \Delta  \\
  @V \widetilde{\varphi_{\Delta}} VV @VV p_2 V \\
 \Delta 
 @> Id  >> \Delta
 \end{CD}. 
\end{equation}

The restriction of the foliation $\cIF_{\blambda}$ to 
$\cM^{\balpha}_{((\cC, \tilde{\bt}), \blambda) /\Delta}(r,n,d) $
determines a vector field $v_{\blambda}$ on $\cM^{\balpha}_{((\cC, \tilde{\bt}), \blambda) /\Delta}(r,n,d) $ such that $ \widetilde{\varphi_{\Delta}}_{\ast}(v_{\blambda}) = \frac{\partial}{\partial t}$ where $t$ is a coordinate of 
$\Delta$. We will show that there exist a unique section 
$\tilde{h}:\Delta \lra \cM^{\balpha}_{((\cC, \tilde{\bt}), \blambda) /\Delta}(r,n,d) $ such that $\tilde{h}(0) = \x $ and $\tilde{h}_{*}(\frac{\partial}{\partial t})= v_{\blambda}$, which gives a $\cIF_{\blambda}$-lift of $h$. 
Such a section $\tilde{h}$  can be locally given by an analytic solution of the 
Cauchy problem of an ordinary 
differential equation associated to the vector field $v_{\blambda}$. 
Such an analytic solution can be locally given by holomorphic functions of 
$t$ on  $\Delta_{\epsilon} = \{ t \in \C \ | \ |t|< \epsilon \}$ for some 
$0 < \epsilon < 1$. This gives a section $ \tilde{h}_{\epsilon}:
\Delta_{\epsilon} \lra 
\cM^{\balpha}_{((\cC, \tilde{\bt}), \blambda) /\Delta_{\epsilon}}(r,n,d) $
which is a $\cIF_{\blambda}$-lift of $h_{\epsilon} = h_{|\Delta_{\epsilon}}$.
Let $\epsilon_1$ be the supremum of  $\epsilon$ such that 
a $\cIF_{\lambda}$ lift of $h_{\epsilon}$ exists.  The above argument shows 
that $\epsilon_1 > 0$. Now we will show that $\epsilon_1 = 1$. 
Assume the contrary, that is, $\epsilon_1 < 1$, and 
let $\tilde{h}_{\epsilon_1}: \Delta_{\epsilon_1} \lra 
\cM^{\balpha}_{((\cC, \tilde{\bt}), \blambda) /\Delta_{\epsilon_1}}(r,n,d)$ 
be the section over $\Delta_{\epsilon_1}$. 

Let $p_1: \cRP^r_{(C_0, \bt_0), \ba} \times 
\Delta \lra \cRP^r_{(C_0, \bt_0), \ba}$ be the first projection and  
consider the morphism 
$$
p_1 \circ \pi_a \circ  
\RH_{\blambda} :\cM^{\balpha}_{((\cC, \tilde{\bt}), \blambda) /\Delta}(r,n,d)\lra \cRP^r_{(C_0, \bt_0), \ba}.
$$
By definition of $\left(\cM^{\balpha}_{((C, \bt), \blambda)}(r,n,d) \right)^{sing}$, the point  $\y = p_1 \circ \pi_a \circ  
\RH_{\blambda}(\x) $ is a singular point of 
$\cRP^r_{(C_0, \bt_0), \ba}$ and let 
$$
\cK_{\Delta, \y} = ( \pi_a \circ  
\RH_{\blambda})^{-1}( \{\y\} \times \Delta) \subset \left(\cM^{\balpha}_{((\cC, 
\tilde{\bt}), \blambda)/\Delta}(r,n,d) \right)^{sing} 
$$ 
denote the exceptional locus 
dominated over $\{ \y \} \times \Delta$. Then restricting 
(\ref{eq:rest-cd}) to $\cK_{\Delta,\y}$,  we have the following
commutative diagram:
\begin{equation}\label{eq:rest-cd2}
\begin{CD}
\cK_{\Delta,\y}  @> \pi_{\ba} \circ \RH_{ \blambda } >> \{ \y \} \times \Delta  \\
  @V \widetilde{\varphi_{\Delta, \y}} VV @VV p_2 V \\
 \Delta @> Id  >> \Delta. 
 \end{CD}   
\end{equation} 
From Theorem \ref{thm:rh}, we see that 
$ \pi_{\ba} \circ \RH_{\blambda}$ is a resolution of singularity of 
 $\cRP^r_{(C_0, \bt_0), \ba} \times \Delta$, hence  
each fiber of $\widetilde{\varphi_{\Delta, \y}}:\cK_{\Delta,\y} \lra \Delta$ is compact. Now from Lemma \ref{lem:sep-fol},  we see that  
$\tilde{h}_{\epsilon_1}(\Delta_{\epsilon}) \subset \cK_{\Delta_{\epsilon_1},\y}$. 
Moreover since $ \widetilde{\varphi_{\Delta, \y}}$ is proper, 
we see that 
$
\tilde{h}_{\epsilon_1}(\overline{\Delta_{\epsilon_{1}}}) \subset \cK_{\overline{\Delta_{\epsilon_1}},\y}
$
where $\overline{\Delta_{\epsilon_1}} =\{t,  |t| \leq \epsilon_1 \} $.
Take and fix $t=b$ such that $|b| = \epsilon_1$. Then 
$$
\tilde{h}_{\epsilon_1}(b) = \y_b \in \cK_{\overline{\Delta_{\epsilon_1}},\y}
\subset \cM^{\balpha}_{((\cC, 
\tilde{\bt}), \blambda)/\overline{\Delta_{\epsilon_1}}}(r,n,d)
$$
Starting from $t=b$ and $\y_b$, we can extend the section $\widetilde{h}_{\epsilon_1}$ over 
$\Delta(b, \epsilon_b) = 
\{ t \in \Delta \  |  \ |t - b| < \epsilon_b \}$ 
with $ 0 < \epsilon_b \leq 1 - \epsilon_1 $. 
Again, from the compactness of the fiber of $\widetilde{\varphi_{\Delta, \y}}:\cK_{\Delta,\y} \lra \Delta$, 
we can show  that the minimum $\epsilon_0$ 
 of $\epsilon_b $ for $|b| = \epsilon_1$  is positive, hence 
 for $\epsilon = \epsilon_1 + \epsilon_0$ the section $\tilde{h}_{\epsilon}$
 exists and this contradicts to the fact that 
 $\epsilon_1 $ is the supremum and  
  $\epsilon_1 < \epsilon$. 
\end{proof}

\begin{Remark}{\rm  
Let us remark that the isomonodromic foliation 
$\cIF_{\blambda}$ on 
$\cM^{\balpha}_{((\cC, \tilde{\bt}), \blambda) /\tilde{T}}(r,n,d)$
descends to a foliation on 
$\cM^{\balpha}_{((\cC, \tilde{\bt}), \blambda) /T'}(r,n,d)$
 under the covering map $\tilde{T} \lra T'$, which 
we also denote  by $\cIF_{\blambda}$. 
Recall that 
the  isomonodromic section (\ref{eq:const})
is the constant  section with respect to  
the isomorphism (\ref{eq:isom}).  Moreover, when the 
base point $ \ast \in T'$  corresponds to $(C_0, \bt_0)$,  
the fundamental group $\pi_{1}(T', *)$ acts on the moduli space 
$ \cRP^r_{(C_0, \bt_0)} $ via the action to the generators of 
 $ \pi_{1}(C_0 \setminus D(\bt_0), \ast') $. Therefore, we can 
define the local isomonodromic sections for $ \cRP^r_{n, T', \ba'} 
\lra T'$, which also defines a local isomonodromic sections for 
$
\left(\cM^{\balpha}_{((\cC, \bt), \blambda)/T'}\right)^{\sharp}  
\lra T'
$. 
Now the set of local isomonodromic sections determines 
a splitting homomorphism 
$ \cV_{\blambda} $ 
like (\ref{eq:splitting}), and it defines an 
isomonodromic foliation 
$$
\cIF_{\blambda} = \cV_{\blambda}(\Theta_{T'}) \subset 
\Theta_{\cM^{\balpha}_{((\cC, \bt), \blambda)/T'}}
$$
which is obviously the descent of the original isomonodromic foliation  
on $\cM^{\balpha}_{((\cC, \bt), \blambda)/\tilde{T'}}$ 
}
\end{Remark}

\end{document}